\expandafter\def\csname ver@fixltx2e.sty\endcsname{}
\documentclass[AMA,STIX1COL]{WileyNJD-v2}

\usepackage[export]{adjustbox}
\usepackage{makecell}

\usepackage{algorithm,algcompatible,amsmath}
\algnewcommand\INPUT{\item[\textbf{Input:}]}%
\algnewcommand\OUTPUT{\item[\textbf{Output:}]}%

\def\Rbb{\mathbb{R}}

\def\Abold{\mathbf{A}}
\def\Bbold{\mathbf{B}}
\def\Cbold{\mathbf{C}}
\def\Dbold{\mathbf{D}}
\def\Ebold{\mathbf{E}}
\def\Fbold{\mathbf{F}}
\def\Gbold{\mathbf{G}}
\def\Hbold{\mathbf{H}}
\def\Ibold{\mathbf{I}}

\def\Kbold{\mathbf{K}}

\def\Mbold{\mathbf{M}}

\def\Pbold{\mathbf{P}}
\def\Qbold{\mathbf{Q}}
\def\Rbold{\mathbf{R}}
\def\Sbold{\mathbf{S}}
\def\Tbold{\mathbf{T}}
\def\Ubold{\mathbf{U}}
\def\Vbold{\mathbf{V}}
\def\Wbold{\mathbf{W}}
\def\Xbold{\mathbf{X}}
\def\Ybold{\mathbf{Y}}
\def\Zbold{\mathbf{Z}}

\def\cbold{\mathbf{c}}

\def\fbold{\mathbf{f}}

\def\qbold{\mathbf{q}}

\def\ubold{\mathbf{u}}

\def\wbold{\mathbf{w}}
\def\xbold{\mathbf{x}}

\def\zbold{\mathbf{z}}

\def\mubold{\boldsymbol{\mu}}

\def\Sigmabold{\mathbf{\Sigma}}

\def\Gammabold{\mathbf{\Gamma}}

\def\zetabold{\boldsymbol{\zeta}}
\def\zerobold{{\bf 0}}
\def\0{\mathbf{\0}}

\DeclareMathAlphabet\mathbfcal{OMS}{cmsy}{b}{n}

\usepackage{amsmath}

\articletype{Article Type}

\received{26 April 2016}
\revised{6 June 2016}
\accepted{6 June 2016}

\begin{document}

\title{Model Reduction Framework with a New Take on Active Subspaces for Optimization Problems with Linearized Fluid-Structure Interaction Constraints}

\author[1]{Gabriele Boncoraglio*}
\author[1,2,3]{Charbel Farhat}
\author[4]{Charbel Bou-Mosleh}

\authormark{Boncoraglio \textsc{et al}}


\address[1]{\orgdiv{Department of Aeronautics and Astronautics}, \orgname{Stanford University},
\orgaddress{Durand Building, 496 Lomita Mall, Stanford, \state{CA} 94305-4035, \country{USA}}}
\address[2]{\orgdiv{Department of Mechanical Engineering}, \orgname{Stanford University},
\orgaddress{Building 530, 440 Escondido Mall, Stanford, \state{CA} 94305-3030, \country{USA}}}
\address[3]{\orgdiv{Institute for Computational and Mathematical Engineering}, \orgname{Stanford University},
\orgaddress{Huang Engineering Center, 475 Via Ortega, Suite 060, Stanford, \state{CA} 94305-4042, \country{USA}}}
\address[4]{\orgdiv{Department of Mechanical Engineering}, \orgname{Notre Dame University-Louaize}, \orgaddress{Zouk Mosbeh, Lebanon}}

\corres{*Gabriele Boncoraglio , Department of Aeronautics and Astronautics, Stanford University, Durand Building, 496 Lomita Mall, Stanford, CA, 94305-4035.\\ \email{gbonco@stanford.edu}}

\abstract[Summary]{In this paper, a new take on the concept of an active subspace for reducing the dimension of the design parameter space in a multidisciplinary analysis and optimization (MDAO) problem
is proposed. The new approach is intertwined with the concepts of adaptive parameter sampling, projection-based model order reduction, and a database of linear, projection-based reduced-order models
equipped with interpolation on matrix manifolds, in order to construct an efficient computational framework for MDAO. The framework is fully developed for MDAO problems with linearized fluid-structure
interaction constraints. It is applied to the aeroelastic tailoring, under flutter constraints, of two different flight systems: a flexible configuration of NASA's Common Research Model; and  NASA's
Aeroelastic Research Wing \#2 (ARW-2). The obtained results illustrate the feasibility of the computational framework for realistic MDAO problems and highlight the benefits of the new approach for
constructing an active subspace in both terms of solution optimality and wall-clock time reduction.}


\maketitle


\section{Introduction}

Multidisciplinary analysis and optimization (MDAO) problems with fluid-structure interaction (FSI) constraints arise in many engineering applications including design and optimal control. Solving such
problems requires repeated evaluations of the FSI constraints when using, for example, the nested analysis and design (NAND) approach. Typically, each evaluation of FSI constraints incurs the solution
of some coupled computational fluid dynamics (CFD) - computational structural dynamics (CSD) problems. In an industrial setting, these CFD-CSD problems are usually large-scale when high-fidelity
simulations are required. Furthermore, the number of solutions of coupled CFD-CSD problems increases when the design parameter space $\mathcal D$ is high-dimensional. This is because the number of
evaluations of the objective function and constraints incurred by an MDAO problem generally increases with the dimension $N_{\mathcal D}$ of $\mathcal D$. For all these reasons, solving highly
parameterized, FSI-constrained, MDAO problems can be very computationally intensive. Consequently, the main objective of this paper is to present an innovative, comprehensive, computational framework
for efficiently solving such problems when the FSI constraints are linearized in the spatial, temporal, or frequency domain. Such linearizations are common when the FSI constraints pertain to stability
or control problems.

The computational cost of a linearized, coupled FSI problem can be reduced by substituting the underlying coupled high-dimensional model (HDM) with a less computationally intensive surrogate model such
as a linear, coupled, projection-based, reduced-order model (PROM). Such a linear PROM can be {\it local} (or pointwise), or {\it global}. In the context of this paper, the label local describes a PROM
that was trained at a single parameter point $\mubold \in {\mathcal D} \subset {\mathbb R}^{N_{\mathcal D}}$ of the design parameter space $\mathcal D$. On the other hand, the label {\it global} refers
to a PROM that was trained at multiple points in $\mathcal D$ so that it is robust and accurate within a large region of this parameter space.

The overhead cost -- also known as the offline cost -- associated with the construction of a local PROM is lower than that associated with the construction of a global PROM, given that it involves
training at a single parameter point. Unfortunately, a local PROM does not usually perform well away from the training point. Therefore, it is not suitable for the solution of an MDAO problem unless it
is reconstructed everytime a new parameter point is queried, which is inefficient. Alternatively, a global PROM may be suitable for MDAO problems, from the robustness and accuracy viewpoints. However,
such a PROM is bound to be of a larger dimension than needed at any parameter point visited online by the optimization trajectory, given that it was trained offline at multiple points of the parameter
space; this also hinders computational efficiency.

In order to overcome the inefficiency issues highlighted above, different approaches have been proposed in the literature. In the context of linear PROMs, two of them are noteworthy.
First, the adaptive approach consisting in progressively updating an initial PROM by continuously training it at a subset of the parameter points previously visited by the optimization trajectory~\cite{zahr2015progressive}. In this approach, which breaks the traditional offline-online framework of model order reduction, the reduced optimization problems are equipped with a nonlinear trust-region based on
a residual error indicator to keep the optimization trajectory in a region of $\mathcal D$ where the PROM is accurate. Hence, whereas this approach is equally applicable to linear PROMs, it was
designed primarily for nonlinear PROMs. The second approach consists in carefully sampling the parameter space $\mathcal D$ at a small set of points, constructing at each sampled parameter point
a local PROM, and storing in a database all constructed local PROMs as well as some related quantities such as the underlying reduced-order bases (ROBs)~\cite{amasallem2010towards}. This approach is natural for linear
PROMs as once the database has been constructed offline, a linear PROM -- as well as its sensitivities -- can be constructed in real-time at any unsampled parameter point queried by the optimization
trajectory using interpolation on matrix manifolds~\cite{amsallem2011online}. This approach has been recently demonstrated for MDAO problems characterized by a small-dimensional
design parameter space -- say $N_{\mathcal D} \le 6$ -- and shown to deliver excellent speedup factors. For high-dimensional design parameter spaces however -- say $N_{\mathcal D} \ge 10$,
building a database of local PROMs suffers from the curse of dimensionality and rapidly becomes impractical if not infeasible. For example, sampling only 2 parameter points along each
dimension of a design parameter space of dimension 10 leads to the construction of $2^{10} = 1,024$ PROMs! Adaptive parameter sampling procedures -- also known as greedy sampling (or simply greedy)
procedures -- and particularly those equipped with a saturation technique~\cite{hesthaven2014efficient} attenuate somehow this issue but insufficiently. Alternatively, this paper addresses this issue by focusing on a
PROM computational framework for MDAO where the concept adaptive parameter sampling, that of a database of linear PROMs, and a new take on the concept of active subspaces (AS) are intertwined.

The remainder of this paper is organized as follows. Section~\ref{sec:MDAO} formulates the problem of interest in sufficient details to keep this paper as self-contained as possible.
Section~\ref{sec:PMOR} reviews and intertwines the concepts of AS, projection-based model order reduction (PMOR), and a database of linear FSI PROMs equipped with interpolation on matrix manifolds,
in the context of MDAO problems with linearized FSI constraints. It also proposes an effective parameter sampling method for constructing an AS. Section~\ref{sec:APP} illustrates the proposed PROM-based computational
framework for MDAO with the aeroelastic tailoring, under flutter constraints, of two different systems: a flexible configuration of NASA's Common Research Model (CRM); and NASA's Aeroelastic Research Wing ARW-2.

\section{MDAO Problem with a Linearized FSI Constraint}
\label{sec:MDAO}

\subsection{Problem Formulation}

The main focus of this paper is on the acceleration of the solution of MDAO problems under a high-dimensional, linearized, FSI constraint, as well as other linear or nonlinear constraints, including
partial differential equation (PDE)-based constraints. Such problems can be written as
    \begin{equation}\label{eq:generalPb}
      \boxed{
        \begin{aligned}
		& \underset{\textbf{q}, \, \zbold, \,\mubold}{\text{min}}
		& & f(\textbf{q}, \zbold, \mubold) \\
		& \text{s.t.} & &  \cbold(\textbf{q}, \zbold, \mubold)) & \leq & \, \boldsymbol{0} \\
		& & &  \Rbold_{L}(\textbf{q},\mubold) & = & \, \boldsymbol{0} \\
		& & &  \Rbold_{NL}(\zbold, \mubold) & = & \, \boldsymbol{0} \\
	\end{aligned}
      }
    \end{equation}
where:
\begin{itemize}
	\item $f(\cdot)$ represents a scalar objective function to minimize.
	\item $\textbf{q} \in \Rbb^{N_{q}}$ denotes a first high-dimensional vector of $N_{q}$ semi-discrete or discrete fluid and structural state variables involved in a high-dimensional, parametric,
		linearized, FSI constraint.
	\item $\zbold \in \Rbb^{N_{z}}$ denotes a second high-dimensional vector of $N_{z}$ semi-discrete or discrete state variables involved in a high-dimensional, parametric, PDE-based constraint.
	\item $\mubold  \in \mathcal D \subset \Rbb^{N_{\mathcal D}}$ denotes a vector of $N_{\mathcal D}$ design optimization parameters and is assumed to be relatively
		large -- say, $N_{\mathcal D} \ge 10$. Each component of this vector may represent some structural material property or an aerodynamic shape parameter.
	\item $\cbold(\cdot)$ represents a general set of linear and nonlinear algebraic constraints.
	\item $\Rbold_{L}(\textbf{q}, \mubold)$ is a high-dimensional, semi-discretized or discretized, parametric PDE that is linear in $\textbf{q}$, and the subscript $L$ designates the
		linear aspect.
	\item $\Rbold_{NL}(\zbold, \mubold)$ is a high-dimensional, semi-discretized or discretized, parametric PDE that is nonlinear in $\zbold$, and the subscript $NL$ designates the
		nonlinear aspect.
\end{itemize}

Specifically, the focus is set on the gradient-based solution of the MDAO problem~\eqref{eq:generalPb} using the NAND approach.

\subsection{Nested Analysis and Design Approach}

Two methodologies for solving MDAO problems such as~\eqref{eq:generalPb} can be considered: Simultaneous Analysis and Design (SAND) and Nested Analysis and Design (NAND). In the SAND approach,
both vectors of discrete state variables $\qbold$ and $\zbold$ as well as the vector of design parameters $\mubold$ are considered optimization parameters for problem~\eqref{eq:generalPb}.
In this case, the MDAO problem has potentially millions of optimization parameters given the high-dimensionality of $\Rbold_{L}(\textbf{q}, \mubold)$ and $\Rbold_{NL}(\zbold, \mubold)$. In the NAND
approach however, the only optimization parameters are the components of $\mubold$: in this case, $\qbold$ and $\zbold$ are recognized as implicit functions of $\mubold$ due to the constraints
$\Rbold_{L}(\textbf{q},\mubold) = \boldsymbol{0}$ and $\Rbold_{NL}(\zbold, \mubold) = \boldsymbol{0}$. It follows that when storage is an issue, the NAND approach is preferrable. Furthermore, the NAND
approach has the additional benefit or enabling the reuse of available solvers for the PDE-based constraints. For these reasons, it is here the method of choice.

Hence, problem~\eqref{eq:generalPb} is rewritten as
\begin{equation}\label{eq:nandPb}
        \begin{aligned}
		& \underset{ \textbf{q}, \, \zbold, \,\mubold}{\text{min}}
		& & f(\textbf{q}, \zbold, \mubold) \\
		& \text{s.t.} & &  \cbold(\textbf{q}, \zbold, \mubold) & \leq & \, \boldsymbol{0} \\
		& & &  \Rbold_{L}(\textbf{q}, \mubold)& = & \, \boldsymbol{0} \\
		& & &  \Rbold_{NL}(\zbold, \mubold)& = & \, \boldsymbol{0} \\
 \end{aligned}
 \rightarrow
 \boxed{
         \begin{aligned}
		 & \underset{\mubold \in {\mathcal D}}{\text{min}}
		& & f(\textbf{q}(\mubold), \zbold(\mubold), \mubold) \\
		& \text{s.t.} & &  \cbold(\textbf{q}(\mubold), \zbold(\mubold), \mubold) \leq \boldsymbol{0} \\
 \end{aligned}
 }
    \end{equation}
Given at iteration $i$ of the NAND approach a queried parameter point $\mubold_i \in {\mathcal D}$, $\Rbold_{L}(\textbf{q},\mubold) = \boldsymbol{0}$ and $\Rbold_{NL}(\zbold, \mubold) = \boldsymbol{0}$ are solved to
obtain $\qbold(\mubold_i)$ and $\zbold(\mubold_i)$, and then the objection function $f$ and the set of algebraic constraints represented by $\cbold$ are evaluated. This process can be written as

 \begin{equation*}\label{eq:nand2}
 \mubold_i \rightarrow
\left\{
\begin{aligned}
	& \Rbold_{L}(\qbold(\mubold_i), \mubold_i) & = & \, \boldsymbol{0} \\
	& \Rbold_{NL}(\zbold(\mubold_i), \mubold_i) & = & \, \boldsymbol{0} \\
 \end{aligned}
 \right.
 \rightarrow
 \left\{
\begin{aligned}
	& \qbold = \qbold(\mubold_i) \\
	& \zbold = \zbold(\mubold_i) \\
 \end{aligned}
 \right.
\end{equation*}

Unfortunately, the number of solutions of the PDE-based constraints generally increases when the dimension of the design parameter space $\mathcal D$ is increased, which is a disadvantage
of the NAND approach and a further motivation for model reducion.

\subsection{Three-Field Linearized Computational Framework for FSI}

The three-field framework pioneered in~\cite{farhat1995mixed} is adopted for modeling a coupled FSI problem such as that represented here by $\Rbold_{L}(\textbf{q},\mubold) = \boldsymbol{0}$.
In this framework, no limiting assumption is made about the behavior of the fluid subsytem whose governing equations of equilibrium are usually formulated in the
Arbitrary Lagrangian Eulerian (ALE) setting and discretized by a CFD approach. Similarly, no limiting assumption is made about the
behavior of the structural subsystem whose governing equations of equilibrium are formulated in the Lagrangian setting and
discretized by a finite element (FE) method. The CFD mesh of the fluid subsystem is viewed as a third subsystem that is assimilated
with a fictitious discrete structural system~\cite{degand2002three}, or the FE discretization of a
fictitious deformable continuous body~\cite{jasak2006automatic} -- so that it may move and/or deform as needed. In either case, it is referred to here as the pseudo-structural subsystem.

Most importantly, the three-field (fluid, fluid mesh, and structure) framework couples the fluid and structural subsystems by
simply enforcing the kinematic transmission conditions expressing the slip or no-slip wall boundary conditions of the
fluid subsystem, and the force transmission conditions expressing equilibrium at the fluid/structure interface.
It is overviewed below with emphasis on the linearization of its semi-discrete form.

\subsubsection{Three-field Formulation of FSI Problems}

The three-field framework for computational FSI outlined above can be summarized by the following three equations

\begin{equation}\label{eq:ALE}
\left\{
	\begin{array} {l c l}
		\displaystyle{\frac{\partial ( \mathcal{J W} )  } {\partial t}} +\mathcal{J } \nabla  \cdot \left(  \mathcal{F(W)} - \frac{\partial x}{\partial t} \mathcal{ W} \right) -\mathcal{J } \nabla  \cdot \mathcal{R(W)} &=& 0 \\[7 pt]
		\displaystyle{\rho \frac{\partial^2 \mathcal{U} }{\partial^2 t}} - div \left( \sigma \left( \epsilon( \mathcal{U} ) \right) \right) - \mathcal{B} & = & 0 \\[7 pt]
		\tilde{\rho} \displaystyle{\frac{\partial^2 x }{\partial^2 t}} - div \left( \widetilde{E} :  \tilde{\epsilon} ( x ) \right)  & = & 0
 \end{array}
 \right.
\end{equation}
where: the first equation is the Navier-Stokes equation written in ALE conservative form -- and can be augmented with
turbulence modeling; the second equation expresses the dynamic equilibrium of the structural subsystem; the third
equation models the fluid mesh as a fictitious, dynamic, deformable body; and all dependences on the parameter point $\boldsymbol{\mu}$ have not been stated
explicitly in order to keep the notation as simple as possible.

In the first of equations~\eqref{eq:ALE}, $\mathcal W$ denotes the vector of conservative fluid variables, $x$ denotes
the position vector of the fluid mesh, $\mathcal{J}$ is the determinant of the jacobian of $x$ with respect to the reference
configuration of the fluid mesh, $t$ denotes time, $\nabla$ denotes the gradient with respect to $x$,
$\cdot$ designates the standard dot product, $\mathcal{F}$ denotes the vector of ALE
convective fluxes, and $\mathcal{R}$ denotes the vector of diffusive fluxes.

In the second of equations~\eqref{eq:ALE}, $\mathcal{U}$ denotes the structural displacement vector, $div$ designates the
divergence operator, $\rho$ denotes the material density, $\sigma$ denotes the stress tensor, $\epsilon$ denotes the strain tensor,
and $\mathcal{B}$ is the vector of body forces acting on the structure.

In the third of equations~\eqref{eq:ALE}, the symbols $\tilde{ }$ and $\widetilde{ }$ designate the fictitious aspect of the
deformable body with which the fluid mesh is assimilated.

As already stated, the first two equations described above are coupled by the kinematic and equilibrium transmission
conditions which have simple expressions that can be found, for example, in~\cite{farhat1998load}: therefore, they are not repeated here.
Similarly, the second and third sets of equations are coupled by the continuity of the structural displacement and velocity fields
across the fluid/structure interface.

Finally, each of equations~\eqref{eq:ALE} is equipped with its own boundary and initial conditions: these are not explicitly stated
here for the sake of simplicity.

\subsubsection{Semi-discretization}

In this work, the fluid subsystem is semi-discretized by a finite volume (FV) method: however, it can be equally approximated by a
FE method. On the other hand, the structural and pseudo-structural subsystems are semi-discretized by the FE method. In this case,
the coupled FSI system~\eqref{eq:ALE} is transformed into

\begin{equation}\label{eq:ALEdiscretized}
\left\{
	\begin{array}{l c l}
\dot{\widehat{\left(\textbf{A}(\textbf{x}) \textbf{w} \right)}} + \textbf{F}(\textbf{w}, \textbf{x}, \dot{\textbf{x}}) -
	\textbf{F}^{\prime}(\textbf{w}, \textbf{x}) &=& \textbf{0}\\
	\textbf{M} \ddot{\textbf{u}} +\textbf{ f}^{\,int} (\textbf{u}, \dot{\textbf{u}}) - \textbf{ f}^{\,ext} (\textbf{u}, \textbf{w}) & = & \textbf{0} \\
	\widetilde{\textbf{K}} \textbf{x} - \widetilde{\textbf{K}}_c \textbf{u} & = & \textbf{0}
 \end{array}
 \right.
\end{equation}
In the first of the above semi-discrete equations: $\Abold \in {\mathbb R}^{N_f \times N_f}$ denotes the diagonal matrix storing the
cell volumes of the FV semi-discretization and $N_f$ denotes the dimension of the semi-discrete fluid HDM;
$\Fbold \in {\mathbb R}^{N_f}$ and $\Fbold^{\prime} \in {\mathbb R}^{N_f}$ denote the vectors of semi-discrete convective and
diffusive fluxes, respectively; $\xbold \in {\mathbb R}^{N_m}$ denotes the semi-discrete position vector of the fluid mesh and
$N_m$ denotes the dimension of the associated HDM; $\wbold \in {\mathbb R}^{N_f}$ denotes the vector of semi-discrete conservative
state variables of the fluid subsystem; and $\dot{( )}$ denotes the first derivative with respect to time.

In the second of equations~\eqref{eq:ALEdiscretized}: $\Mbold \in {\mathbb R}^{N_s \times N_s}$ denotes the usual mass matrix and
$N_s$ denotes the dimension of the semi-discrete structural HDM; $\fbold^{\,int} \in {\mathbb R}^{N_s}$ and
$\fbold^{\,ext} \in {\mathbb R}^{N_s}$ denote the vectors of semi-discrete internal and external forces, respectively;
$\ubold \in {\mathbb R}^{N_s}$ denotes the vector of semi-discrete structural state variables (displacements and rotations); and
$\ddot{( )}$ denotes the second derivative with respect to time.

The third of the above equation is a quasi-static semi-discretization of the third of equations~\eqref{eq:ALE}. In this equation,
$\widetilde{\Kbold}$ denotes the pseudo-structural stiffness matrix of the fluid mesh, and $\widetilde{\textbf{K}}_c$ is a
transfer matrix that describes the effect of a structural motion on the fluid mesh motion.

\subsubsection{Linearized CFD-based FSI}

Here, attention is focused on those FSI constraints where for all practical purposes, the fluid can be assumed to be inviscid
($\Fbold^{\prime} = \textbf{0}$ in~\eqref{eq:ALEdiscretized}), and the structure undergoes small rotations and deformations.
These assumptions are realistic for many FSI applications such as those pertaining to flutter, stability, and control. In this case, equations~\eqref{eq:ALEdiscretized}
can be linearized about a static, fluid-structure equilibrium point $(\wbold_0, \xbold_0, \ubold_0)$ characterized by
$\dot{\wbold}_0 = \textbf{0}$, $\dot{\xbold}_0 = \textbf{0}$, and $\dot{\ubold}_0 = \textbf{0}$, following
the approach described in~\cite{lesoinne2001linearized}.
Furthermore, $\xbold$ can be eliminated from the three-way coupled system using the relationship
$\xbold = \widetilde{\textbf{K}}^{-1} \widetilde{\textbf{K}}_c \textbf{u} = \Tbold \ubold$, which follows from the third
of equations~\eqref{eq:ALEdiscretized}. The resulting linearized system of equations can be written as
\begin{equation}\label{eq:ALElinear}
\left\{
\begin{aligned}
& \overline{\Abold}   \dot{\bar{\wbold}} + \Hbold  \bar{{\wbold}}
      + \left( \textbf{E} + \textbf{C} \right) \Tbold \dot{\bar{ \textbf{u}}} +  \textbf{B}
	\, \Tbold {\bar{ \textbf{u}}} &=& \hspace{2mm}\textbf{0}\\
& \Mbold \ddot{\bar{\ubold}} + \Dbold \dot{\bar{\ubold}}
      + \Kbold \bar{\ubold} &=& \hspace{2mm} \Pbold \bar{\wbold}
 \end{aligned}
 \right.
\end{equation}
In the first of the above equations: $\bar{\wbold}$, $\bar{\ubold}$, and $\dot{\bar{\ubold}}$ denote perturbations of
$\wbold$, $\ubold$, and $\dot \ubold$ around the equilibrium point $(\wbold_0, \xbold_0, \ubold_0)$, respectively;
$\overline{\Abold} = \Abold(\xbold_0)$,
${\Hbold} = \displaystyle{\frac{\partial \textbf{F}}{\partial \wbold}} (\wbold_0, \xbold_0, \dot{\xbold}_0 )$,
${\Ebold} = \displaystyle{\frac{\partial \textbf{A}}{\partial \textbf{x}}} (\xbold_0)\,\wbold_0 $,
${\Cbold} = \displaystyle{\frac{\partial \textbf{F}}{\partial \dot{\xbold}}} (\wbold_0, \xbold_0, \dot{\xbold}_0 )$, and
${\Bbold} = \displaystyle{\frac{\partial \textbf{F}}{\partial \textbf{x}}} (\wbold_0, \xbold_0, \dot{\xbold}_0 )$. The
matrices ${\Hbold}$, ${\Ebold}$, ${\Cbold}$, and ${\Bbold}$ arise from the first-order expansion
in Taylor series of the semi-discrete fluid equation around the aforementioned static equilibrium point. In this expansion, the terms
${\Abold} \dot\wbold_0$, ${\Ebold} \dot \xbold_0$,
$\left( \displaystyle{\frac{\partial \Abold}{\partial \xbold}} \bar{\xbold} \right) \bigg\rvert_0 \dot{\wbold}_0$,
$\left( \displaystyle{\frac{\partial \Ebold}{\partial \wbold}} \bar{\wbold} \right) \bigg\rvert_0 \dot{\xbold}_0$, and
$\Fbold (\wbold_0, \xbold_0, \dot{\xbold}_0 )$ vanish because $\dot{\wbold}_0 = \textbf{0}$ and $\dot{\xbold}_0 = \textbf{0}$.

In the second of equations~\eqref{eq:ALElinear},
${\Dbold}  = \displaystyle{\frac{\partial \textbf{ f}^{\,int} }{\partial {\dot{ \ubold}}}}
( \ubold_0, \dot{ \ubold}_0)$,
${\Kbold} = \displaystyle{\frac{\partial \textbf{ f}^{\,int} }{\partial {\ubold}}} (\ubold_0, \dot{\ubold}_0) -
\frac{\partial \textbf{ f}^{\,ext} }{\partial {\ubold}} (\ubold_0,\wbold_0)$, and
${\Pbold} = \displaystyle{\frac{\partial \textbf{ f}^{\,ext} }{\partial \wbold}} (\ubold_0, \wbold_0)$.
The specific form of this equation results from the linearization about a fluid-structure equilibrium point, and therefore a point
characterized by
$- \Mbold \ddot{\ubold}_0 - \fbold^{\,int}(\ubold_0, \dot{\ubold}_0) + \fbold^{\,ext} (\ubold_0,\wbold_0)  = \boldsymbol{0}$.

In the remainder of this paper, the symbols $\bar{ }$ and $\overline{\phantom{A}}$ are dropped in order to keep the notation as simple
as possible. Hence, the linearized FSI system \eqref{eq:ALElinear} is rewritten as follows

\begin{equation}\label{eq:ALElinear1}
\left\{
\begin{aligned}
& \Abold {\dot \wbold} + \Hbold  {\wbold}
      + \Rbold {\dot\ubold} + \Gbold {\ubold} &=& \hspace{2mm}\textbf{0}\\
& \Mbold \ddot\ubold + \Dbold \dot\ubold
      + \Kbold \ubold &=& \hspace{2mm} \Pbold \wbold
 \end{aligned}
 \right.
\end{equation}
where: $\Abold \in {\mathbb R}^{N_f\times N_f}$, $\Hbold \in {\mathbb R}^{N_f \times N_f}$,
$\Rbold = (\Ebold+\Cbold)\,\Tbold \in {\mathbb R}^{N_f \times N_s}$, $\Gbold = \Bbold \Tbold \in {\mathbb R}^{N_f \times N_s}$,
$\wbold \in {\mathbb R}^{N_f}$, $\ubold \in {\mathbb R}^{N_s}$; and $\Mbold \in {\mathbb R}^{N_s \times N_s}$,
$\Dbold \in {\mathbb R}^{N_s \times N_s}$, $\Kbold \in {\mathbb R}^{N_s \times N_s}$, and $\Pbold \in {\mathbb R}^{N_s \times N_f}$.

Let
\begin{equation}\label{eq:MISS1}
	\mathbfcal{A}=
\left[\begin{array}{ccc} \textbf{A} & \textbf{0}_{N_f, N_s} & \textbf{0}_{N_f, N_s}
                             \\ \textbf{0}_{N_s, N_f}  & \textbf{M} &  \textbf{0}_{N_s, N_s}
\\ \textbf{0}_{N_s,N_f} & \textbf{0}_{N_s, N_s} & \textbf{M} \end{array}\right] \in \Rbb^{N_q \times N_q}, \quad
			     \mathbfcal{B} = \left[\begin{array}{ccc} \textbf{H} & \textbf{R} & \textbf{G}
                             \\  -\textbf{P} & \textbf{D} & \textbf{K}
			     \\ \ \textbf{0}_{N_s, N_f} & - \textbf{M} & \textbf{0}_{N_s, N_s} \end{array}\right] \in \Rbb^{N_q \times N_q}, \quad
			     \hbox{and} \quad \textbf{q} = \left[\begin{array}{c} \textbf{w} \\  \dot{\textbf{{u}}} \\
			     \textbf{u} \end{array}\right] \in \Rbb^{N_q}
\end{equation}
where \begin{equation*}
	N_q = N_f + 2 N_s
\end{equation*}

The linearized FSI system~\eqref{eq:ALElinear1} can be rewritten in more compact form as

  \begin{equation}\label{eq:linearPDE}
	  \widetilde{\Rbold}_{L}(\textbf{q}(\mubold),\mubold) = \mathbfcal{A}(\mubold)\, \dot{\textbf{q}}(\mubold) + \mathbfcal{B}(\mubold)\, \textbf{q}(\mubold) = \boldsymbol{0}
      \end{equation}

      where all dependences on the parameter point $\mubold$ have been recalled.

\section{Active Subspace and Database of Linearized FSI PROMs}
\label{sec:PMOR}

Solving FSI-constrained MDAO problems in a high-dimensional design parameter space $\mathcal D$ can be computationally intensive, even when the FSI constraints are linearized. Here, the computational
expense associated with the enforcement of the parametric, linearized, FSI constraint $\widetilde{\Rbold}_{L}(\textbf{q}(\mubold),\mubold) = \boldsymbol{0}$~\eqref{eq:linearPDE} is reduced by substituting
this constraint with a linear, $\mubold$-parametric, FSI PROM, and treating the dependence of this PROM on the parameter point $\mubold$ using the concept of a database ${\mathcal DB}$ of local
PROMs~\cite{amasallem2010towards}. When $\mathcal D$ is relatively high-dimensional, the feasibility of this concept is achieved by significantly reducing the computational cost associated with the
offline construction of the database ${\mathcal DB}$, as follows. First, a low-dimensional space of design parameter generalized coordinates ${\mathcal G} \subset \Rbb^{n_{\mathcal G}}$
is generated -- with $n_{\mathcal G} \ll N_{\mathcal D}$ -- using a new take on the concept of an AS~\cite{lukaczyk2014active}. Then, ${\mathcal G}$ is
sampled using a greedy procedure, and an FSI PROM is constructed at each design parameter point $\mubold_i \in \mathcal D$ associated with each generalized coordinates parameter point $\mubold_{r_i} \in \mathcal G$
sampled in ${\mathcal G}$, and stored in the database ${\mathcal DB}$. Finally, at each unsampled parameter point $\mubold_{\star} \in \mathcal D$ queried by the chosen optimization procedure, a local
PROM counterpart of~\eqref{eq:linearPDE} is constructed in real-time by interpolation on matrix manifolds~\cite{amsallem2011online}, and applied to the real-time solution of the reduced-order counterpart of
the linearized, FSI constraint~\eqref{eq:linearPDE}.

The key elements of the approach outlined above for reducing the computational cost associated with the enforcement of the high-dimensional, linearized, FSI constraint~\eqref{eq:linearPDE} in order to
accelerate the solution of the MDAO problem~\eqref{eq:nandPb} are discussed below.

\subsection{Active Subspace Approximation}
\label{sec:AS}

\subsubsection{Background}
\label{sec:BACK}

In a gradient-based iterative procedure (or algorithm) for the solution of {\it unconstrained} optimization problems, the vector of optimization parameters $\mubold$ is typically updated at each iteration $k$ as follows
\begin{equation}\label{eq:updateOPT}
\mubold^{k+1} = \mubold^k - \eta^k \Cbold^k \nabla f(\mubold^k)
\end{equation}
where $\eta^k$ is a scalar parameter often referred to as the step size, $\Cbold^k \in \Rbb^{N_{\mathcal D} \times N_{\mathcal D}}$ is an algorithm-dependent matrix, and
$\nabla f(\mubold^k) \in \Rbb^{N_{\mathcal D}}$ is the gradient of the scalar objective function $f$ with respect to the vector of optimization parameters $\mubold$ evaluated at $\mubold^k$.
From~\eqref{eq:updateOPT},
it follows that $\mubold^{k+1}$ is computed in the subspace generated by the set of gradient vectors ${\mathcal S}^k = \left \{\nabla f(\mubold^1), \ldots, \nabla f(\mubold^k)\right \}$.
Following the standard ideas of model reduction, the AS technique seeks to approximate each iterate $\mubold^k$ in a low-rank matrix representation of ${\mathcal S}^{k-1}$.
This can be written as
\begin{equation}\label{eq:subspaceAssumption}
\mubold^k \approx \Vbold_{\mu} \, \mubold_r^k
\end{equation}
where $\Vbold_{\mu} \in \Rbb^{N_{\mathcal D} \times n_{\mathcal G}}$ is a basis of dimension $n_{\mathcal G} \ll N_{\mathcal D}$, and $\mubold_{r} \in \Rbb^{n_{\mathcal G}}$ is the vector of
generalized coordinates of $\mubold$ -- or in AS parlance, the vector of active parameters. Specifically, $\Vbold_{\mu}$ is considered to be a global ROB for $\mubold$, and therefore the
approximation~\eqref{eq:subspaceAssumption} is assumed to be accurate in the entire design parameter space $\mathcal D$. Typically, this ROB is constructed using the proper orthogonal decomposition (POD)
method of snapshots~\cite{sirovich1987turbulence} or its singular value decomposition (SVD)-based variant. This calls for sampling the gradients of $f(\cdot)$ at various parameter points $\mubold_i \in {\mathcal D}$,
collecting the sampled values in a snapshot matrix
  \begin{equation}\label{eq:snapX}
  \Sbold =
\begin{bmatrix}
	\nabla f (\mubold_{1}),   \ldots ,  \nabla f (\mubold_{N_{\Sbold}})  \in R^{N_{{\mathcal D}} \times N_{\Sbold}}
\end{bmatrix}
  \end{equation}
where $N_{\Sbold}$ denotes the number of aforementioned sampled gradients, and compressing $\Sbold$ using SVD in order to obtain $\Vbold_{\mu}$. This process is summarized in Algorithm~\ref{alg:POD}.

  \begin{algorithm}
    \caption{Compression of the matrix of gradient snapshots.}
    \label{alg:POD}
  \begin{algorithmic}[1]
	  \INPUT Snapshot matrix $\Sbold \in \Rbb^{N_{\mathcal D} \times N_{\Sbold}}$, and tolerance $\varepsilon$.
	  \OUTPUT $n_{\mathcal G}$, and ROB $\Vbold_{\mu}$ of dimension $n_{\mathcal G}$.
	  \STATE Compute the thin SVD of $\Sbold$: $\Sbold = \Ubold \Sigmabold \Dbold$, where $\Ubold = \left[ \ubold_1 \,\,\,\, \ubold_2 \,\,\,\, \ldots \,\,\,\, \ubold_{rank}  \right]$,
	  $rank$ denotes the rank of $\Sbold$, \break and
	  $\Sigmabold = \left [\begin{array} {c c c c c }
		                       \sigma_1 & 0        & 0      & \ldots           & 0 \\
	                               0        & \sigma_2 & 0      & \ldots           & 0 \\
	                               0        & 0        & \ddots & 0                & 0 \\
				       0        & \ldots   & 0      & \sigma_{rank-1}  & 0\\
				       0        & 0        & \ldots & 0                & \sigma_{rank}
	  \end{array}
	  \right ]$\,, where $\sigma_1 \ge \sigma_2 \ge \ldots \ge \sigma_{rank} > 0$
	  \STATE $n_{\mathcal G}$ = minimum integer for which $\displaystyle{\frac{\sum\limits_{j = n_{\mathcal G}+1}^{rank}\sigma_j^2}{\sum\limits_{j=1}^{rank}\sigma_j^2}} \le \varepsilon$
	  \STATE $\Vbold_{\mu} = \left[ \ubold_1 \,\,\,\, \ubold_2  \,\,\,\, \ldots \,\,\,\, \ubold_{n_{\mathcal G}}  \right]$
  \end{algorithmic}
\end{algorithm}

The sampling underlying the computation of the snapshot matrix~\eqref{eq:snapX} is usually performed using a nonadaptive algorithm. For example, the Latin Hypercube Sampling (LHS)~\cite{mckay1979comparison} method
can be used for this purpose. Typically, the number of samples $N_g$ is logarithmically increased with the dimension $N_{\mathcal D}$ of the design parameter space $\mathcal D$. For example,
reference~\cite{constantine2014computing} recommends

    \begin{equation}\label{eq:numberGRAD}
	    N_{\Sbold} = \alpha \beta \log N_{\mathcal D}
  \end{equation}
where $\alpha$ is a free parameter referred to as the oversampling factor, and $\beta$ is another free parameter. While the values of $\alpha$ and $\beta$ are problem-dependent, reference~\cite{constantine2014computing} suggests
a value between 2 and 10 for $\alpha$, and a value larger than $n_{\mathcal G} + 1$ for $\beta$.

Substituting the AS approximation~\eqref{eq:subspaceAssumption} into~\eqref{eq:nandPb} leads to the following MDAO problem which features a reduced number of optimization parameters

\begin{equation}\label{eq:nandPb1}
        \begin{aligned}
		& \underset{\mubold_r \in {\mathcal G}}{\text{min}}
  		& & f\left(\textbf{q}(\Vbold_{\mu} \, \mubold_r), \zbold(\Vbold_{\mu} \, \mubold_r), \Vbold_{\mu} \, \mubold_r\right) \\
  		& \text{s.t.} & &  \cbold\left(\textbf{q}(\Vbold_{\mu} \, \mubold{_r}), \zbold(\Vbold_{\mu} \, \mubold_r), \Vbold_{\mu} \, \mubold_r\right ) \leq \boldsymbol{0} \\
 \end{aligned}
  \rightarrow
  \boxed{
         \begin{aligned}
		 & \underset{\mubold_r \in {\mathcal G}}{\text{min}}
 		& & \tilde{f}\left(\textbf{q}(\mubold_r), \zbold(\mubold_r), \mubold_r\right) \\
 		& \text{s.t.} & &  \tilde{\cbold}\left(\textbf{q}(\mubold_r), \zbold(\mubold_r), \mubold_r\right) \leq \boldsymbol{0} \\
  \end{aligned}
  }
\end{equation}

The main advantage of solving the MDAO problem~\eqref{eq:nandPb1}, which is based on the AS approximation~\eqref{eq:subspaceAssumption}, instead of solving
the original MDAO problem~\eqref{eq:nandPb} is two-fold:
\begin{itemize}
	\item First, and most importantly, mitigating the curse of dimensionality associated with the construction of a database of local, linear FSI PROMs, when $\mathcal D$ is high-dimensional,
		since in this case parameter sampling is performed in a space of design parameter generalized coordinates ${\mathcal G}$ of much lower dimension than the design parameter space
		$\mathcal D$ on which it is based.
	\item Second, lower the overall computational cost of the iterative optimization procedure by working with a smaller number of optimization parameters, without necessarily sacrificing the
		optimality of the computed solution.
\end{itemize}

As outlined above however, the construction of the AS suffers from two main drawbacks:
\begin{itemize}
	\item The guideline~\eqref{eq:numberGRAD} for estimating the number $N_{\Sbold}$ of gradients to sample at various parameter points $\mubold_i \in {\mathcal D}$ and collect in the snapshot matrix
		$\Sbold$ is an adhoc estimate. Furthermore, because $\alpha$ and $\beta$ are free, problem-dependent algorithmic parameters, this estimate may still call for sampling a very large
		number of parameter points in $\mathcal D$, which defeats the main purpose of the AS concept.
	\item Using a nonadaptive algorithm such as LHS for sampling the gradients of $f(\cdot)$ does not necessarily lead to an appropriate AS for the solution of the MDAO problem~\eqref{eq:nandPb1},
		particularly when $N_{\Sbold}$ itself is limited by the curse of dimensionality. For example, it is shown in Section~\ref{sec:APP} that for one of the two applications considered in this paper, the
		solution of the MDAO problem~\eqref{eq:nandPb1} based on an AS constructed as outlined above is not an optimal solution of the orginal MDAO problem~\eqref{eq:nandPb}.
\end{itemize}

For this reason, a new take on the concept of an AS for the solution on an MDAO problem such as~\eqref{eq:nandPb1} is proposed next.

\subsubsection{New Take on the Concept of an AS and its Construction}
\label{sec:ASNT}

In order to address the drawbacks of the original concept of an AS for the solution of a highly parameterized MDAO problem highlighted in Section~\ref{sec:BACK}, an alternative approach for defining
and constructing an AS is proposed here. This approach is based on the observation that when solving a {\it constrained} MDAO problem such as~\eqref{eq:nandPb}, the update of each iterate $\mubold^k$
performed by an optimization procedure can be written as
 \begin{equation}\label{eq:updateOPT1}
\mubold^{k+1} = \mubold^k + \Delta \mubold^k
\end{equation}
where $\Delta \mubold^k \in \Rbb^{N_{\mathcal D}}$ is determined by the optimization procedure at its $k$-th iteration -- and therefore is optimization-procedure-dependent. Note that
expression~\eqref{eq:updateOPT1} is more general than its counterpart~\eqref{eq:updateOPT}. Specifically:
\begin{itemize}
\item When solving an unconstrained optimization problem using a gradient-based procedure, $\Delta \mubold^k = -\eta^k \Cbold^k \nabla f(\mubold^k)$, where the scalar $\eta^k$ and the matrix
	$\Cbold^k$ are algorithm-dependent. For example, $\Cbold^k=\Ibold$ for the steepest descent algorithm,  $\Cbold^k=\nabla^2f(\mubold^k)$ for Newton's method, and $\Cbold^k=\Bbold^k$ for the
		Broyden-Fletcher-Goldfarb-Shanno (BFGS) algorithm, where $\Bbold^k$ is an approximation of the Hessian $\nabla^2f(\mubold^k)$. When using instead a non gradient-based optimization
		procedure such as a genetic algorithm, $\Delta \mubold^k$ is determined by some stochastic process and/or a probability model.
\item When solving a constrained optimization problem, $\Delta \mubold^k$ depends on whether the optimization procedure is of the interior point or active set type. For example, if the optimization
	problem is solved using a sequential linear programming (SLP) or sequential quadratic programming (SQP) algorithm, $\Delta \mubold^k$ is determined by the iterative solution of a sequence of
		problems formulated for the purpose of satisfying the Karush-Kuhn-Tucker (KKT) conditions.
\end{itemize}

Hence, depending on the chosen optimization procedure, $\Delta \mubold^k$ can have different forms pertaining to the gradient or Hessian of $f(\cdot)$. For this reason, an alternative
concept of an AS for the solution of a highly parameterized MDAO problem that accounts for how $\Delta \mubold^k$ is computed by the chosen optimization procedure is proposed here. Specifically, this
concept is based on the specific form of the increment $\Delta \mubold$ rather than on $\nabla f(\mubold)$. Consequently, this alternative concept of an AS calls for constructing the ROB $\Vbold_{\mu}$
by collecting and compressing the snapshot matrix
\begin{equation}\label{eq:snapX1}
	{\bf {\mathbb S}} =
\begin{bmatrix}
	\mubold^0,	\Delta \mubold^0,   \ldots ,  \Delta \mubold^{(N_{{\mathbb S}-2})}
\end{bmatrix} \in R^{N_{{\mathcal D}} \times N_{\mathbb S}}
\end{equation}
where $\mubold^0$ is the same initial guess as for the solution of the original MDAO problem~\eqref{eq:nandPb}, and $N_{\mathbb S}$ denotes the total number of snapshots collected in $\bf {\mathbb S}$.

Furthermore, since in the context of this work the main purpose of the AS is to mitigate the curse of dimensionality faced by an adaptive procedure such as a greedy procedure
(for example, see~\cite{farhat2019feasible}) for sampling the design
parameter space $\mathcal D$, it is not compelling to construct $\bf{\mathbb S}$~\eqref{eq:snapX1} by sampling the same parameter space $\mathcal D$ using a nonadaptive algorithm! For this reason, it is
proposed here to construct $\bf{\mathbb S}$~\eqref{eq:snapX1} by following instead the trajectory of the solution of the dramatically less computationally intensive version of the original MDAO
problem~\eqref{eq:nandPb}, where the FSI constraint~\eqref{eq:linearPDE} is dropped -- that is, by solving the auxiliary optimization problem
\begin{equation}\label{eq:nandPb2}
         \begin{aligned}
		 & \underset{\mubold \in {\mathcal D}}{\text{min}}
		& & f( \zbold(\mubold), \mubold) \\
		& \text{s.t.} & &  \cbold( \zbold(\mubold), \mubold) \leq \boldsymbol{0} \\
 \end{aligned}
\end{equation}
using the same optimization procedure and same initial guess $\mubold^0$ as for the solution of the original MDAO problem~\eqref{eq:nandPb}, and collecting in $\bf{\mathbb S}$~\eqref{eq:snapX1} the
snapshots $\left \{ \Delta \mubold^k \right \}_{k=0}^{N_{{\mathbb S}-2}}$.

If the aforementioned optimization procedure does not involve the computation of the true Hessian of the objective function, the cost of sampling $\Delta \mubold^k = \mubold^{k+1}-\mubold^{k}$ is roughly
equal to that of sampling $\nabla f(\mubold^k)$. If on the other hand the optimization procedure involves the computation of $\nabla^2f(\mubold^k)$, the additional cost associated with sampling
$\Delta \mubold^k$ instead of $\nabla f(\mubold^k)$ is that incurred by, for example, enforcing the KKT conditions if the formulation of problem~\eqref{eq:nandPb2} incorporates constraints, or performing
the matrix-vector multiplication $\Cbold^k \nabla f(\mubold^k)$ if it does not. In either case, the additional cost is affordable and offsetted by the
following advantages of the proposed concept of an AS for the solution of highly parameterized MDAO problems:
\begin{itemize}
	\item The sampling procedure associated with its construction is adaptive, whereas that associated with the original concept of an AS outlined in Section~\ref{sec:BACK} is typically nonadaptive
		(as in the case of the LHS method highlighted above): therefore for the same number of samples, the sampling procedure associated with the proposed concept of an AS can be expected to
		deliver a more effective ROB $\Vbold_{\mu}$ than its counterpart associated with the original concept of an AS.
	\item The number of snapshots needed for the construction of an effective AS of the type advocated here is automatically determined by the problem-independent convergence criterion governing
		the solution of the auxiliary optimization problem~\eqref{eq:nandPb2}. On the other hand, that needed for the construction of an effective AS of the original type reviewed in
		Section~\ref{sec:BACK} is governed by the adhoc and problem-dependent criterion~\eqref{eq:numberGRAD}, and therefore is determined in practice by a far less convenient trial and error
		approach.
\end{itemize}

\subsection{Parametric PMOR for FSI and Sensitivities}
\label{sec:PPMOR}

PMOR reduces the order (dimension, or size) of an HDM such as the linearized, FSI constraint~\eqref{eq:linearPDE} by:
\begin{itemize}
	\item Performing subspace approximations, which in this case can be written as
		\begin{equation}\label{eq:subspaceAssumption1}
			\wbold = \Vbold_{w} \, \wbold_r, \quad \ubold = \Vbold_{u}\, \ubold_r \quad \Rightarrow \quad
			\qbold = \Vbold_{q} \, \qbold_r = \left[\begin{array}{c c c} \Vbold_{w}     & \boldsymbol{0}_{N_f, n_s} & \boldsymbol{0}_{N_f, n_s} \\
				\boldsymbol{0}_{N_s, n_f} & \Vbold_{u}     & \boldsymbol{0}_{N_s, n_s} \\
			\boldsymbol{0}_{N_s, n_f} & \boldsymbol{0}_{N_s, n_s} & \Vbold_{u}\end{array}\right]
								 	       \left[\begin{array}{c} \textbf{w}_r \\  \dot{\textbf{{u}}}_r \\  \textbf{u}_r \end{array}\right]
		\end{equation}
		where $\Vbold_w \in \Rbb^{N_f\times n_f}$ is a {\it right} fluid ROB with $n_f \ll N_f$, $\Vbold_u \in \Rbb^{N_s\times n_s}$ is a {\it right} structural ROB with $n_s \ll N_s$ -- and
		therefore $\Vbold_q \in \Rbb^{(N_q \times n_q}$ is a fluid-structure ROB with $n_q \ll N_q$ -- $\wbold_r \in \Rbb^{n_f}$,  $\ubold_r \in \Rbb^{n_s}$,
		and $\qbold_r \in \Rbb^{n_q}$ where
		\begin{equation*}
			n_q = n_f + 2 n_s \ll N_q
		\end{equation*}
		Substituting the subspace approximation $\qbold = \Vbold_{q} \, \qbold_r$~\eqref{eq:subspaceAssumption1} in the parametric, linearized, FSI constraint~\eqref{eq:linearPDE} and
		using the AS approximation~\eqref{eq:subspaceAssumption} leads to
		\begin{equation}\label{eq:PROM_unc}
			\mathbfcal{A}(\mubold_r) \Vbold_q \, \dot{\textbf{q}}_r(\mubold_r)  + \mathbfcal{B}(\mubold_r) \Vbold_q \, \textbf{q}_r(\mubold_r) =   \boldsymbol{0}
		\end{equation}
		which is an overdetermined system of equations.
	\item Performing a Galerkin or Petrov-Galerkin projection on each of the three block equations embedded in system~\eqref{eq:PROM_unc}, in order to transform it into a square
		system. Typically, a Galerkin projection -- which does not require the computation of a left ROB -- is justified for the structural subsystem. However, a Petrov-Galerkin projection
		may be preferred for the fluid subsystem for numerical stability reasons~\cite{amsallem2012stabilization}. It follows that in general, the left fluid-structure ROB can be written as
		\begin{equation}\label{eq:leftFSIROB}
			\Wbold_q = \left[\begin{array}{c c c} \Wbold_{w} & \boldsymbol{0}_{N_f, n_s} & \boldsymbol{0}_{N_f, n_s} \\
				\boldsymbol{0}_{N_s, n_f} & \Vbold_{u}     & \boldsymbol{0}_{N_s, n_s} \\
			\boldsymbol{0}_{N_s, n_f} & \boldsymbol{0}_{N_s, n_s} & \Vbold_{u}\end{array}\right]
		\end{equation}
		where $\Wbold_{w} \in \Rbb^{N_f\times n_f}$ is a left fluid ROB. Projecting next the overdetermined system~\eqref{eq:PROM_unc} onto the subspace represented by the left fluid-structure ROB
		$\Wbold_q$~\eqref{eq:leftFSIROB} transforms this system into the square counterpart
		\begin{equation*}\label{eq:PROM1}
			\left(\Wbold_q^T \mathbfcal{A}(\mubold_r) \Vbold_q\right) \, \dot{\textbf{q}}_r(\mubold_r)  + \left(\Wbold_q^T \mathbfcal{B}(\mubold_r) \Vbold_q\right) \, \textbf{q}_r (\mubold_r) = \boldsymbol{0}
		\end{equation*}
		where the superscript $T$ designates the transpose of a quantity. The above reduced system is a parametric, linear, FSI PROM. Using~\eqref{eq:MISS1}, ~\eqref{eq:subspaceAssumption},
		~\eqref{eq:subspaceAssumption1}, and~\eqref{eq:leftFSIROB}, it can be rewritten in a more compact form as follows
		\begin{eqnarray}\label{eq:PROM2}
			\widetilde{\Rbold}_{L_r}(\textbf{q}(\mubold_r),\mubold_r) &=& \mathbfcal{A}_r (\mubold_r) \, \dot{\textbf{q}}_r(\mubold_r)  + \mathbfcal{B}_r (\mubold_r) \, \textbf{q}_r(\mubold_r) = \boldsymbol{0} \nonumber\\
			\hbox{where} \quad \mathbfcal{A}_r (\mubold_r) &=& \Wbold_q^T {\mathbfcal A} (\mubold_r) \Vbold_q = \left [\begin{array}{c c c} \Wbold_w^T \Abold(\mubold_r)\Vbold_w & \boldsymbol{0}_{n_f, n_s} & \boldsymbol{0}_{n_f, n_s}\\
			\boldsymbol{0}_{n_s, n_f} & \Vbold_u^T\Mbold(\mubold_r)\Vbold_u & \boldsymbol{0}_{n_s, n_s} \\ \boldsymbol{0}_{n_s, n_f} & \boldsymbol{0}_{n_s, n_s} & \Vbold_u^T \Mbold(\mubold_r) \Vbold_u \end{array} \right ] \in \Rbb^{n_q \times n_q} \\
			\hbox{and} \quad \mathbfcal{B}_r (\mubold_r) &=& \Wbold_q^T {\mathbfcal B} (\mubold_r) \Vbold_q = \left [\begin{array}{c c c} \Wbold_w^T \Hbold(\mubold_r) \Vbold_w & \Wbold_w^T \Rbold(\mubold_r) \Vbold_u &
				\Wbold_w^T \Gbold(\mubold_r) \Vbold_u \\ -\Vbold_u^T \Pbold(\mubold_r) \Vbold_w & \Vbold_u^T \Dbold(\mubold_r) \Vbold_u & \Vbold_u^T \Kbold(\mubold_r) \Vbold_u \\ \boldsymbol{0}_{n_s, n_f} & -\Vbold_u^T \Mbold(\mubold) \Vbold_u & \boldsymbol{0}_{n_s, n_s}
			\end{array} \right ] \in \Rbb^{n_q \times n_q} \nonumber
		\end{eqnarray}
\end{itemize}

The parametric, linearized, FSI constraint considered this work is a flutter constraint, as flutter remains one of the most important considerations in aircraft design. Consequently,
{\it for a given $\mubold_r \in {\mathcal G}$} -- and therefore, a given parameter point $\mubold = \Vbold_{\mubold}\mubold_r \in {\mathcal D}$ -- the fluid and structural ROBs are computed as follows:

\begin{itemize}
	\item The right fluid ROB $\Vbold_w$ is computed as described in~\cite{lieu2006reduced}. For this purpose, and only for this purpose, the first of the linearized equations~\eqref{eq:ALElinear} is rewritten in
	the frequency domain by assuming a periodic solution of the form $\bar \wbold = \bar \wbold_a e^{I\kappa t}$ and a periodic excitation of the form $\bar \ubold = \bar \ubold_ae^{I\kappa t}$,
		where the subscript $a$ designates the amplitude, $I$ denotes the pure imaginary number satisfying $I^2 = -1$, $\kappa$ denotes the (aeroelastic) reduced frequency of interest, and $t$
		denotes as usual time. Next, a few -- say $N_m$ -- ground-based natural mode shapes of the structural subsystem are considered, together with a reduced frequency band of interest. For
		each considered ground-based mode shape $\bar \ubold_{a_m}$, $m = 1, \ldots N_m$, a sweep is performed in the reduced frequency band of interest, which for this purpose is sampled into
		$N_l$ points to obtain $N_m N_l$ excitation inputs of the form $\bar \ubold_{a_m}e^{I\kappa_lt}$. For each of these excitation inputs, a fluid amplitude solution $\bar \wbold_{a_m}(k_l)$
		is computed by prescribing $\bar \ubold = \bar \ubold_{a_m}e^{I\kappa_lt}$ in the aforementioned frequency domain version of the first of the linearized equations~\eqref{eq:ALElinear} and
		solving this equation. Next, the real and imaginary parts of each computed solution $\bar \wbold_{a_m}(k_l)$ are collected in a matrix of fluid solution snapshots. Finally, the
		$2N_mN_l$ collected fluid snapshots are compressed using SVD to construct the desired fluid ROB $\Vbold_w$ (for example, see Algorithm~\ref{alg:POD}). In summary, this procedure for
		constructing $\Vbold_w$ amounts to training the fluid ROB for different deformed configurations of the structural subsystem flapping at multiple frequencies, in a relevant frequency band.
	\item The left fluid ROB $\Wbold_w$ is computed in two steps as follows. First, this ROB is temporarily set to $\Wbold_w = \Vbold_w$ and the eigenvalues of the reduced fluid matrix
		$\Wbold_w^T\Hbold\Vbold_w$ (see~\eqref{eq:PROM2}) are computed in real-time. If all eigenvalues of this reduced matrix turn out to be positive, $\Wbold_w^T\Hbold\Vbold_w$ is stable
		and the fluid subsystem is reduced using a Galerkin projection ($\Wbold_w = \Vbold_w$). On the other hand, if any eigenvalue of  $\Wbold_w^T\Hbold\Vbold_w$ turns out to be negative,
		$\Wbold_w^T\Hbold\Vbold_w = \Vbold_w^T\Hbold\Vbold_w$ is unstable. In this case, $\Wbold_w$ is updated by adding to it 1 or 2 columns that are constructed using the procedure described
		in~\cite{amsallem2012stabilization}, which provably guarantees the stability of $\Wbold_w^T\Hbold\Vbold_w$, and the fluid is reduced by a Petrov-Galerkin projection $(\Wbold_w \ne \Vbold_w)$.
	\item The right structural ROB $\Vbold_u$ is chosen as a collection of low frequency, ground-based, natural mode shapes, which is justified by the fact that a flutter response is typically
		dominated by these modes. From~\eqref{eq:PROM2}, it follows that if these mode shapes are mass-orthonormalized,
		\begin{eqnarray}
			\label{eq:MISS2}
			\mathbfcal{A}_r (\mubold_r) &=& \Wbold_q^T {\mathbfcal A} (\mubold_r) \Vbold_q = \left [\begin{array}{c c c} \Wbold_w^T \Abold(\mubold_r)\Vbold_w & \boldsymbol{0}_{n_f, n_s} & \boldsymbol{0}_{n_f, n_s}\\
			\boldsymbol{0}_{n_s, n_f} & \Ibold_{n_s, n_s} & \boldsymbol{0}_{n_s, n_s} \\
			\boldsymbol{0}_{n_s, n_f} & \boldsymbol{0}_{n_s, n_s} & \Ibold_{n_s, n_s} \end{array} \right ] \in \Rbb^{n_q \times n_q} \\
			\hbox{and} \quad \mathbfcal{B}_r (\mubold_r) &=& \Wbold_q^T {\mathbfcal B} (\mubold_r) \Vbold_q = \left [\begin{array}{c c c}
				\Wbold_w^T \Hbold(\mubold_r) \Vbold_w & \Wbold_w^T \Rbold(\mubold_r) \Vbold_u & \Wbold_w^T \Gbold(\mubold_r) \Vbold_u \\ -\Vbold_u^T \Pbold(\mubold_r) \Vbold_w & \Vbold_u^T \Dbold(\mubold_r) \Vbold_u & \boldsymbol{\Omega}_r^2(\mubold_r) \\
				\boldsymbol{0}_{n_s, n_f} & -\Ibold_{n_s, n_s} & \boldsymbol{0}_{n_s, n_s} \end{array} \right ] \in \Rbb^{n_q \times n_q} \nonumber
		\end{eqnarray}
		where $\boldsymbol{\Omega}_r^2$ is the diagonal matrix storing the squares of the natural angular frequencies associated with the chosen collection of low frequency, ground-based, natural
		mode shapes.
\end{itemize}

From~\eqref{eq:PROM2}, it also follows that the parametric, linear, FSI PROM associated with the parametric, linearized, FSI constraint~\eqref{eq:linearPDE} can be described by the following tuple of FSI
reduced-order operators
\begin{equation}\label{eq:tuple}
	{\mathcal T}_r^{FSI}(\mubold_r) = \left\{{\mathbfcal A}_r(\mubold_r), {\mathbfcal B}_r(\mubold_r) \right\} = \left\{\Wbold_q^T {\mathbfcal A} (\mubold_r) \Vbold_q, \Wbold_q^T {\mathbfcal B} (\mubold_r) \Vbold_q \right\}
\end{equation}
However, the solution of the MDAO problem~\eqref{eq:nandPb1} using a gradient-based optimization algorithm and the parametric, PROM counterpart~\eqref{eq:PROM2} of the parametric, FSI
constraint~\eqref{eq:linearPDE} also requires the computation of the sensitivities $\displaystyle{\frac{\partial \mathbfcal{A}_r}{\partial \mu}}(\mubold_r)$ and
$\displaystyle{\frac{\partial \mathbfcal{B}_r}{\partial \mu}(\mubold_r)}$ at each point $\mubold_{r_i} \in {\mathcal G}$ (and therefore parameter point
$\mubold_i = \Vbold_{\mubold}\mubold_{r_i} \in {\mathcal D}$)
visited by the optimization algorithm. For this reason, the extended tuple
\begin{equation}\label{eq:exttuple}
	\overline{\mathcal T}_r^{FSI}(\mubold_r) = \left\{{\mathbfcal A}_r(\mubold_r), {\mathbfcal B}_r(\mubold_r),
	\left\{\displaystyle{\frac{\partial \mathbfcal{A}_r(\mubold_r)}{\partial \mu[j]}}, \displaystyle{\frac{\partial \mathbfcal{B}_r(\mubold_r)}{\partial \mu[j]}} \right\}_{j=1}^{N_{\mathcal D}}  \right\}
\end{equation}
where $\mubold[j]$ denotes the $j$-th component of the parameter point $\mubold \in \mathcal D \subset {\mathbb R}^{N_{\mathcal D}}$ is also introduced.



\subsection{Database of Consistent Pointwise Linearized FSI PROMs}

The linear, FSI PROM~\eqref{eq:PROM2} is $\mubold$-parametric. As already stated, the approach chosen here for treating the parameter dependence of such a PROM consists in: constructing offline a
database ${\mathcal DB}$ of local tuples of FSI reduced-order operators ${\mathcal T}_r(\mubold_{r_i}) = \left\{{\mathbfcal A}_r(\mubold_{r_i}), {\mathbfcal B}_r(\mubold_{r_i}) \right\}$, $i = 1, \ldots,
N_{\mathcal DB}$; and interpolating these online on appropriate matrix manifolds to generate in real-time the tuple ${\mathcal T}_r(\mubold_{r_{\star}})$ at each parameter point
$\mubold_{\star} = \Vbold_{\mubold}\mubold_{r_{\star}}$ queried by the optimization algorithm, but where ${\mathcal T}_r(\mubold_{r_{\star}})$ is not available in ${\mathcal DB}$. This approach requires
first addressing two important issues:
\begin{itemize}
	\item Selecting the number $N_{\mathcal DB}$ and locations in ${\mathcal G}$ of the points $\mubold_{r_i}$ -- and therefore locations in ${\mathcal D}$ of the parameter points
		$\mubold_i = \Vbold_{\mubold}\mubold_{r_i}$ -- where to construct the local tuples of reduced-order operators defining the local, linear PROMs of interest.
	\item Ensuring that all constructed tuples stored in ${\mathcal DB}$ are consistent in the sense defined in~\cite{amsallem2011online} and examplified as well as explained below.
\end{itemize}

The first issue highlighted above defines the {\it parameter sampling} problem. It is common to many local and global approaches for treating the parameter dependence of a PROM -- that is, for
{\it training} a PROM in ${\mathcal G}$ (here) or in $\mathcal D$ (in general).
Among all elements of the approach described in this paper for reducing the computational cost associated with the enforcement of the high-dimensional, linearized,
FSI constraint~\eqref{eq:linearPDE} in order to accelerate the solution of the MDAO problem~\eqref{eq:nandPb1} (here) or~\eqref{eq:nandPb} (in general),
it is the most vulnerable to the curse of dimensionality -- specifically, the dimension
$N_{\mathcal D}$ of $\mathcal D$. This issue is treated in Section~\ref{sec:FEASAM}.

The second issue highlighted above can be illustrated as follows. Consider, for example, the linear, structural PROM associated with the linear, FSI PROM~\eqref{eq:PROM2} and describable by the
sub-tuple of ${\mathcal T}_r^{FSI}(\mubold_r)$~\eqref{eq:tuple}
\begin{eqnarray*}
	{\mathcal T}_r^{S}(\mubold_r) &=& \left\{\Mbold_r(\mubold_r), \Dbold_r(\mubold_r), \Kbold_r(\mubold_r), \Pbold_r(\mubold_r) \right\} \\
	\hbox{where} \quad \Mbold_r(\mubold_r) &=& \Vbold_u^T(\mubold_r) \Mbold(\mubold_r) \Vbold_u(\mubold_r), \quad \Dbold_r(\mubold_r) = \Vbold_u^T(\mubold_r) \Dbold(\mubold_r) \Vbold_u(\mubold_r)\\
	\Kbold_r(\mubold_r) &=&  \Vbold_u^T(\mubold_r) \Kbold(\mubold_r) \Vbold_u(\mubold_r), \quad \Pbold_r(\mubold_r) = \Vbold_u^T(\mubold_r) \Pbold(\mubold_r) \Vbold_w(\mubold_r)
\end{eqnarray*}
and the local aspect of the right structural ROB is emphasized by the notation $\Vbold_u(\mubold_r)$. Furthermore, consider the case where this ROB is made of low frequency, ground-based, natural mode
shapes of the structural subsystem, and the design parameter space $\mathcal D$ includes material properties such as the structural Young modulus and/or structural density. It is known that in this case
(for example, see~\cite{amsallem2011online}), the ordering of the natural mode shapes is material property dependent. Consequently, it is possible that for $\mubold_r = \mubold_{r_1}$
($\mubold = \Vbold_{\mubold}\mubold_{r_1} = \mubold_1$), the first 3 columns of $\Vbold_u(\mubold_{r_1})$ correspond to -- for example -- the first bending mode, second-bending mode, and first torsion
mode of the structural subsystem, respectively, while for $\mubold_r = \mubold_{r_2}$ ($\mubold = \Vbold_{\mubold}\mubold_{r_2} = \mubold_2$), they correspond to the first bending mode, first torsion
mode, and second bending mode, respectively. In this case, ${\mathcal T}_r^{S}(\mubold_{r_1})$ and ${\mathcal T}_r^{S}(\mubold_{r_2})$ are said to be inconsistent, because their underlying ROBs
$\Vbold_u(\mubold_{r_1})$ and $\Vbold_u(\mubold_{r_2})$ are inconsistent in the sense that algebraic manipulations performed on these ROBs -- and therefore, on their associated structural PROMs -- such as
linear combinations and interpolations are not physically meaningful, unless the columns of $\Vbold_u(\mubold_{r_1})$ and $\Vbold(\mubold_{r_2})$ are first re-ordered to be physically consistent (for
example, first bending, second bending, and then first torsion, etc., for both $\Vbold_u(\mubold_{r_1})$ and $\Vbold(\mubold_{r_2})$). Note however that the ordering of the ROBs is only one source of
inconsistency among many others.

To enforce consistency between the pre-computed tuples ${\mathcal T}_r^{FSI}(\mubold_{r_i}) = \left\{{\mathbfcal A}_r(\mubold_{r_i}), {\mathbfcal B}_r(\mubold_{r_i}) \right\} =
\left\{\Wbold_q^T(\mubold_{r_i}) {\mathbfcal A} (\mubold_{r_i}) \Vbold_q (\mubold_{r_i}), \Wbold_q^T (\mubold_{r_i}) {\mathbfcal B} (\mubold_{r_i}) \Vbold_q (\mubold_{r_i}) \right\}$,
$i = 1, \ldots, N_{\mathcal{DB}}$, it is first noted that any ROB such as $\Vbold_q$ remains a ROB representing the same subspace approximation after post-multiplication by an orthogonal matrix
$\Qbold \in \Rbb^{n_q \times n_q}$ -- that is, after ``rotation'' of this ROB. In other words, $\Vbold_q $ and $\Vbold_q \Qbold$ represent the same subspace when $\Qbold$ is an
orthogonal matrix (for example, $\Qbold^T\Qbold = \Ibold_{n_q, n_q}$). This illustrates the following facts:
\begin{itemize}
	\item There is no unique ROB associated with a given subspace approximation.
	\item Two different ROBs such as $\Vbold_q$ and $\Vbold_q\Qbold$ associated with the same subspace approximation define two different generalized coordinates systems for performing the same
		approximation.
\end{itemize}

It follows that one methodology for enforcing consistency between several pre-computed PROMs is to ensure that their underlying ROBs are constructed for the same generalized coordinates system~\cite{amsallem2011online, amsallem2016real}. This methodology consists in selecting one of the pre-computed ROBs as a reference basis, for example, $\Vbold_q(\mubold_{r_{ref}}) = \Vbold_q(\mubold_{r_1})$, and finding for
each other pre-computed ROB $\Vbold_q(\mubold_{r_i})$ and associated tuple the orthogonal transformation matrix $\Qbold_i$ that solves the following minimization problem

\begin{equation}\label{eq:Procrustes}
	\underset{\textbf{Q}\, \in\, {\mathcal O}(n_q)}{\text{min}} \| \Vbold_q(\mubold_{r_{ref}})  - \Vbold_q(\mubold_{r_i})\Qbold\|_F
\end{equation}
where ${\mathcal O}(n_q)$ denotes the set of orthogonal matrices of size $n_q$, and the subscript $F$ designates the Frobenius norm.

There are two reasons for $\Vbold_q(\mubold_{r_i})$ to be different from $\Vbold_q(\mubold_{r_{ref}})$: $\mubold_{r_i} \ne \mubold_{r_{ref}}$; and $\Vbold_q(\mubold_{r_{ref}})$ and
$\Vbold_q(\mubold_{r_i})$ may define two different generalized coordinates systems. By solving the minimization problem~\eqref{eq:Procrustes}, the orthogonal matrix $\Qbold_i$ is determined so that
$\mubold_{r_i} \ne \mubold_{r_{ref}}$ is the only reason why $\Vbold_q(\mubold_{r_i})\Qbold_i \ne \Vbold_q(\mubold_{r_{ref}})$, and therefore the 2 ROBs are consistent in the sense defined above.

Problem~\eqref{eq:Procrustes} is known as the orthogonal Procustes problem~\cite{schonemann1966generalized}. Most importantly, it has an analytical solution that can be computed in real-time using
Algorithm~\ref{alg:rotationV}. This solution, $\Qbold_i$, is also a congruence transformation matrix that can be used to enforce the consistency of each tuple  ${\mathcal T}_r^{FSI}(\mubold_{r_i})$,
$i=1, \ldots, N_{\mathcal DB}$, with the reference tuple ${\mathcal T}_r^{FSI}(\mubold_{r_{ref}})$ by transforming this tuple into
\begin{equation}\label{eq:tupleROT}
	\mathcal{C}_r^{FSI}(\mubold_{r_i}) = \left\{\mathbfcal{A}_r^C(\mubold_{r_i}), \mathbfcal{B}_r^C(\mubold_{r_i}) \right\} = \left\{ \Qbold_i^T \Vbold_q^T (\mubold_{r_i})  \mathbfcal{A}(\mubold_{r_i}) \Vbold_q (\mubold_{r_i}) \Qbold_i,  \Qbold_i^T \Vbold_q^T (\mubold_{r_i}) \mathbfcal{B} \Vbold_q (\mubold_{r_i}) \Qbold_i \right\} = \left\{ \Qbold_i^T {\mathbfcal A}_r(\mubold_{r_i}) \Qbold_i,
	\Qbold_i^T {\mathbfcal B}_r(\mubold_{r_i}) \Qbold_i \right \}
\end{equation}

Hence, after ${\mathcal DB}$ is constructed, and each time it is updated, its consistency is enforced by selecting arbitrarily a reference basis and applying Algorithm~\ref{alg:rotationV}
to perform congruence transformations on its pre-computed tuples of FSI reduced-operators.

\begin{algorithm}
\caption{Congruence transformation of a tuple of reduced-order operators for enforcing consistency.}
\label{alg:rotationV}
\begin{algorithmic}[1]
	\INPUT Reference point $\mubold_{r_{ref}}$, associated fluid-structure ROB $ \Vbold_q (\mubold_{r_{ref}})$, and another ROB $ \Vbold_q (\mubold_{r_{i}})$ where
	$\mubold_{r_i} \ne \mubold_{r_{ref}}$.
\OUTPUT Optimal congruence transformation matrix $\Qbold_i$.
\STATE Compute SVD of $\Vbold_q(\mubold_{r_{i}})^T\Vbold_q(\mubold_{r_{ref}}) =\Ubold \Sigmabold \Zbold^T $
\STATE Compute $\Qbold_i = \Ubold \Zbold^T$
\end{algorithmic}
\end{algorithm}


\subsection{Feasible Adaptive Parameter Sampling}\label{par:as}
\label{sec:FEASAM}

When the design parameter space $\mathcal D$ is high-dimensional, the offline pre-computation of a database of linear, FSI PROMs may not be feasible. This, because it may require sampling a large number of parameter
points $\mubold_i \in \mathcal D$ where to pre-compute linear FSI PROMs, in order to enable the accurate interpolation at a queried but unsampled parameter point of a tuple of FSI reduced-operators such as~\eqref{eq:tuple}.
In this case, the new take on the AS approach described in Section~\ref{sec:ASNT} can be used to find a lower-dimensional subspace in which to solve the FSI-constrained MDAO problem~\eqref{eq:nandPb1}.
After this subspace is determined and represented by the ROB $\Vbold_{\mu}$ as explained in Section~\ref{sec:AS}, parameter sampling can be performed in the lower-dimensional space of generalized
coordinates $\mathcal G$. For this purpose, two different approaches can be pursued:
\begin{itemize}
\item {\it A priori sampling.} This appropach samples points in the target space either randomly, or according to a nonadaptive, pre-designed scheme. Examples include the full factorial sampling, random sampling, and
LHS methods. Such sampling methods are not optimal, because they have no explicit awareness of where the interpolant will be inaccurate. For this reason, they may require a larger than necessary number of sampled points
in order to deliver the expected accuracy at interpolation time: as such, they may lead to unaffordable databases of linear PROMs when $N_{\mathcal D}$ is relatively high.
\item {\it Adaptive sampling}. This alternative approach is typically iterative. It requires the availability of an error estimator or indicator for the PROM. It samples points in regions of the target space where the current
instance of the database is assessed by the error estimator/indicator to be inaccurate. Consequently, it avoids over-sampling: in this sense, it reduces the size $N_{\mathcal DB}$ of the needed database and makes its
construction affordable. For this reason, this parameter sampling approach is chosen here for sampling $\mathcal G$.
\end{itemize}

Specifically, the indirect sampling of the design parameter space $\mathcal D$ is performed here by directly sampling the parameter space $\mathcal G$ using an iterative, adaptive sampling algorithm equipped with:
\begin{itemize}
\item A pre-selected set of candidate sample points $\Xi_r$ that is large enough to faithfully represent $\mathcal G$.
\item A residual-based error indicator $e(\mubold_r, {\mathcal DB})$.
\end{itemize}
This algorithm can be described as a greedy procedure that samples at each $i$-th iteration the parameter point $\mubold_{r_i} \in \Xi_r$ where some norm of $e(\mubold_r, {\mathcal DB})$ is maximized, and builds at this point the
tuple of FSI reduced-operators ${\mathcal T}_r^{FSI}(\mubold_{r_i})$~\eqref{eq:tuple}. Hence, this algorithm leads to an incremental construction of $\mathcal DB$ that can be written as
\begin{equation}\label{eq:greedy}
	\mathcal{DB}_i = \left\{  \mathcal{DB}_{i-1}, \Vbold_q(\mubold_{r_i}), \mathcal{T}_r(\mubold_{r_i}) \right\}, \quad \text{ where } \quad \mubold_{r_i} =  \arg \max_{\mubold_r \, \in \, \Xi_r} e(\mubold_r, \mathcal{DB}_{i-1})
\end{equation}
$\mathcal{DB}_i$ is the instance of the database of linear FSI PROMs $\mathcal DB$ at iterations $i$, and the storage in the database $\mathcal DB$ of each pre-computed local ROB $\Vbold_q(\mubold_{r_i})$
is justified below.

In general, the MDAO problem of interest may contain bounding constraints on the design optimization parameters that form $\mubold$ -- for example, $\cbold\left(\textbf{q}(\Vbold_{\mu} \, \mubold{_r}), \zbold(\Vbold_{\mu} \, \mubold_r), \Vbold_{\mu} \, \mubold_r\right ) \leq \boldsymbol{0}$ in~\eqref{eq:nandPb1} may contain bounding constraints on $\mubold$. In this case, each sampled parameter point of the form $\mubold_i = \Vbold_{\mu} \mubold_{r_i}$ must satisfy the {\it feasibility} constraint
    \begin{equation}\label{eq:bound}
	    \mubold_{\text{lb}} \leq \Vbold_{\mu} \mubold_{r_i} \leq \mubold_{\text{ub}}
  \end{equation}
where $\mubold_{\text{lb}} \in \Rbb^{N_{\mathcal D}}$ and $\mubold_{\text{ub}} \in \Rbb^{N_{\mathcal D}}$ are the lower and upper bounds on $\mubold$, respectively. The case where no upper or lower bounding constraint on $\mubold$ is
specified can be simply accommodated by setting the entries of $\mubold_{\text{ub}}$ and/or $\mubold_{\text{lb}}$, as needed, to $\pm \infty$ (very large numbers in practice).

It follows that {\it all} candidate parameter points $\mubold_{r_i} \in \Xi_r$ must also satisfy the constraint ~\eqref{eq:bound}. This can be achieved by constructing $\Xi_r$ using
Algorithm~\ref{alg:feasibleAS} which samples points in $\mathcal G$ that are feasible in $\mathcal D$ by solving the feasibility problem~\eqref{eq:FEAS}. It can be initialized using a single design
parameter point in the center of the AS, or a small number of design parameter points randomly chosen in the AS.

  \begin{algorithm}
    \caption{Construction of a set of candidate parameter points $\Xi_r$ using a feasible active subspace algorithm.}
    \label{alg:feasibleAS}
  \begin{algorithmic}[1]
	  \INPUT Bounding vectors $\mubold_{\text{lb}}$ and $\mubold_{\text{ub}}$, scalar constant values $c_1$ and $c_2$, cardinal $N$ of $\Xi_r$, and AS ROB $\Vbold_{\mu}$.
    \OUTPUT Set of feasible candidate points in the AS, $\Xi_r$, and set of feasible points in the full space $\mathcal D$, $\Xi$.
	  \STATE Sample $N$ points $\{ \mubold_{r_i} \in {\mathcal G} \}_{i=1}^{N}$, using a uniform tensor product design of experiment

	  $$c_1 diag\left( sign( \Vbold_{\mu} )^T \Vbold_{\mu} \right) \leq { \mubold_{r_i} } \leq c_2 diag\left( sign( \Vbold_{\mu} )^T \Vbold_{\mu} \right), \quad \forall i = 1, \ldots, N~\cite{lukaczyk2014active}$$
      or
	  $$c_1 \leq { \mubold_{r_i} } \leq c_2, \quad \forall i = 1, \ldots, N$$
       \FOR{$i=1:N$}
       \STATE Solve for $\mubold_{r_i}$ the feasibility problem
	  \begin{equation}\label{eq:FEAS}
               {\mubold_i}^{\star}  = \hspace{3mm}
        \begin{aligned}
            & \underset{\mubold \, \in \, \Rbb^{N_{{\mu}}}}{\text{min}}
            & & \boldsymbol{0}^T \mubold \\
            & \text{s.t.} & &    \Vbold_{\mu}^T  \mubold =  {\mubold_{r_i}} \\
            & & &  \mubold_{\text{lb}} \leq \mubold \leq
				 \mubold_{\text{ub}}  \\
         \end{aligned}
	  \end{equation}

    \STATE If the feasibility problem has a nontrivial solution, update $\Xi_r$ and $\Xi$
	  $$ \Xi_r = \left\{  \Xi_r \cup \mubold_{r_i} \right\}, \quad \Xi = \left\{  \Xi \cup \mubold_i^{\star} \right\}$$

\ENDFOR
  \end{algorithmic}
\end{algorithm}

As for the residual-based error indicator, which is needed for sampling at each iteration $i$ the parameter point $\mubold_{r_i} = \arg \max_{\mubold_r \, \in \, \Xi_r} e(\mubold_r, \mathcal{DB}_{i-1})$,
a practical and cost-effective choice in the context of the solution of the MDAO problem~\eqref{eq:nandPb1} is
\begin{equation}\label{eq:approxE}
	e(\mubold_r, {\mathcal DB}) = \left\|{\widetilde \Rbold}_L\left(\widetilde{\Vbold}_q (\Vbold_{\mu}\mubold_r)\,\qbold_r (\mubold_r)\right)\right\|_2
\end{equation}
Here, ${\widetilde \Rbold}_L$ is the residual~\eqref{eq:linearPDE} associated with the linearized FSI system~\eqref{eq:ALElinear1}, $\Vbold_q \left(\Vbold_{\mu}\mubold_r\right)\,\qbold_r (\mubold_r)$ is
the reconstructed, high-dimensional, linearized FSI solution at the reconstructed design parameter point $\Vbold_{\mu}\mubold_r$, and the tilde notation designates that this ROB is not available
in the content of the database ${\mathcal DB}_{i-1}$ that is available at the beginning of the $i$-th iteration of the greedy sampling procedure. Since
$\widetilde{\Vbold}_q \left(\Vbold_{\mu}\mubold_{r_i}\right)$ is expected to satisfy the same orthogonality procedure satisfied by every previously sampled ROB
$\Vbold_q \left(\Vbold_{\mu}\mubold_{r_j}\right)$, $j = 1, \ldots, i-1$, it can be rigorously approximated by interpolation on a Grassmann manifold $G_{N_q, n_q}$. Alternatively, it can be approximated
faster using constant extrapolation from the nearest, previously sampled parameter point $\mubold_{r_j}$. In either case, this justifies the storage in $\mathcal DB$ of the pre-computed ROBs noted
in~\eqref{eq:greedy}.

The evaluation of ${\widetilde \Rbold}_L$ at the reconstructed FSI solution $\widetilde{\Vbold}_q \left(\Vbold_{\mu}\mubold_r\right)\,\qbold_r (\mubold_r)$ is typically inexpensive compared to
the solution of the problem $\widetilde{\Rbold}_{L}\left(\textbf{q}(\mubold),\mubold\right) = \mathbfcal{A}(\mubold)\, \dot{\textbf{q}}(\mubold) + \mathbfcal{B}(\mubold)\, \textbf{q}(\mubold) =
\boldsymbol{0}$. Nevertheless, since $\Xi_r$ must be large enough to faithfully represent $\mathcal G$, the evaluation of the error indicator~\eqref{eq:approxE} at a large number of candidate parameter
points may still be overwhelming for some applications. In this case, the error indicator~\eqref{eq:approxE} can be applied only to a subset of $\Xi_r$ that is re-selected at the beginning of each
iteration of the greedy sampling procedure.

\subsection{Interpolation in the AS on a Matrix Manifold}
\label{sec:interpolation}

After a database ${\mathcal DB}$ of consistent tuples of FSI reduced-operators~\eqref{eq:tupleROT} -- which is simply referred to in the remainder of this paper as a consistent database ${\mathcal DB}$ --
is constructed, the linear FSI problem~\eqref{eq:linearPDE} can be solved in real-time at each unsampled parameter point $\mubold_{r_{\star}} \in \mathcal G$ (and therefore $\mubold_{\star}
\in \mathcal D$) queried by the optimization procedure chosen for solving the MDAO problem~\eqref{eq:nandPb1} in two steps as follows:
\begin{itemize}
		\item Interpolate in real-time the content of the consistent database ${\mathcal DB}$ on appropriate matrix manifolds~\cite{amsallem2011online} to compute a linear, FSI
			PROM of the form given in~\eqref{eq:tuple} at the queried parameter point $\mubold_{\star} = \Vbold_{\mu}\mubold_{r_\star}$. For this purpose: note that the content
			of ${\mathcal DB}$, which consists of tuples of FSI reduced-operators of the form given in~\eqref{eq:tupleROT} and ROBs (see~\eqref{eq:greedy}), can be organized block-by-block
			according to~\eqref{eq:MISS2} and~\eqref{eq:subspaceAssumption1}, respectively; and therefore, perform matrix interpolation block-by-block on an appropriate matrix manifold
			$\mathcal M$.
		\item Apply the interpolated tuple ${\mathcal T}_r^{FSI}(\mubold_{\star})$ to solve in real-time the linearized, FSI constraint~\eqref{eq:linearPDE}.
\end{itemize}
Here, the appropriate matrix manifold $\mathcal M$ can be chosen as follows:
\begin{itemize}
	\item The manifold of invertible real matrices of size $n_f$, $GL(n_f, \mathbb R)$, for the submatrices of the form $\Abold_r(\mubold_r) = \Wbold_w^T\Abold(\mubold_r)\Vbold_w \in
		{\mathbb R}^{n_f\times n_f}$ of a matrix of the form ${\mathbfcal A}_r(\mubold_r)$~\eqref{eq:MISS2}, and the submatrices of the form $\Hbold_r(\mubold_r) = \Wbold_w^T\Hbold(\mubold_r)
		\Vbold_w \in {\mathbb R}^{n_f\times n_f}$ of a matrix of the form ${\mathbfcal B}_r(\mubold_r)$~\eqref{eq:MISS2}.
	\item The manifold of symmetric positive definite (SPD) matrices of size $n_s$, $SPD(n_s)$, for the submatrices of the form $\Dbold_r(\mubold_r) = \Vbold_u^T\Dbold(\mubold_r) \Vbold_u \in
		{\mathbb R}^{n_s\times n_s}$ and $\Kbold_r(\mubold_r) = \Vbold_u^T\Kbold(\mubold_r) \Vbold_u = \boldsymbol{\Omega}_r^2(\mubold_r) \in {\mathbb R}^{n_s\times n_s}$ of a matrix of the form
		${\mathbfcal B}_r(\mubold_r)$~\eqref{eq:MISS2}.
	\item The manifold of $n_f\times n_s$ real matrices, $\mathbb{R}^{n_f\times n_s}$, for the submatrices of the form $\Rbold_r(\mubold_r) = \Wbold_w^T\Rbold(\mubold_r)\Vbold_u \in
		{\mathbb R}^{n_f\times n_s}$ and the submatrices of the form $\Gbold_r(\mubold_r) = \Wbold_w^T\Gbold(\mubold_r)\Vbold_u \in {\mathbb R}^{n_f\times n_s}$ of a matrix of the form
		${\mathbfcal B}_r(\mubold_r)$~\eqref{eq:MISS2}.
	\item The manifold of $n_s\times n_f$ real matrices, $\mathbb{R}^{n_s\times n_f}$, for the submatrices of the form $\Pbold_r(\mubold_r) = \Vbold_u^T\Pbold(\mubold_r)\Vbold_w \in
		{\mathbb R}^{n_s\times n_f}$ of a matrix of the form ${\mathbfcal B}_r(\mubold_r)$~\eqref{eq:MISS2}.
\end{itemize}
As for the interpolation procedure, it can be summarized in 3 steps as follows (see~\cite{amsallem2011online}):
\begin{enumerate}
	\item Choose a reference block $\Xbold_{ref} = \Xbold(\mubold_{r_{ref}})$ among the matrix blocks $\Xbold_j = \Xbold(\mubold_{r_j})$, $j = 1, \ldots, N_{\mathcal DB}$, to be interpolated on the
		matrix manifold $\mathcal M$.
	\item Apply the logarithm map to transfer each matrix block $\Xbold_j$ to the tangent space to $\mathcal{M}$ at $\Xbold_{ref}$, $\mathcal{T}_{\Xbold}\mathcal{M}$ -- that is,
		compute the logarithm of each matrix $\Xbold_j$, $\Gammabold_j = \Gammabold(\mubold_{r_j}) = \text{Log}_{\Xbold_{ref}} \Xbold_j$, $j = 1, \ldots, N_{\mathcal DB}$.
	\item In $\mathcal{T}_{\Xbold}\mathcal{M}$, interpolate in the AS the computed matrix logarithms $\Gammabold_j = \Gammabold(\mubold_{r_j})$, $j = 1, \ldots, N_{\mathcal DB}$, entry-by-entry,
		using weighted sums of radial basis functions applied directly to the parameter points $\mubold_{r_j}$, $j = 1, \ldots, N_{\mathcal DB}$~\cite{YCetAL}. Let
		$\Gammabold_{\star} = \Gammabold(\mubold_{r_{\star}})$ denote the result of this interpolation.
	\item Apply the exponential map to $\Gammabold_{\star}$ in order to bring this matrix to the manifold $\mathcal M$ -- that is, compute the exponential of the matrix $\Gammabold_{\star}$,
		$\Xbold_{\star} = \Xbold(\mubold_{r_{\star}}) = \text{Exp}_{\Xbold_{ref}} \Gammabold_{\star}$.

\end{enumerate}

A related interpolation method proposed in~\cite{YCetAL} can be used to interpolate the sensitivity matrices defining the extended tuples~\eqref{eq:exttuple} without having to pre-compute and store
in $\mathcal DB$ the additional quantities $\left\{\displaystyle{\frac{\partial \mathbfcal{A}_r(\mubold_r)}{\partial \mu[j]}}\right\}_{j=1}^{N_{\mathcal D}}$ and
$\left\{\displaystyle{\frac{\partial \mathbfcal{B}_r(\mubold_r)}{\partial \mu[j]}} \right\}_{j=1}^{N_{\mathcal D}}$ characterizing such extended tuples -- which are needed for the fast solution
of the MDAO problem~(\ref{eq:nandPb1}) using a gradient-based optimization method.

Table~\ref{ta:manifoldInterpolation} gives the expressions of the logarithm and exponential maps for the various matrix manifolds identified above. Algorithm~\ref{al:manifoldInterpolation} summarizes
the method of interpolation in the AS on a matrix manifold.

   \begin{table}[!ht]
        \caption{Logarithm and exponential maps for some matrix manifolds of interest.}
        \centering
        \begin{tabular}{|l|c|c|c|}
          \hline\Xhline{2\arrayrulewidth}
		Manifold & $\mathcal{R}^{m\times n}$ & $GL(n, {\mathbb R})$ & $SPD(n)$ \\ [0.5ex]
          \Xhline{2\arrayrulewidth}
          $\text{Log}_{\Xbold}(\Ybold)$      &  $\Ybold-\Xbold$      &  $\text{log}(\Ybold\Xbold^{-1})$  & $\text{log}(\Xbold^{-1/2}\Ybold\Xbold^{-1/2})$ \\ [1ex]
          $\text{Exp}_{\Xbold}(\Gammabold)$  &  $\Xbold+\Gammabold$  &  $\text{exp}(\Gammabold)\Xbold$   & $\Xbold^{1/2}\text{exp}(\Gammabold)\Xbold^{1/2}$ \\ [1ex]
          \hline\Xhline{2\arrayrulewidth}
        \end{tabular}
        \label{ta:manifoldInterpolation}
      \end{table}

 \begin{algorithm}[!ht]
      \caption{Interpolation in the AS on a matrix manifold $\mathcal M$.}
      \label{al:manifoldInterpolation}
	 \textbf{Input:} $N_\mathcal{DB}$ PROM matrices $\Xbold_j = \Xbold(\mubold_{r_j})$ belonging to $\mathcal M$, $\mubold_{r_j} \in \mathcal G$, $j = 1, \ldots, N_{\mathcal DB}$, and queried but unsampled parameter point
	 ${\mubold}_{r_{\star}} \in {\mathcal G}$.\\
	 \textbf{Output:} Interpolated PROM matrix $\Xbold_{\star} = \Xbold(\mubold_{r_{\star}})$.
      \begin{algorithmic}[1]
        \STATE Choose a reference PROM matrix $\Xbold_{ref} = \Xbold(\mubold_{r_{ref}})$
        \FOR {$j=1, \ldots, N_\mathcal{DB}$}
           \STATE Compute $\Gammabold_j = \text{Log}_{\Xbold_{ref}}\Xbold_j$
        \ENDFOR
        \STATE Interpolate in the AS the matrices $\Gammabold_j$, $j = 1, \ldots, N_{\mathcal DB}$, using weighted sums of radial basis functions to obtain $\Gammabold_{\star} = \Gammabold(\mubold_{r_{\star}})$.
        \STATE Compute $\Xbold_{\star} = \Xbold(\mubold_{r_{\star}}) = \text{Exp}_{\Xbold_{ref}} \Gammabold_{\star}$.
      \end{algorithmic}
      \end{algorithm}

\section{Applications}
\label{sec:APP}

In this section, the computational framework described in this paper -- which intertwines the proposed alternative AS concept with adaptive parameter sampling, PMOR, and the concept of a database of
linear PROMs equipped with interpolation on matrix manifolds -- is applied to the solution of MDAO problems with flutter constraints. For this purpose, two different aeroelastic systems are considered:
a flexible configuration of NASA's CRM; and NASA's ARW-2. In both cases, air is modeled as a perfect gas, and the flow is assumed to be inviscid.

All HDM, PMOR, and PROM computations discussed below are performed in double precision floating point arithmetic on a Linux cluster, using AERO
Suite~\cite{farhat2003application, geuzaine2003aeroelastic}. Specifically, all fluid and structural ROBs and all aeroelastic PROMs are computed as described in Section~\ref{sec:PPMOR}.

In~\cite{YCetAL}, the authors considered a computational framework for linear PMOR that is related to that presented in this paper, but without any AS concept. They applied it to MDAO problems similar
to those considered here, but with fewer optimization parameters. They discussed in detail the computational benefits of the framework and reported three orders of magnitude speedup factors due to PMOR.
Because the main difference between the computational framework described in~\cite{YCetAL} and its counterpart presented in this paper is the incorporation in the process of the AS concept presented in
Section~\ref{sec:ASNT} in order to enable computational feasibility for a larger number of optimization parameters, only the speedup factors due to the considered AS concept are reported herein. These
speedup factors can be simply multiplied by those due to PMOR in the absence of an AS to assess the potential of the overall computational framework described in this paper for accelerating the solution
of MDAO problems with linearized FSI constraints.

\subsection{Flutter Constraint Using Linearized CFD-based Aeroelastic PROMs}\label{sec:dampROM}

In aeronautics, flutter constraints are typically formulated in terms of a lower bound on the modal damping ratios of the linear dynamical system~\eqref{eq:PROM_unc}, or more practically here, its
reduced-order counterpart~\eqref{eq:PROM2}. The latter linear, aeroelastic PROM can be re-written as
\begin{equation}\label{eq:FDR}
	\dot{\textbf{q}}_r(\mubold_r)  + \mathbfcal{N}_r (\mubold_r) \, \textbf{q}_r(\mubold_r) = \boldsymbol{0},\qquad \hbox{where} \quad
	\mathbfcal{N}_r(\mubold_r) = \mathbfcal{A}_r^{-1} (\mubold_r) \mathbfcal{B}_r (\tilde{\mubold}_r)
\end{equation}
and $\mathbfcal{A}_r(\mubold_r)$ and $\mathbfcal{B}_r(\mubold_r)$ are given in~\eqref{eq:MISS2}.

The eigenvalue problem associated with~\eqref{eq:FDR} is
\begin{equation}\label{eq:EP}
	\mathbfcal{N}_r (\mubold_r) \, \hat{\qbold}_j({\mubold}_r) = \lambda_j({\mubold}_r) \, \hat{\qbold}_j({\mubold}_r), \qquad j=1, \ldots, n_q
\end{equation}
where $\lambda_j({\mubold}_r) \in \mathbb C$ denotes the $j$-th eigenvector and $\hat{\qbold}_j({\mubold}_r) \in {\mathbb C}^{n_q}$ denotes its associated eigenvector (note that if
an eigenpair $\left(\lambda_j(\mubold_r), {\hat\qbold}_j(\mubold_r) \right)$ is complex-valued, its complex conjugate is also an eigenpair of~\eqref{eq:EP}).

Then, the modal damping ratios of the linear, aeroelastic PROM~\eqref{eq:PROM2} are given by
\begin{equation*}\label{eq:dampRatioDef}
\zeta_{j} =-\frac{\lambda_j^R}{\sqrt{\left(\lambda_j^R\right)^2+\left(\lambda_j^I\right)^2}}, \qquad j=1, \ldots, n_q
\end{equation*}
where $\lambda_j^R$ and $\lambda_j^I$ denote the real and imaginary parts of the complex eigenvalue $\lambda_j$, respectively, and the flutter constraints are typically formulated as
\begin{equation*}
	\zeta_{j}(\mubold_r) \ge \zeta_{\text{lb}}, \qquad \qquad j=1, \ldots, n_q
\end{equation*}
where $\zeta_{\text{lb}}$ is a regulation-specified lower bound.

\subsection{Aeroelastic Tailoring of a Flexible CRM Configuration}
\label{sec:CRM}

First, the following sixth-dimensional MDAO problem for a flexible configuration of NASA's CRM, whose geometrical and material descriptions are summarized in Table~\ref{tab:GMP}, is considered
\begin{equation}\label{eq:optimizationCRM}
    \begin{aligned}
      & \underset{\mubold\, \in \, \mathcal{D} \, \subset \, \mathbb{R}^{6}}{\text{maximize}}  & & \frac{L(\mubold)}{D(\mubold)}
          \\ & \text{subject to} & & \sigma_{\text{VM}}(\mubold) \leq \sigma_{\text{ub}}
                    \\ & & & \mubold_{\text{lb}} \leq \mubold \leq \mubold_{\text{ub}}
                    \\ & & & \zetabold(\mubold) \geq \textbf{1} \zeta_{\text{lb}}
    \end{aligned}
\end{equation}
where:
\begin{itemize}
	\item $L(\mubold)$ and $D(\mubold)$ denote the parametric lift and drag, respectively, associated with the following flight conditions:
		\begin{itemize}
			\item Level flight at the altitude of $h = 10,668$ m, where the free-stream pressure and density are $P_\infty = 23835.89\hspace{3pt} \text{Pa}$ and
				$\rho_\infty = 0.3796275\hspace{3pt} \text{kg}/ \text{m}^3$, respectively.
			\item Free-stream Mach number $M_{\infty} = 0.85$, angle of attack AoA $= 1.823^{\circ}$, and angle of sideslip AoS $= 0^{\circ}$.
		\end{itemize}
	\item $\sigma_{\text{VM}}(\mubold)$ denotes the parametric von Mises stress at any point of the structural model of the CRM, and its upper bound is set to $\sigma_{\text{ub}}= 1.331\times10^9$ Pa.
	\item $\mubold_{\text{lb}}$ and $\mubold_{\text{ub}}$ are lower and upper bounds for the vector of optimization parameters, respectively, and define box constraints
		for $\mubold$.
	\item $\zetabold(\mubold)$ is the parametric vector of damping ratios $\zeta_j$, $\textbf{1} \zeta_{\text{lb}}  > \zerobold$ is its regulation-specified lower bound,
		$\textbf{1} \in \mathbb{R}^{6}$ and $\zerobold  \in \mathbb{R}^{6}$ are vectors of ones and zeros, respectively, the lower bound coefficient $\zeta_{\text{lb}}$
		is set to $\zeta_{\text{lb}} = 4.75 \times 10^{-3}$, and $\zetabold(\mubold) \geq \textbf{1} \zeta_{\text{lb}}$ defines a flutter avoidance constraint.
\end{itemize}

\begin{table}[!ht]
	\caption{Geometrical and material descriptions of a flexible CRM configuration.}
        \centering
	\begin{tabular}{|l|l|l|}
         \hline\Xhline{2\arrayrulewidth}
		                           & Item                             &  Value                 \\ [0.5ex]
          \hline\Xhline{2\arrayrulewidth}
		Geometry                   &                                  &                        \\ [1ex]
		\hline

		\hspace{5pt} Wing          & span                             &   58.76 m              \\
					   & root chord                       &   11.92 m              \\
		                           & tip chord                        &    2.736 m             \\ [1ex]
          \hline\Xhline{2\arrayrulewidth}
		Materials                  &                                  &                        \\ [1ex]
		\hline

		\hspace{5pt}Skin, except flaps &                              &                        \\
		\hspace{5pt}(various composite materials) &                   &                        \\[2ex]
		\hspace{5pt}Stiffeners (aluminum) & E                         & $73.1\times10^9$ Pa  \\
					   & $\rho$                           & $2.78\times 10^{3}$ $\text{kg}/\text{m}^3$    \\
					   & $\nu$                            & 0.3                 \\
          \Xhline{2\arrayrulewidth}
        \end{tabular}
	\label{tab:GMP}
\end{table}

Because the focus is set here on level flight at a zero angle of sideslip, the considered CRM configuration is a symmetric one. Hence, only half of its geometry is modeled.

\begin{figure}[!ht]
\begin{center}
\includegraphics[scale=0.45]{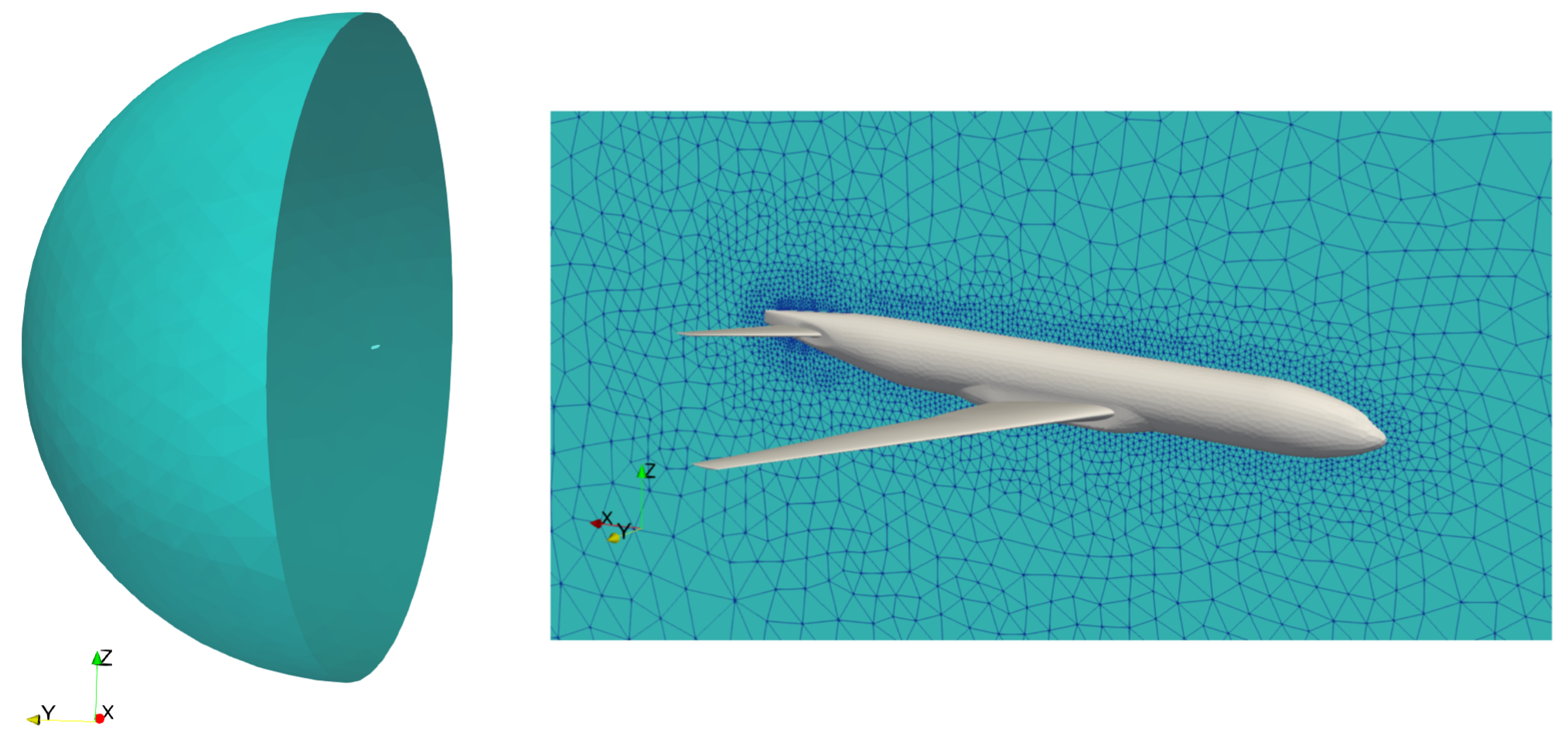}
\end{center}
\caption{CRM: hemispherical computational fluid domain with a symmetry plane (left); and partial view of its discretization by an inviscid, unstructured, body-fitted mesh (right).}
\label{fig:CRM-CFD}
\end{figure}

The computational fluid domain is chosen to be a hemisphere with a symmetry plane along the middle of the fuselage. It is discretized by a three-dimensional (3D), unstructured, body-fitted CFD mesh
with $740,248$ grid points (see Figure~\ref{fig:CRM-CFD}). The parametric lift and drag are obtained by postprocessing the computed flow solution.

\begin{figure}[!ht]
	\begin{center}
		\includegraphics[scale=0.40]{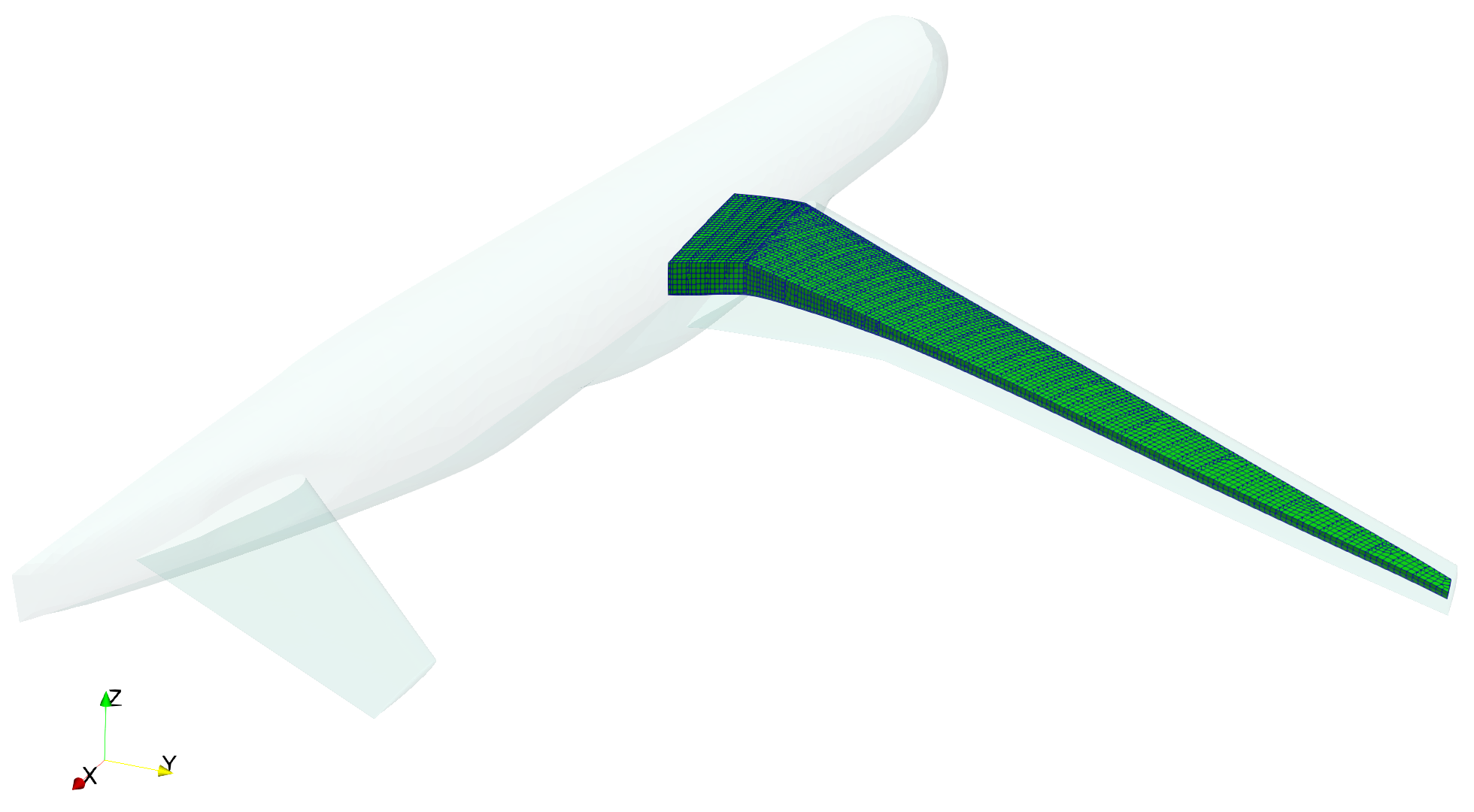}
		\includegraphics[scale=0.40]{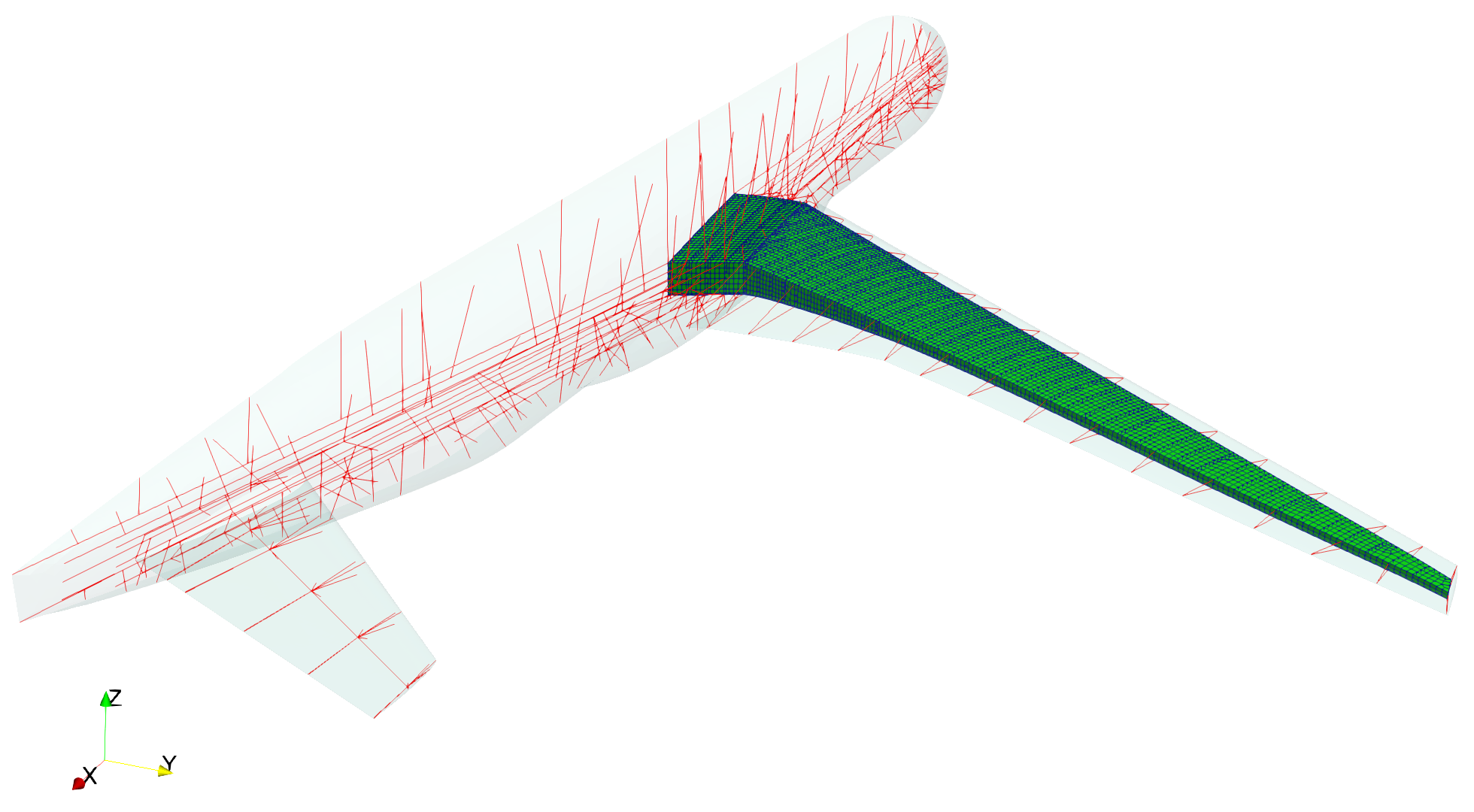}
		\includegraphics[scale=0.40]{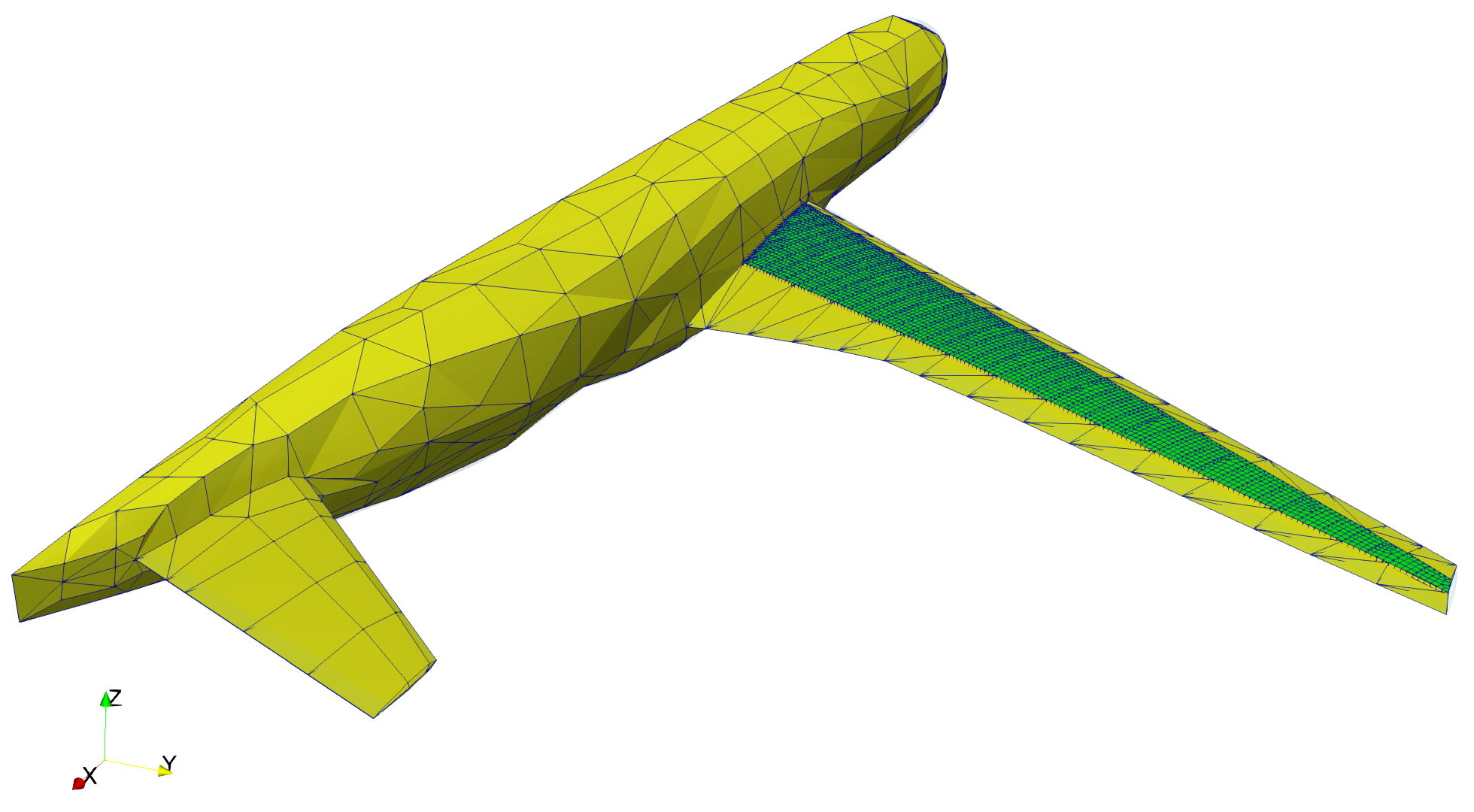}
	\end{center}
	\caption{CRM -- FE structural respresentations of: wing (top); fuselage and control surfaces (middle); and skin (bottom).}
	\label{fig:CRM-CSD}
\end{figure}

The parametric von Mises stress field is computed using the FE structural representation of half of the flexible CRM configuration shown in Figure~\ref{fig:CRM-CSD}, which has 58,888 degrees of freedom
(dofs). In this representation, the FE structural modeling of the wing is that described in~\cite{Kenway2014c}. The other components of the FE structural representation shown in Figure~\ref{fig:CRM-CSD}
were developed by the authors and incorporated in the final FE structural model to facilitate two-way coupled aeroelastic computations.

The flutter avoidance constraint is the most CPU intensive constraint in the formulation of the MDAO problem~\eqref{eq:optimizationCRM}. It is approximated using a database of linear FSI PROMs.

The six design optimization parameters stored in $\mubold$ are organized in two groups: three optimization parameters $\{\mubold[1], \mubold[2], \mubold[3] \}$ that are shape design parameters,
and three optimization parameters $\{ \mubold[4], \mubold[5], \mubold[6] \}$ that are structural design parameters.

\begin{figure}[!ht]
      \begin{center}
        \includegraphics[width = 3.25truein]{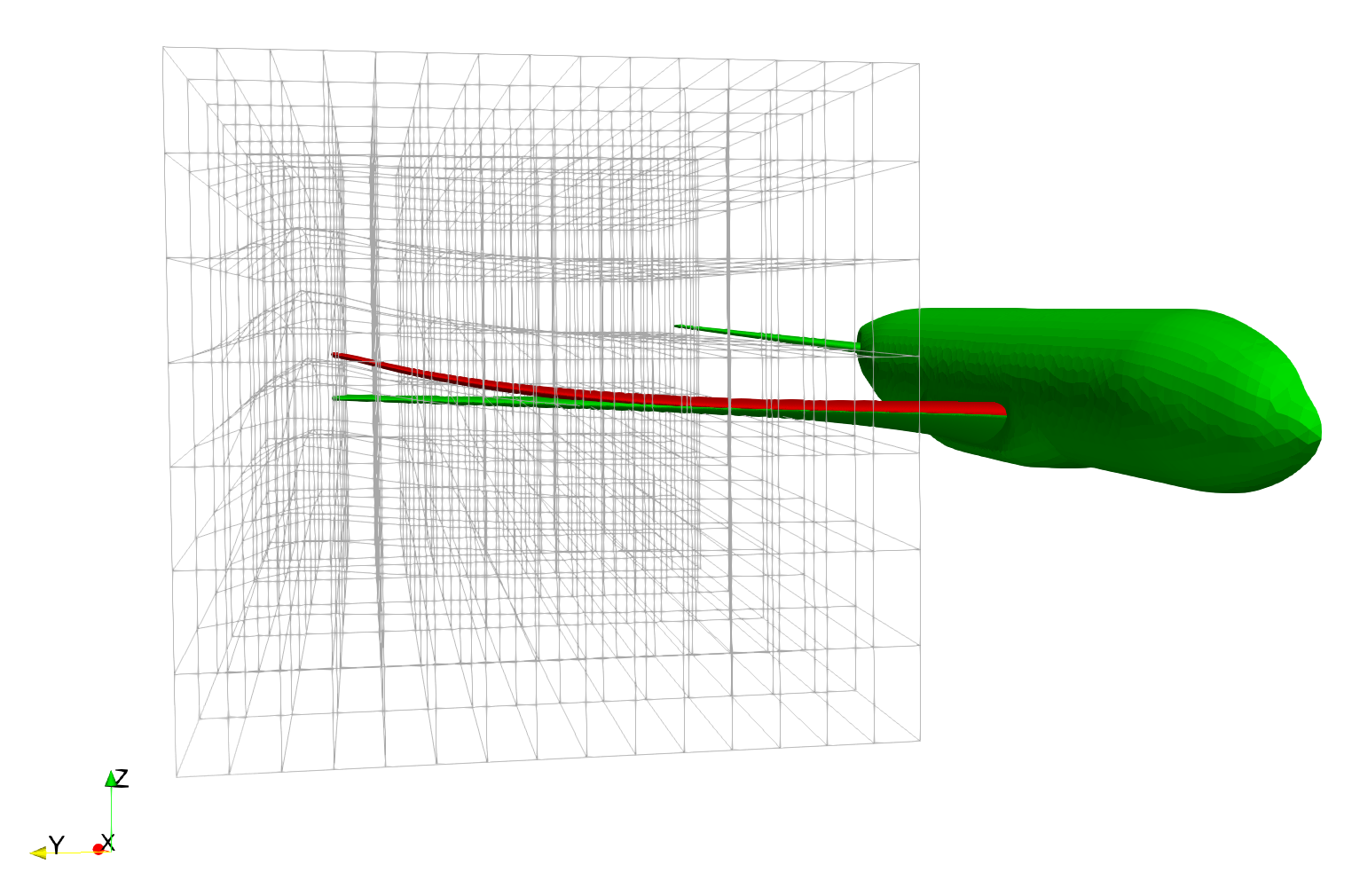}
         \includegraphics[width = 3.25truein]{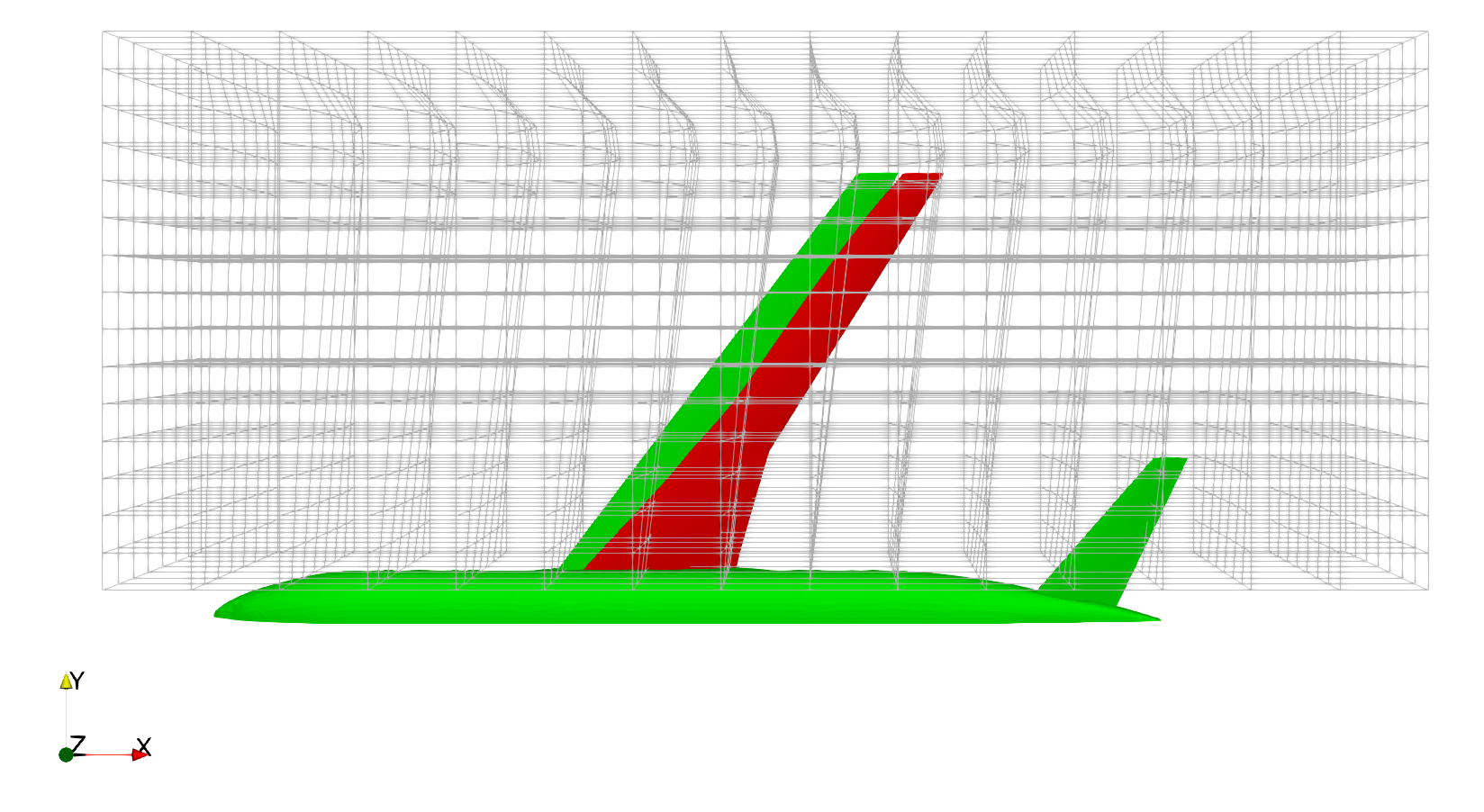}
      \end{center}
      \caption{CRM -- parameters controling: dihedral angle of the wing (left); and sweep angle of the wing (right).}
      \label{fig:blenderCRM}
\end{figure}

\begin{figure}[!ht]
      \begin{center}
      \hspace{-10mm}
        \includegraphics[width = 3.30truein]{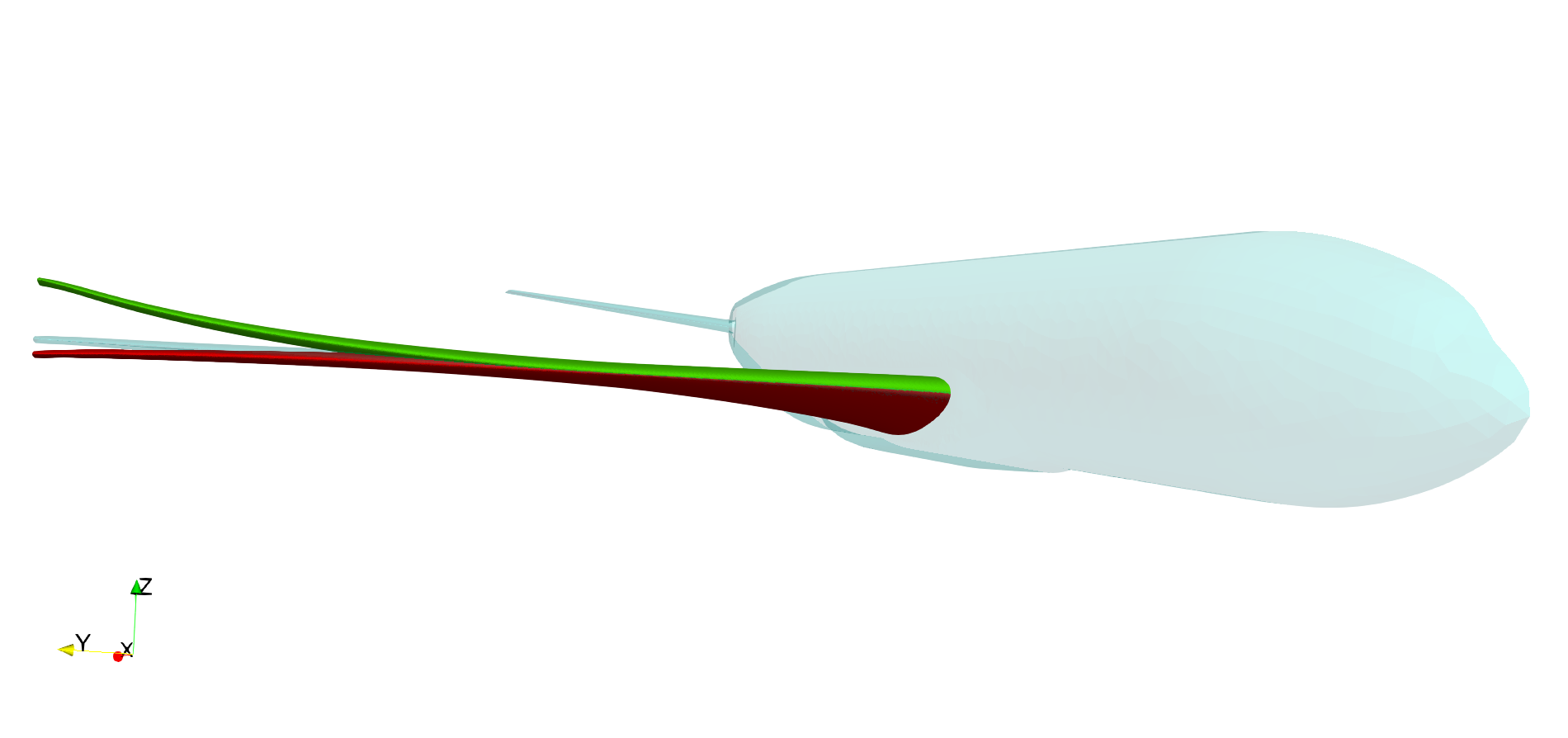}
        \includegraphics[width = 3.30truein]{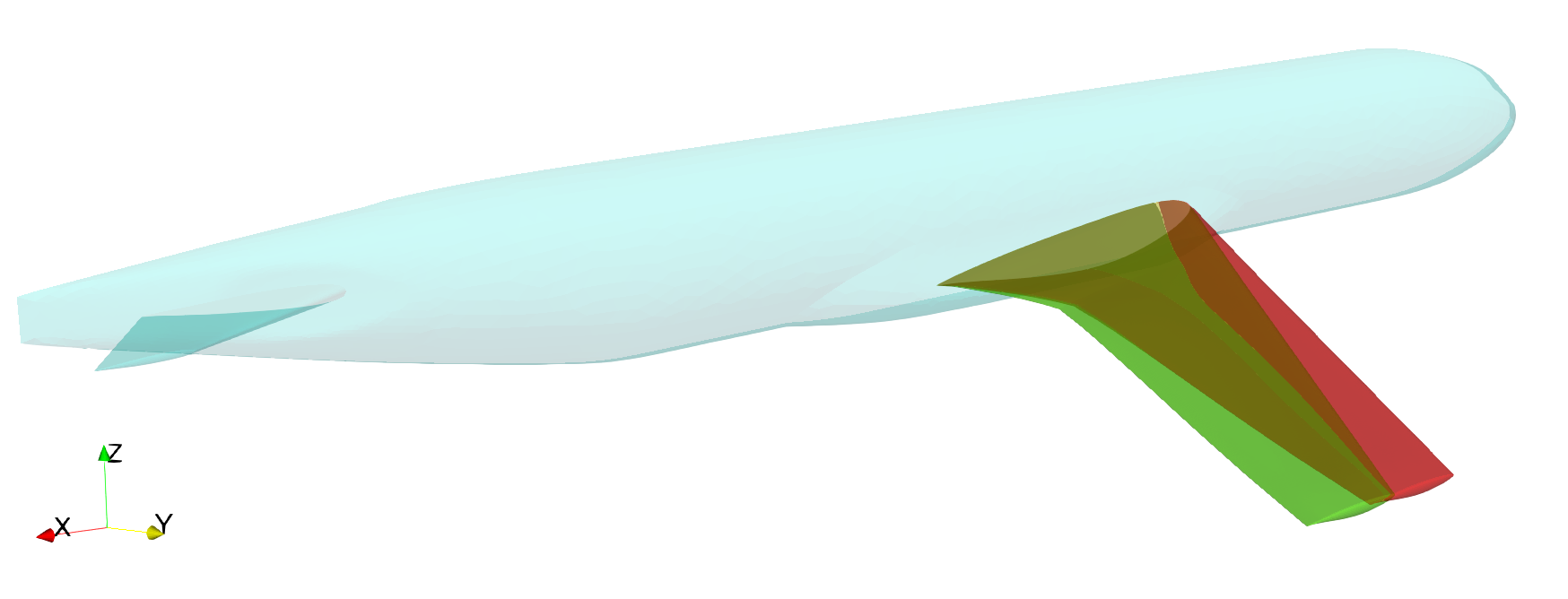}
      \end{center}
	\caption{CRM -- wing shapes associated with: upper (green) and lower (red) bounds for the dihedral angle (left); and upper (green) and lower (red) bounds for the sweep angle (right).}
      \label{fig:boundsCRM}
\end{figure}

The three shape design parameters are chosen as follows: $\mubold[1]$ and $\mubold[2]$ to control the dihedral angle of the wing; and $\mubold[3]$ to control its sweep angle.
All shape changes dictated by the optimization procedures are effected using the open-source 3D computer graphics software Blender~\cite{Blender}. Figure \ref{fig:blenderCRM} illustrates the concept
of lattice-based deformations of Blender to effect shape changes~\cite{anderson2012parametric}, and Figure~\ref{fig:boundsCRM} graphically depicts the effects on the shape of the wing of the
specified lower and upper bounds for the dihedral and sweep angles.

\begin{figure}[!ht]
      \begin{center}
        \includegraphics[width = 3.30truein]{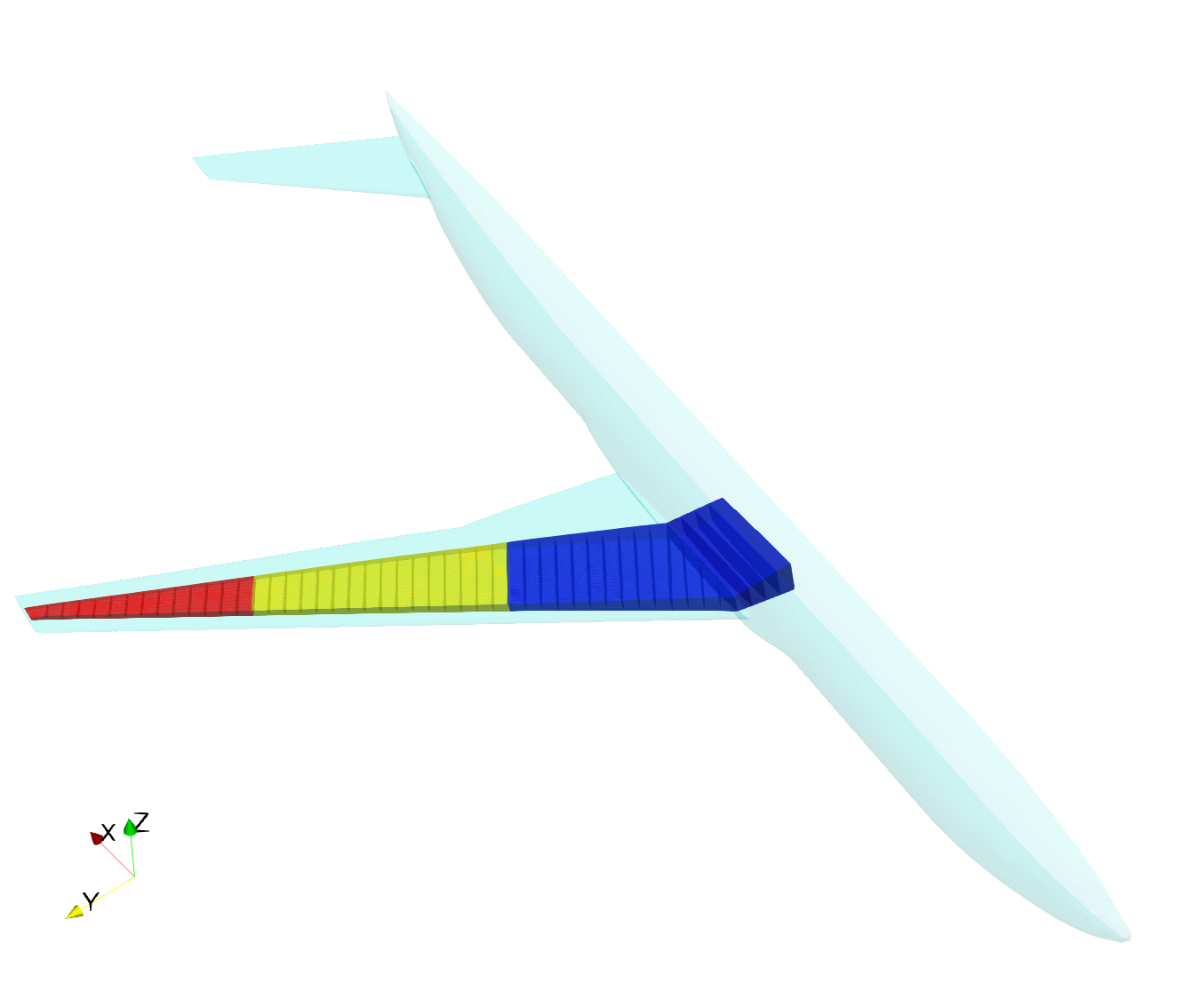}
        \includegraphics[width = 3.30truein]{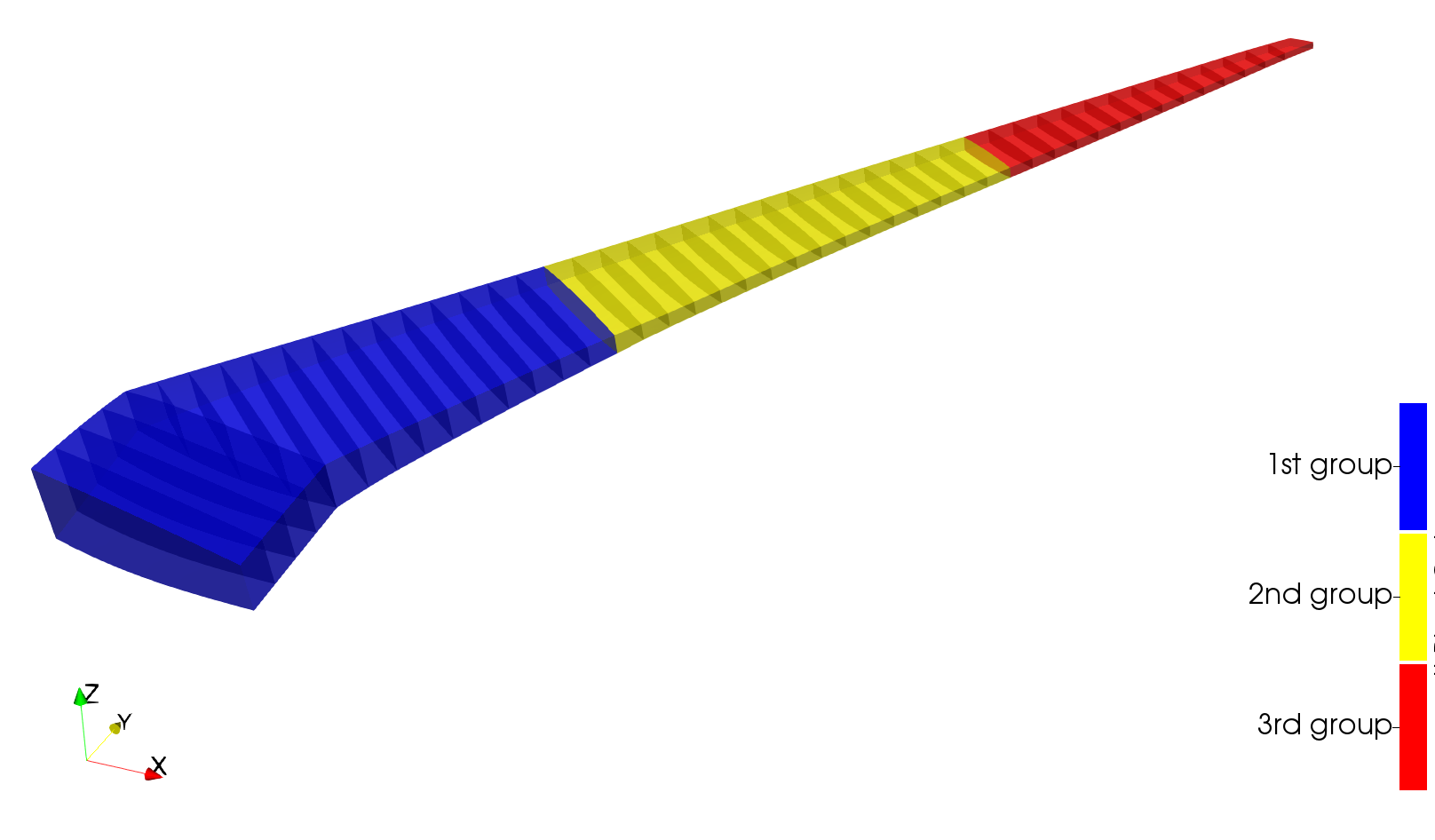}
      \end{center}
	\caption{CRM -- organization of the stiffeners of the wing in three different groups represented by three different colors.}
      \label{fig:CRMstiffenersGroup}
\end{figure}

The three structural design parameters are chosen as three thickness increments for three different groups of stiffeners (see Figure \ref{fig:CRMstiffenersGroup}).
Each increment is defined as a percentage of the initial thickness to which it is applied, and is constrained to have a magnitude less than 10\% (box constraint).

\clearpage
\subsubsection{Offline Construction of Databases of FSI PROMs}
\label{sec:OC}

Three databases of FSI PROMs are built for the solution of the MDAO problem~\eqref{eq:optimizationCRM}:
\begin{itemize}
\item ${\mathcal DB}_1$, which is constructed by directly sampling the design parameter space $\mathcal D$ using a greedy procedure similar to that described in Section~\ref{sec:FEASAM}, but applied
	directly to $\mathcal D$. The greedy procedure samples in this case $N_{\mathcal DB} = 119$ feasible design parameter points.
\item ${\mathcal DB}_2$, which is constructed by building first the lower-dimensional space of design parameter generalized coordinates ${\mathcal G}^{\,\text{cl}}$ using the classical AS method
	described in Section~\ref{sec:BACK} and its associated empirical formula~\eqref{eq:numberGRAD}, then indirectly sampling $\mathcal D$ by directly sampling $\mathcal G$. In this case, the free
	parameters of formula~\eqref{eq:numberGRAD} are chosen as $\alpha = 10$ and $\beta = 6$, which yields $N_{\Sbold} = 47$. Hence, 47 gradients of the objective function
	are computed at 47 feasible parameter points that are randomly sampled in $\mathcal D$ using the LHS method. The computation of each gradient consumes 0.3 hour wall-clock time
	on a Linux cluster with 240 processors. Next, the computed snapshots are compressed into a ROB $\Vbold_{\mu}$ of dimension $n_{\mathcal G}^{\text{cl}} = 3$ (see below).
\item ${\mathcal DB}_3$, which is constructed by building first ${\mathcal G}^{\,\text{al}}$ using the alternative AS method described in Section~\ref{sec:ASNT}, then indirectly sampling $\mathcal D$ by
	directly sampling $\mathcal G$, using the greedy procedure described in Section~\ref{sec:FEASAM}. In this case, the snapshots of the form $\Delta \mubold$ are adaptively computed by solving the
	auxiliary optimization problem
			\begin{equation}\label{eq:optimizationCRM1}
				\begin{aligned} & \underset{\mubold \, \in \, \mathcal{D} \, \subset \, \mathbb{R}^{6}}{\text{maximize}}  & & \frac{L(\mubold)}{D(\mubold)}  \\
					& \text{subject to} & & \sigma_{\text{VM}}(\mubold) \leq \sigma_{\text{ub}} \\ & & & \mubold_{\text{lb}} \leq \mubold \leq \mubold_{\text{ub}}
				\end{aligned}
			\end{equation}
which is obtained by removing from the MDAO problem~\eqref{eq:optimizationCRM} the CPU intensive FSI constraint. This auxiliary
problem is solved in 16 iterations using the sequential least-squares programming method SLSQP, which uses the Han-Powell
quasi-Newton method with a BFGS update and an $L_1$ test function in the step-length algorithm. Hence, the solution of the
optimization problem~\eqref{eq:optimizationCRM1} generates 17 snapshots. The computation of each of the snapshots consumes about
0.34 hour wall-clock time -- roughly the same amount of time as the computation of a gradient snapshot in the classical
AS method -- on the same Linux cluster. Next, the snapshots are compressed into a ROB $\Vbold_{\mu}$ of dimension
$n_{\mathcal G}^{\text{al}} = 3$ (see below).
\end{itemize}

Figure~\ref{fig:svdCRM} plots the distributions of the singular values of both matrices $\Sbold$~\eqref{eq:snapX} and
$\mathbb S$~\eqref{eq:snapX1}. For the snapshot matrix $\Sbold$ associated with the classical AS method, there is no clear cutoff
singular value for constructing a ROB $\Vbold_{\mu}$, as the singular values are shown to continuously decay. For the snapshot matrix
$\mathbb S$ associated with the alternative AS method however, Figure~\ref{fig:svdCRM}-right shows that the last three singular
values are much smaller than the first three ones. This suggests constructing a ROB $\Vbold_{\mu}$ as the set of the first three
columns of the matrix $\Ubold$ arising from the SVD of $\mathbb S$ -- that is, constructing a ROB $\Vbold_{\mu}$ of dimension
$n_{\mathcal G}^{\text{al}} = 3$. Hence, the dimension of $\Sbold$ is also set to $n_{\mathcal G}^{\text{cl}} = 3$, in order to
enable various meaningful performance comparisons.

\begin{figure}[!ht]
\begin{center}
\includegraphics[scale=0.325]{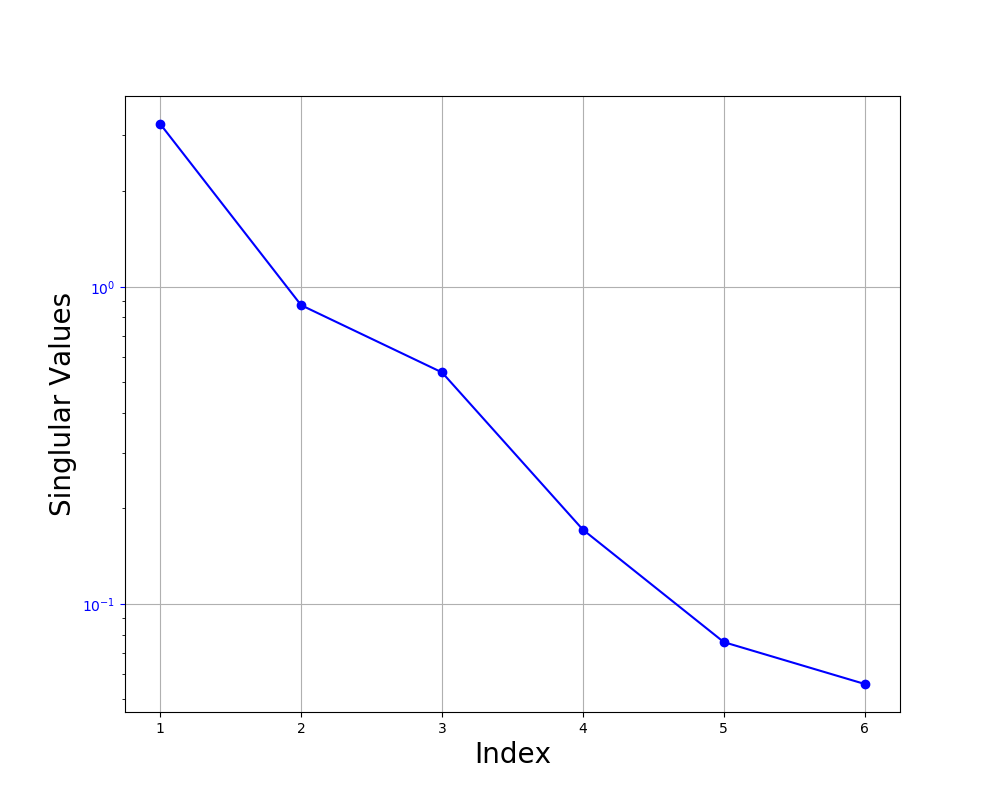}
\includegraphics[scale=0.325]{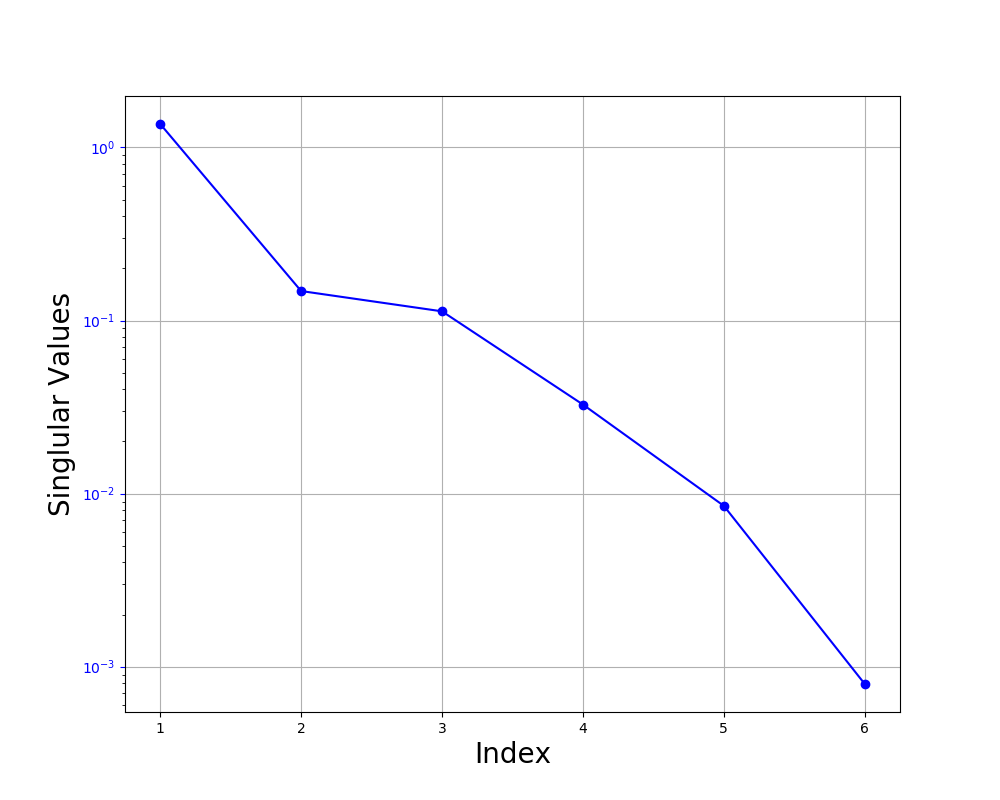}
\end{center}
\caption{Singular values of: the snapshot matrix $\Sbold$ computed using the classical AS method (left); and the counterpart matrix $\mathbb S$ computed using the alternative AS method (right).}
\label{fig:svdCRM}
\end{figure}

\begin{figure}[!ht]
\begin{center}
\includegraphics[scale=0.41]{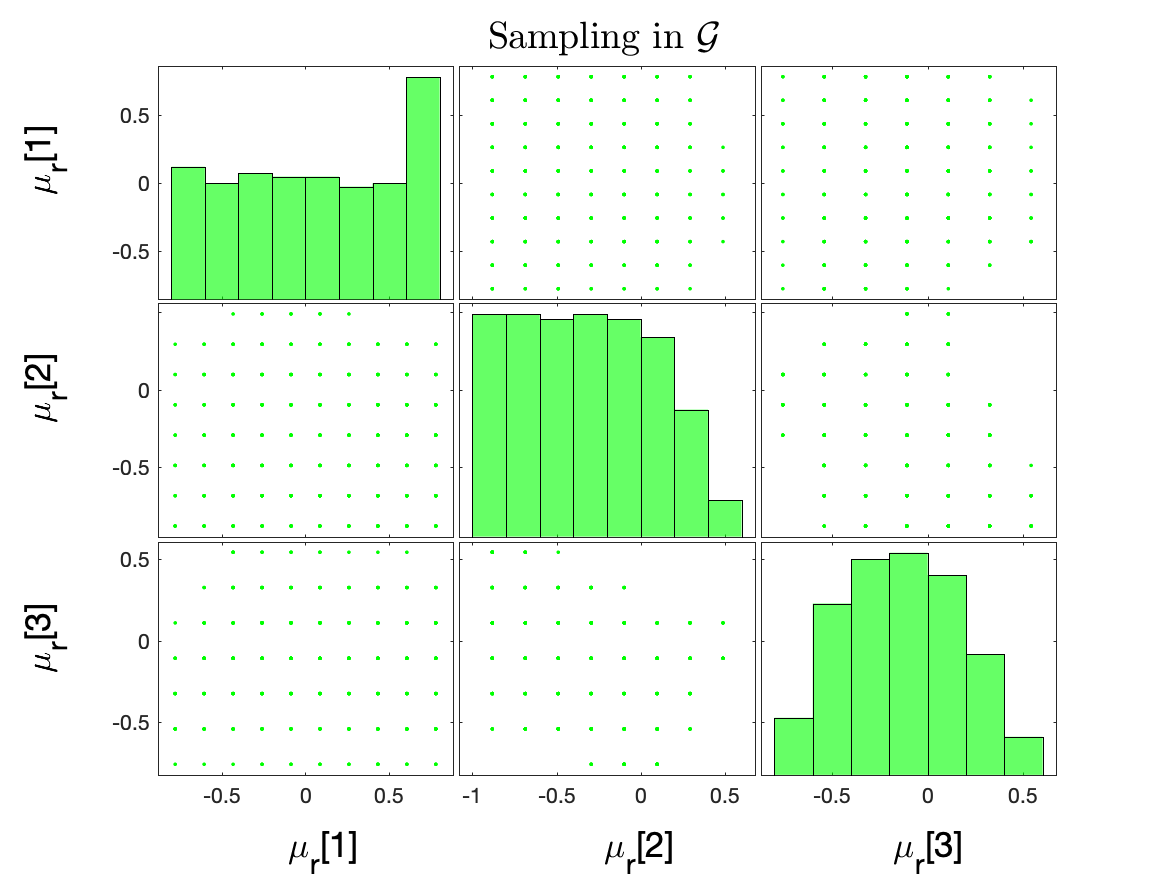}
\includegraphics[scale=0.41]{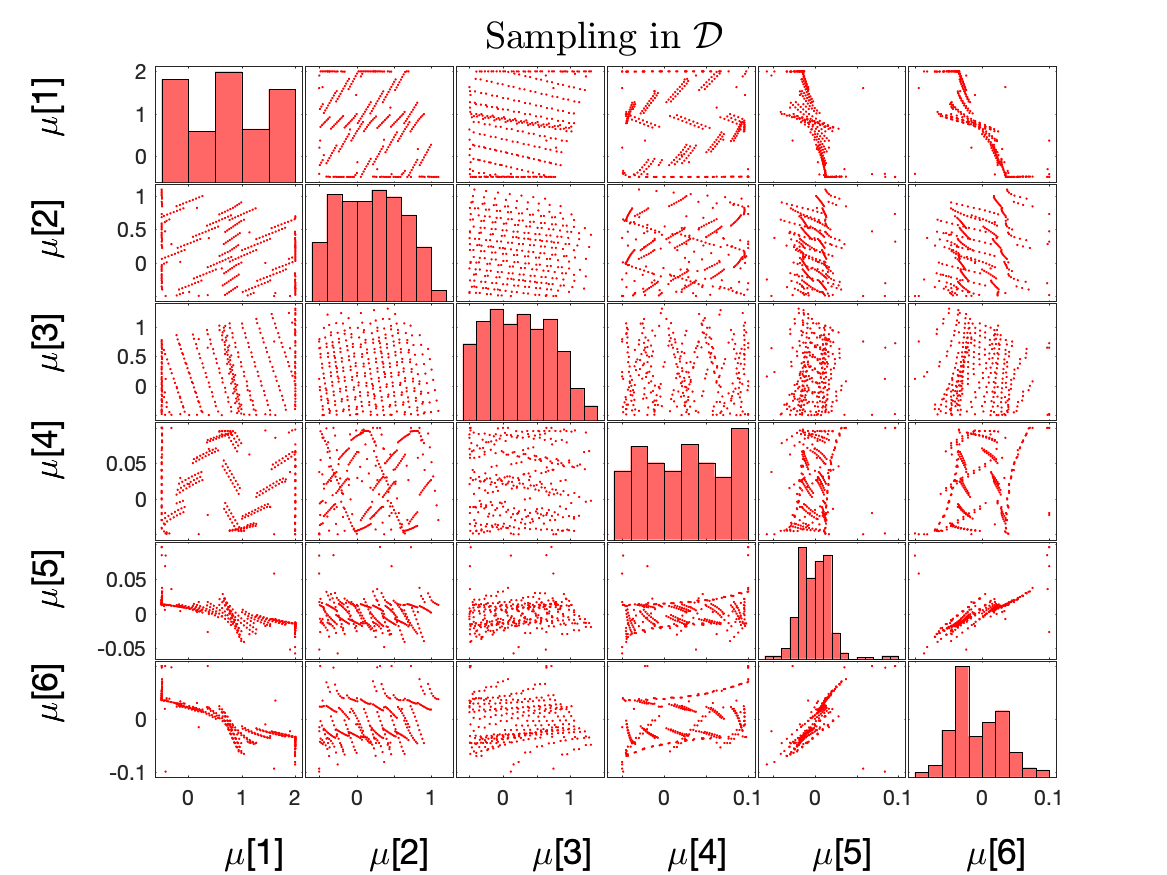}
\end{center}
\caption{Set $\Xi_r \in {\mathcal G}^{\,\text{al}}$ containing 323 feasible candidate points determined by the feasible AS algorithm as well as the alternative approach for constructing an AS (left), and
associated set $\Xi \in \mathcal D$ containing the same number of feasible candidate points (right).}
\label{fig:samplingCRM}
\end{figure}

The concept of indirectly sampling $\mathcal D$ by directly sampling the lower-dimensional space $\mathcal G$ using the feasible AS algorithm (Algorithm~\ref{alg:feasibleAS}) is illustrated in
Figure~\ref{fig:samplingCRM}, for the case of the problem considered herein and ${\mathcal DB}_3$. First, a set of 1,000 preliminary points is generated in ${\mathcal G}^{\,\text{al}}$ using
a uniform tensor product in this parameter space. Next, Algorithm~\ref{alg:feasibleAS} is applied to select among these points a set $\Xi_r$ of 323 feasible candidate points
(see Figure~\ref{fig:samplingCRM}-left), and determine the associated set of feasible points $\Xi \in \mathcal D$ (see Figure~\ref{fig:samplingCRM}--right). Then, the greedy algorithm described
in Section~\ref{sec:FEASAM} is directly applied to $\Xi_r$ -- and indirectly to $\Xi$ -- to perform the final parameter sampling and adaptively construct the database of FSI PROMs ${\mathcal DB}_3$.
Note that the feasible parameter points in $\Xi$ are not uniformly spread in $\mathcal D$, but only in some regions of this design parameter space that depend on the ROB $\Vbold_{\mu}$
(see Figure~\ref{fig:samplingCRM}). The plot in the $i$-th diagonal of Figure~\ref{fig:samplingCRM}-left shows the distribution of the $i$-th component $\mubold_r[i]$ of the candidate parameter points
in $\Xi_r$. Its counterpart in Figure~\ref{fig:samplingCRM}-right shows the distribution of the $i$-th component $\mubold[i]$ of the corresponding candidate parameter points in $\Xi$. The plot in the
$(i,j)$ off-diagonal of Figure~\ref{fig:samplingCRM}-right shows the distribution of the contribution of the $j$-th generalized coordinate $\mubold_r[j]$ to the $i$-th component $\mubold[i]$ of the
candidate parameter point $\mubold$ in $\Xi$.

In each constructed database, each pre-computed linear FSI PROM is built as described in Section~\ref{sec:PPMOR}, with $n_f = 100$ and $n_s = 6$. Hence, each pre-computed linear FSI PROM has
the dimension $n_q = n_f + 2n_s = 112$.

Table (\ref{ta:databaseCRM}) reports on the performance of the offline phase for problem~\eqref{eq:optimizationCRM} of the
computational framework for MDAO described in this paper, depending on whether an AS is used or not, and which method is chosen to
construct it. The following observations are noteworthy:
\begin{itemize}
	\item Even when using the same greedy procedure equipped with the same convergence tolerance, the database of linear FSI
		PROMs ${\mathcal DB}_3$ ends up being less populated than ${\mathcal DB}_1$. This is because in this case,
		${\mathcal DB}_3$ is constructed to be as accurate as ${\mathcal DB}_1$, but in a smaller parameter space.
		Similarly, ${\mathcal DB}_2$ is constructed to be as accurate as ${\mathcal DB}_1$ -- albeit using a different
		construction procedure -- but in a smaller parameter space.
	\item The computational cost of the greedy procedure is roughly proportional to the number of sampled parameter points,
		which itself appears to be proportional to the dimension of the sampled parameter space.
	\item Furthermore, even after accounting for the computational overhead associated with constructing a representative basis
		$\Vbold_{\mu}$, an AS speeds up the construction of a database of PROMs. It remains to asses next however, the
		performance of each database constructed using an AS relative to the performance of ${\mathcal DB}_1$, which is built
		without an AS.
\end{itemize}

\begin{table}[!ht]
\caption{CRM: performance results for the offline phase}
\centering
\begin{tabular}{|l |r|r|r|r|r|r|r|}
\hline\Xhline{2\arrayrulewidth}
Method & $N_{\mathcal D}$ & $N_{\Sbold}$ or $N_{\mathbb S}$ & Wall-clock time & $N_{\mathcal DB}$ & Wall-clock time & Wall-clock time & Speed-up\\[0.5ex]
	&                  &                       & computing $\Vbold_{\mu}$ &                  & greedy procedure & total offline   & relative to no AS\\[0.5ex]
\hline\Xhline{2\arrayrulewidth}
w/o AS         &    6     & N.A. &  N.A.     & 119 & 120.87 hrs & 120.87 hrs & N.A. \\ [1ex]
Classical AS   &    3     & 47   & 14.06 hrs &  61 &  61.95 hrs &  76.01 hrs & 1.59 \\ [1ex]
Alternative AS &    3     & 17   &  5.83 hrs &  62 &  62.97 hrs &  68.80 hrs & 1.76 \\ [1ex]
\hline\Xhline{2\arrayrulewidth}
\end{tabular}
\label{ta:databaseCRM}
\end{table}

\subsubsection{Online Prediction and Aeroelastic Optimization}

\begin{figure}[!ht]
      \begin{center}
        \includegraphics[scale=0.70]{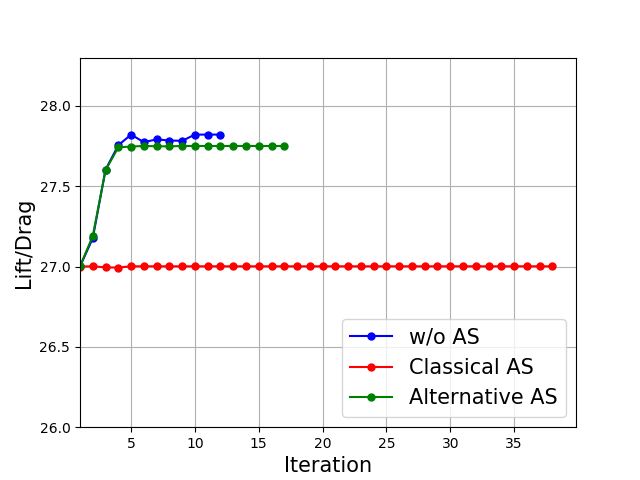}
      \end{center}
      \caption{CRM: convergence histories of the objective function.}
      \label{fig:opta}
\end{figure}

\begin{figure}[!ht]
      \begin{center}
        \includegraphics[scale=0.54]{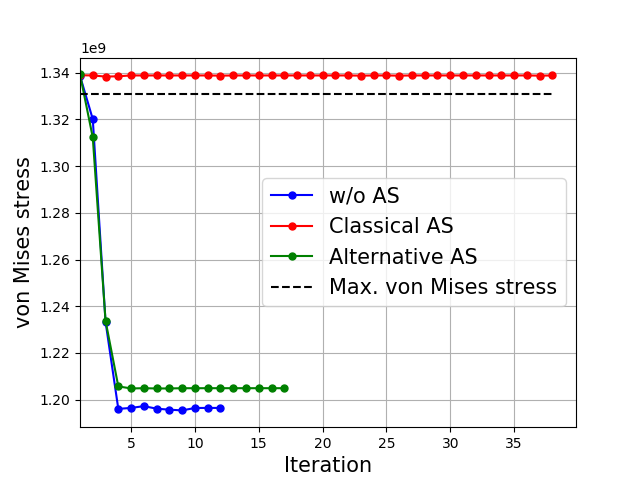}
        \includegraphics[scale=0.54]{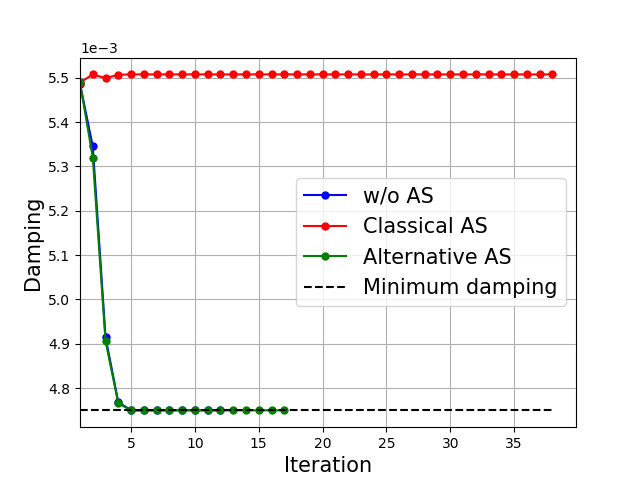}
      \end{center}
      \caption{CRM: convergence histories of the constraints.}
      \label{fig:optb}
\end{figure}

Figure~\ref{fig:opta} and Figure~\ref{fig:optb} display the convergence histories of the objective function and constraints of the
MDAO problem~\eqref{eq:optimizationCRM}, respectively. They highlight that, at least for this CRM MDAO problem, the classical
AS method leads to an ineffective optimization strategy. Indeed, both figures reveal that in the case of the classical AS method,
the PROM-based computational framework described in this paper fails to improve the objective function beyond its initial value.
On the other hand, the same computational framework equipped with the alternative AS method leads to almost the same optimal
objective function and constraints as its counterpart framework equipped with direct sampling of the entire design parameter space, and
in almost the same number of iterations. In particular, note that during the first few iterations, the alternative AS method
leads to the same values of the optimized objective function and constraints as directly sampling the full space ${\mathcal D}$.
This is because during the first few iterations, the flutter constraint is inactive, the auxiliary MDAO
problem~\eqref{eq:optimizationCRM1} is in this case identical to the original MDAO problem~\eqref{eq:optimizationCRM}, and thereore
the basis $\Vbold_{\mu}$ based on the snapshot matrix $\mathbb S$~\eqref{eq:snapX1} is optimal.

Table~\ref{ta:parCRM}) shows the optimal solution of the MDAO problem~\eqref{eq:optimizationCRM} obtained using the three different
PROM databases ${\mathcal DB}_1$, ${\mathcal DB}_2$, and ${\mathcal DB}_3$. The reader can observe that the optimal solution
obtained using the PROM database ${\mathcal DB}_2$ -- that is, the classical AS method for sampling -- is far from its counterpart
obtained using the PROM database ${\mathcal DB}_1$ -- that is, by directly sampling the design parameter space $\mathcal D$. On the
other hand, the local optimal solution obtained using the PROM database ${\mathcal DB}_3$ -- that is, the aternative AS method for
sampling -- is relatively close to its counterpart obtained using the database of linear FSI PROMs ${\mathcal DB}_1$.

\begin{table}[!ht]
\caption{CRM: locally optimal computed solutions}
\centering
\begin{tabular}{| l | r | r | r | r | r | r | r | }
\hline\Xhline{2\arrayrulewidth}
                      & $\mubold[1]$ & $\mubold[2]$ & $\mubold[3]$ & $\mubold[4]$ & $\mubold[5]$ & $\mubold[6]$ \\ [1ex]
\hline\Xhline{2\arrayrulewidth}
Initial configuration &  0.0 &   0.0 &    0.0 &  0.0  &   0.0  &  0.0      \\ [1ex]
\hline
Specified lower bound & -0.5 &  -0.5 &   -0.5 & -0.05 &  -0.1  & -0.1     \\ [1ex]
Specified upper bound &  2.0 &   2.0 &    2.0 &  0.1  &   0.1  &  0.1     \\ [1ex]
\hline\Xhline{2\arrayrulewidth}
Optimal configuration using ${\mathcal DB}_1$&    -0.50 &  -0.16 &  1.47 &   0.07 &   0.05 &   0.06          \\[1ex]
Optimal configuration using ${\mathcal DB}_2$   & $1.6\times10^{-4}$  & $-1.1\times10^{-5}$ & $-2.3\times10^{-3}$ & $-4.9\times10^{-2}$  & $-3.0\times10^{-3}$ & $5.0\times10^{-3}$ \\[1ex]
Optimal configuration using ${\mathcal DB}_3$ &              -0.32 & -0.15              & 1.35               & 0.10                & 0.05               & 0.08 \\
\hline\Xhline{2\arrayrulewidth}
\end{tabular}
\label{ta:parCRM}
\end{table}

In summary, for this MDAO problem with six design optimization parmeters only, the alternative subspace method speeds up the offline
computations by almost a factor two, while leading to a reasonably close local optimum solution.

Next, another MDAO problem with a larger number of design optimization parameters is considered.

\subsection{Aeroelastic Design Optimization of a Parametric ARW-2 Configuration}

\begin{figure}[!ht]
	\begin{center}
      \includegraphics[scale=0.50]{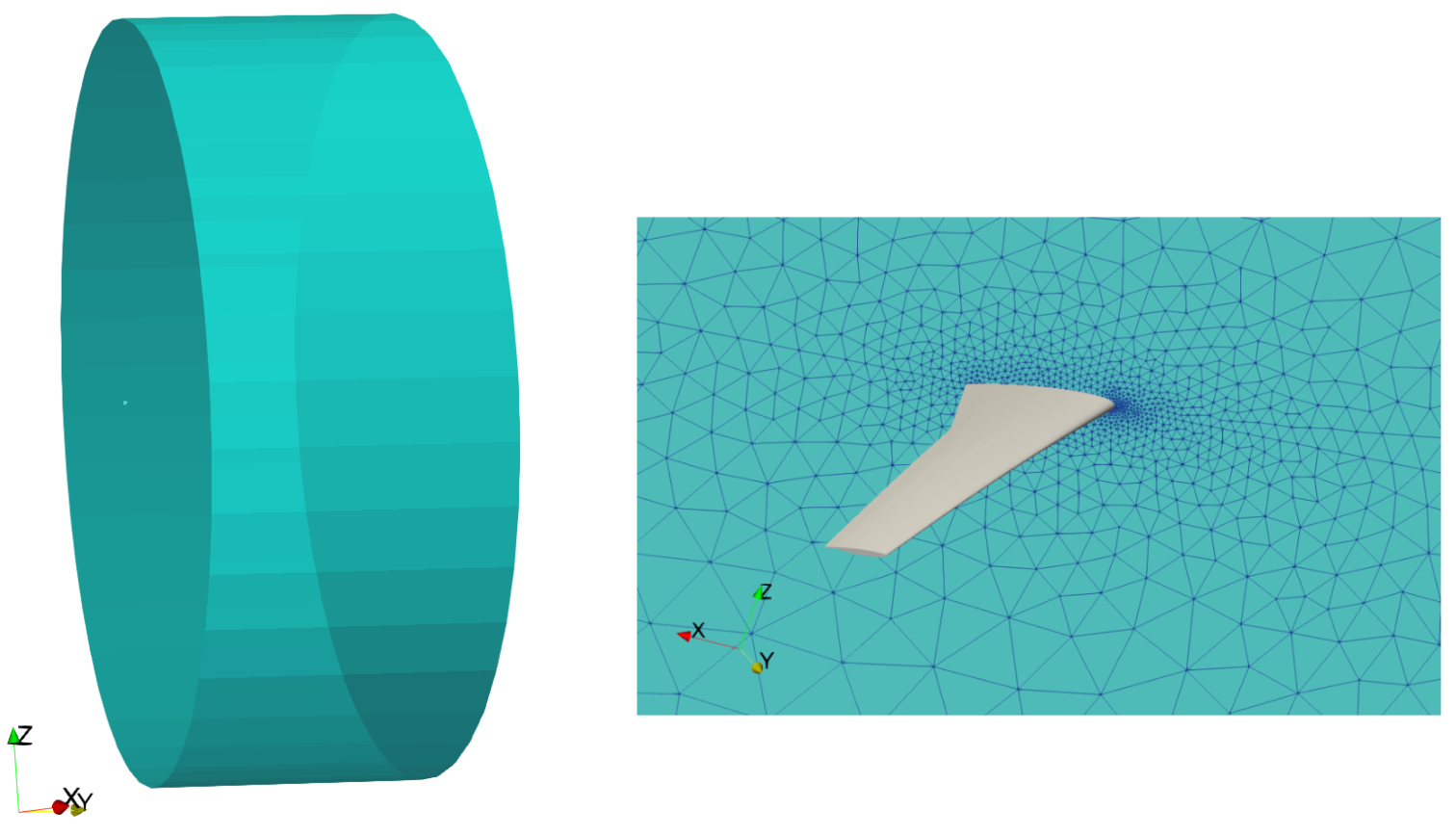}
	\end{center}
	\caption{ARW-2: cylindrical computational fluid domain with a symmetry plane (left); and discretization by an inviscid, unstructured, body-fitted mesh (right).}
	\label{fig:ARW-CFD}
\end{figure}

\begin{figure}[!ht]
	\begin{center}
		\includegraphics[scale=0.45]{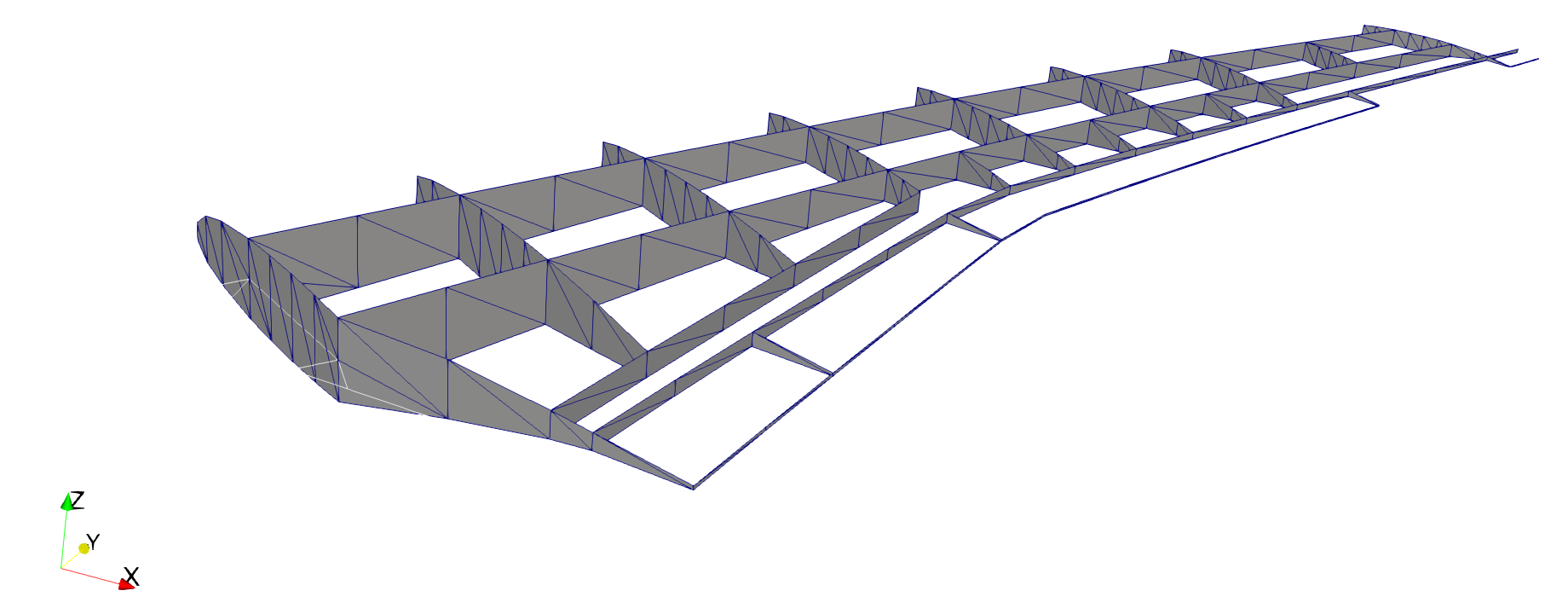}
		\hspace{3mm}
		\includegraphics[scale=0.40]{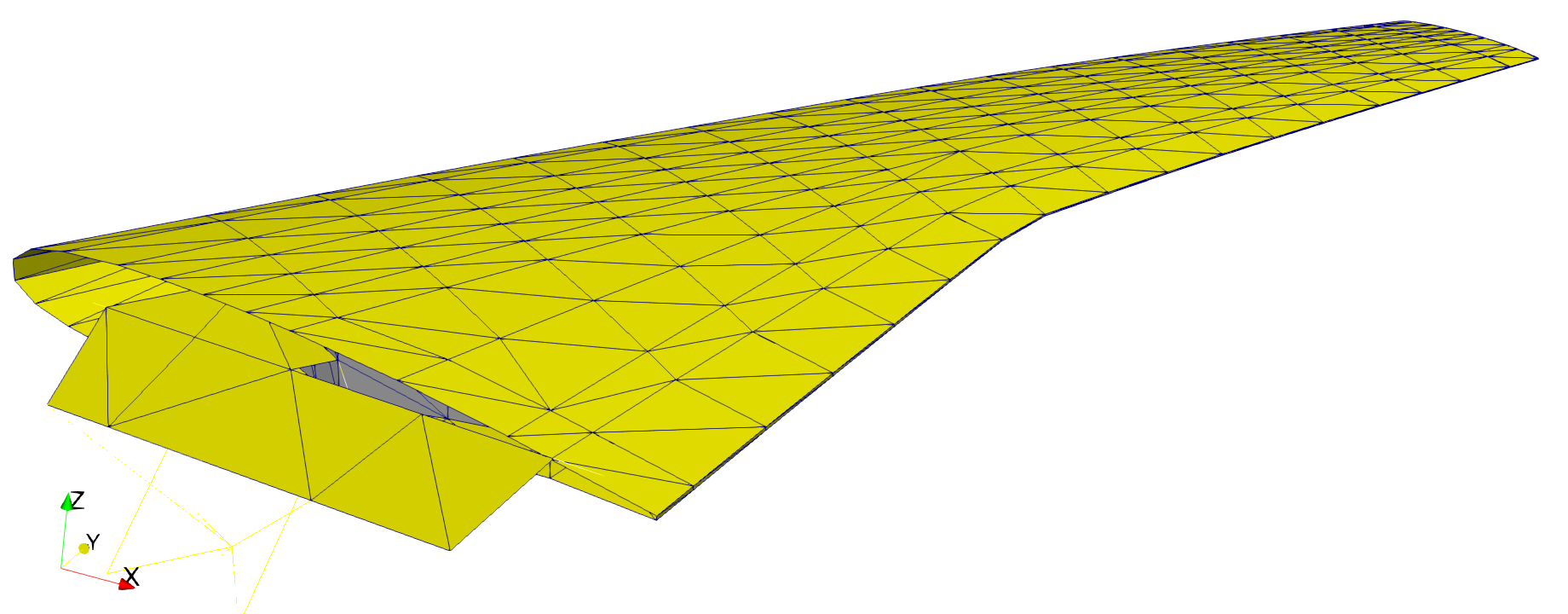}
	\end{center}
	\caption{ARW-2 -- FE respresentations of: spars and ribs (top); and skin (bottom).}
	\label{fig:ARW2-CSD}
\end{figure}

Next, the following fifteenth-dimensional MDAO problem for NASA's ARW-2, whose geometrical and material descriptions are summarized in Table~\ref{tab:GMP1}, is considered
    \begin{equation}\label{eq:optimizationARW2}
    \begin{aligned}
      & \underset{\mubold~\in~\mathcal{D}~\subset~\mathbb{R}^{15}}{\text{maximize}}  & & \frac{L(\mubold)}{D(\mubold)}
          \\ & \text{subject to} & & W(\mubold) \leq W_{\text{ub}}
                    \\ & & & \sigma_{\text{VM}}(\mubold) \leq \sigma_{\text{ub}}
                    \\ & & & \mubold_{\text{lb}} \leq \mubold \leq \mubold_{\text{ub}}
                    \\ & & & \zetabold(\mubold) \geq \textbf{1} \zeta_{\text{lb}}
    \end{aligned}
    \end{equation}
where:
\begin{itemize}
	\item $L(\mubold)$ and $D(\mubold)$ denote as before the parametric lift and drag, respectively, associated here with the following flight conditions:
		\begin{itemize}
			\item Level flight at the altitude of $h = 4,000$ ft, where the free-stream pressure and density are $P_\infty = 12.7\hspace{3pt} \text{lb/in}^2$ and
				$\rho_\infty = 1.0193\times 10^{-7}\hspace{3pt} \text{lb}\cdot \text{s}^2/\text{in}^4$, respectively.
			\item Free-stream Mach number $M_\infty = 0.8$, angle of attack AoA $= 0^{\circ}$, and angle of sideslip AoS $= 0^{\circ}$.
		\end{itemize}
	\item $W(\mubold)$ denotes the parametric weight of the ARW-2, and $W_{\text{ub}}$ is its upper bound set to $W_{\text{ub}}= 400$ lbs.
	\item The upper bound for the parametric von Mises stress field is set to $\sigma_{\text{ub}} = 2.3\times10^4$ psi.
	\item $\textbf{1} \in \mathbb{R}^{15}$ and $\zerobold  \in \mathbb{R}^{15}$ are vectors of ones and zeros, respectively, and the lower bound coefficient $\zeta_{\text{lb}}$
		is set to $\zeta_{\text{lb}} = 3.9 \times 10^{-4}$.
\end{itemize}
\begin{table}[!ht]
	\caption{Geometrical and material descriptions of the ARW-2.}
        \centering
	\begin{tabular}{|l|l|l|}
         \hline\Xhline{2\arrayrulewidth}
		                           & Item                             &  Value                 \\ [0.5ex]
          \hline\Xhline{2\arrayrulewidth}
		Geometry                   &                                  &                        \\ [1ex]
		\hline
		\hspace{5pt} Wing          & span                             &   104.9 in              \\
					   & root chord                       &   40.2 in              \\
		                           & tip chord                        &    12.5 in             \\ [1ex]
          \hline\Xhline{2\arrayrulewidth}
		Materials                  &                                  &                        \\ [1ex]
		\hline
		\hspace{5pt}Skin, except flaps &                              &                        \\
		\hspace{5pt}(various composite materials) &                   &                        \\[2ex]
		\hspace{5pt}Stiffeners (aluminum) & E                         & $1.03\times10^7$ psi  \\
					   & $\rho$                           & $2.6\times 10^{-4}$ $lbf \cdot s^2/\text{in}^4$    \\
					   & $\nu$                            & 0.32                 \\
          \Xhline{2\arrayrulewidth}
        \end{tabular}
	\label{tab:GMP1}
\end{table}
The computational fluid domain is chosen to be a cylinder surrounding the wing. It is discretized by a 3D, unstructured, body-fitted CFD mesh with $63,484$ grid points (see Figure~\ref{fig:ARW-CFD}).
As in the previous example, the parametric lift and drag are obtained by postprocessing the computed flow solution. The detailed FE structural model of the wing includes representations of the spars,
ribs, hinges, and control surfaces (see Figure~\ref{fig:ARW2-CSD}). It contains a total of $2,556$ dofs. The parametric weight and von Mises stress field are computed using this FE representation.
The flutter avoidance constraint, which is also in this case the most CPU intensive constraint, is approximated using a database of linear FSI PROMs.

The 15 design optimization parameters stored in $\mubold$ are organized in two groups: eight shape design parameters
$\{\mubold[1], \mubold[2], \mubold[3],\mubold[4],\mubold[5],\mubold[6],\mubold[7],\mubold[8] \}$, and seven structural design parameters
$\{ \mubold[9], \mubold[10], \mubold[11],\mubold[12],\mubold[13],\mubold[14],\mubold[15] \}$.
The eight shape design parameters are chosen as follows: $\mu^1$ acts on the stretching of the wingspan; $\mu^2$ and $\mu^3$ modify the twist angle of the wing; $\mu^4$ and $\mu^5$ modify its dihedral
angle; $\mu^6$ acts on its backsweep angle; $\mu^7$ modifies the wingspan taper; and $\mu^8$ controls the chord stretching of the wing.

Again, all shape changes dictated by the optimization procedures are effected using Blender~\cite{Blender}.

\begin{figure}[!ht]
      \begin{adjustbox}{minipage=\linewidth,scale=1}
      \begin{center}
        \includegraphics[scale=0.200]{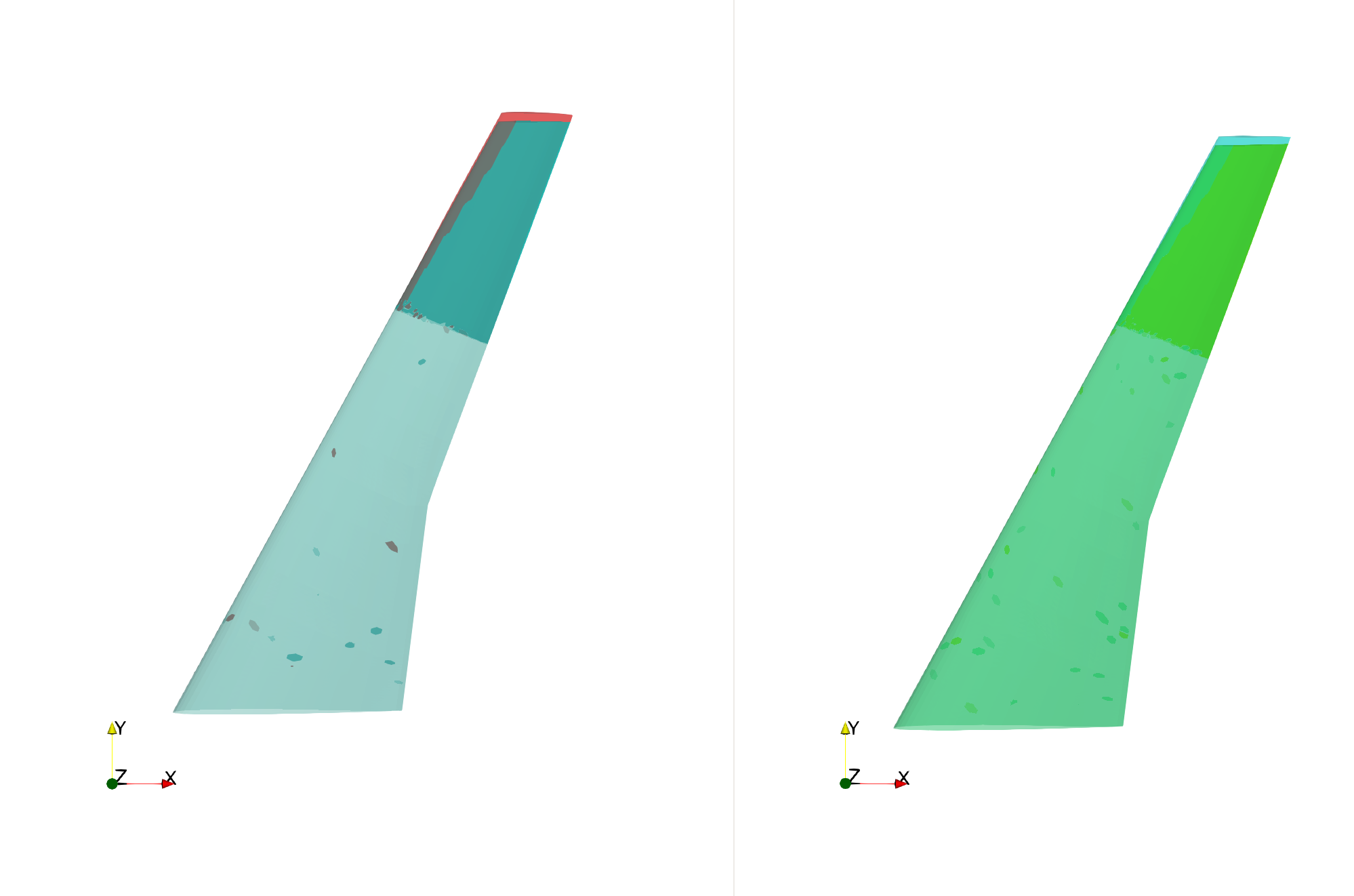}\\
         \includegraphics[scale=0.350]{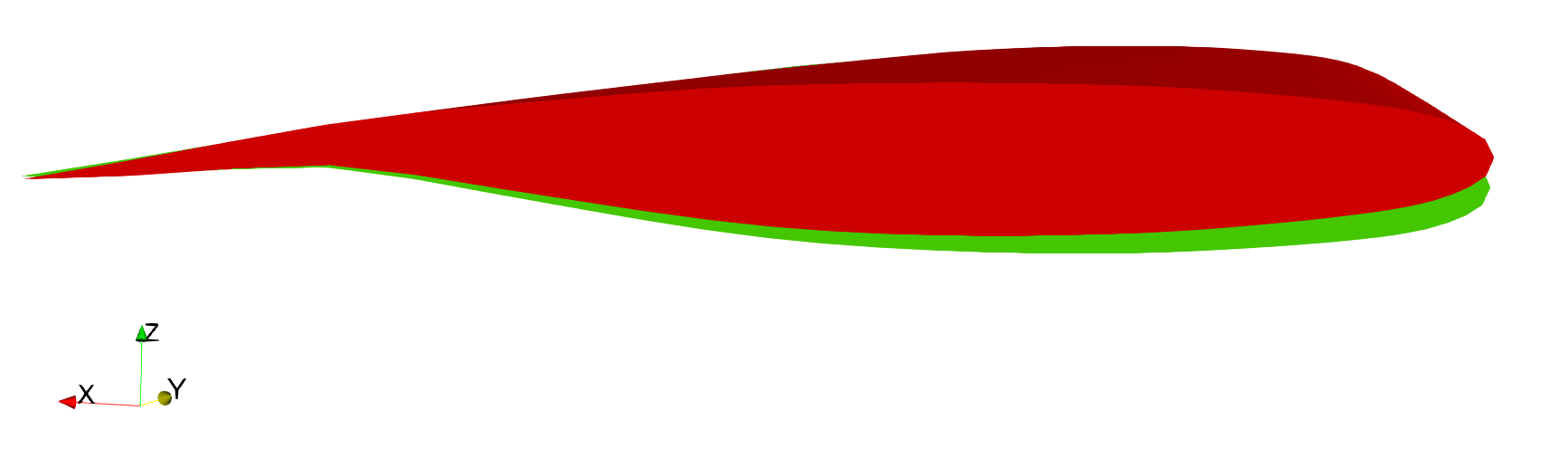}
      \end{center}
      \end{adjustbox}
       \caption{ARW-2 -- wing shapes associated with: upper (green) and lower (red) bounds on the wingspan stretching (top); and upper (green) and lower (red) bounds on the twist angle (bottom).}
     \label{fig:bd1}
    \end{figure}

\begin{figure}[!ht]
         \begin{adjustbox}{minipage=\linewidth,scale=1}
      \begin{center}
        \includegraphics[scale=0.500]{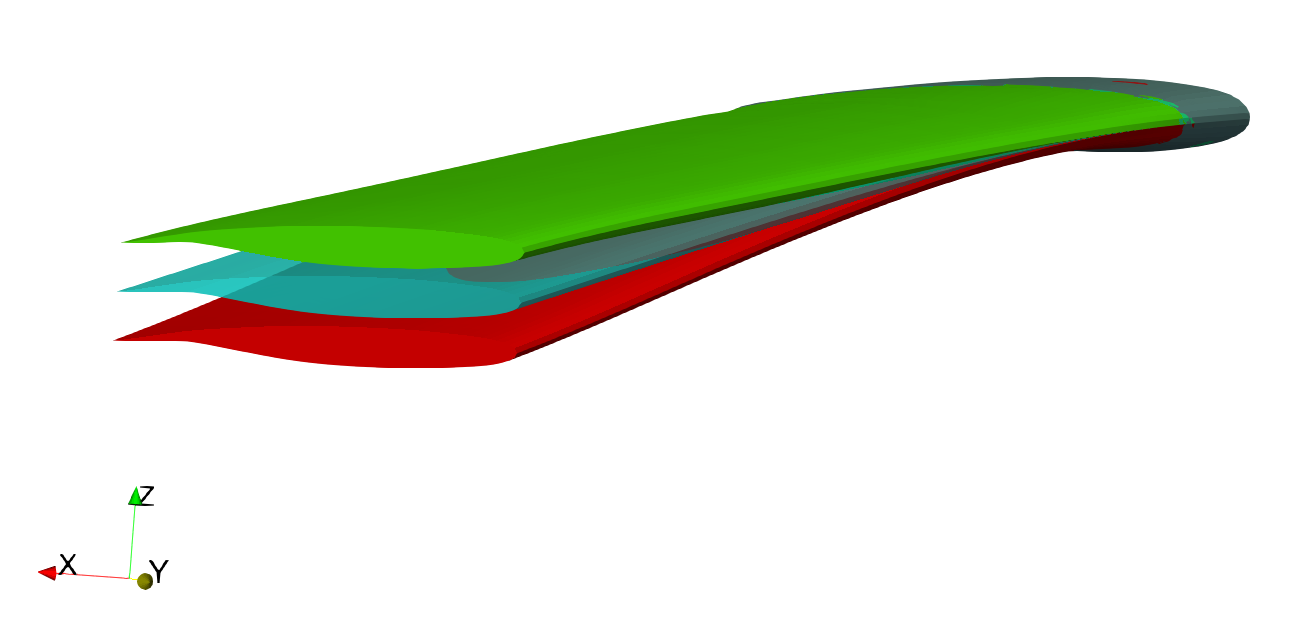}\\
         \includegraphics[scale=0.300]{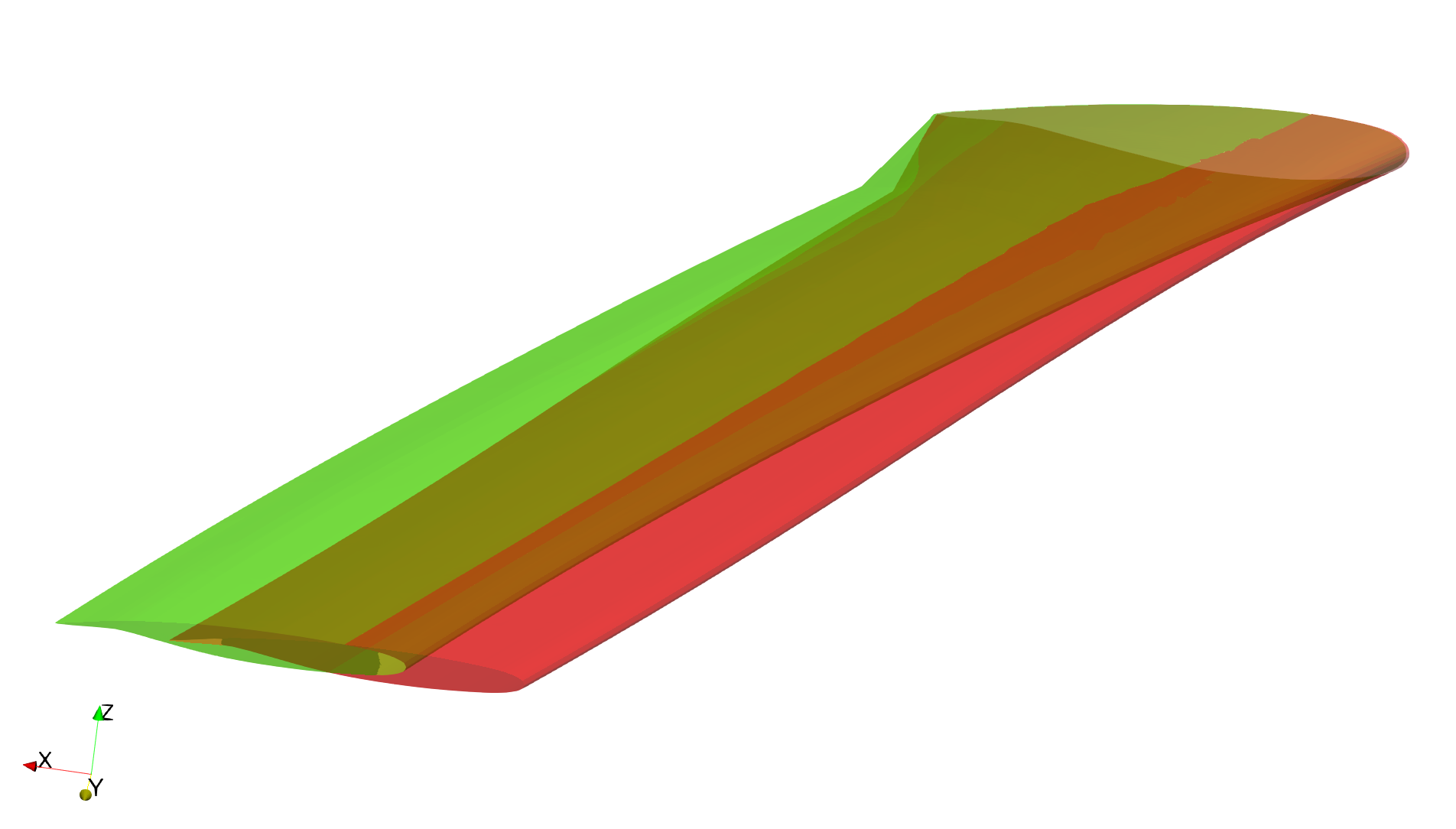}
      \end{center}
        \end{adjustbox}
         \caption{ARW-2 -- wing shapes associated with: upper (green) and lower (red) bounds on the dihedral angle (top); and upper (green) and lower (red) bounds on the sweep angle (bottom).}
        \label{fig:bd2}
    \end{figure}

\begin{figure}[!ht]
        \begin{adjustbox}{minipage=\linewidth,scale=1}
      \begin{center}
        \includegraphics[scale=0.225]{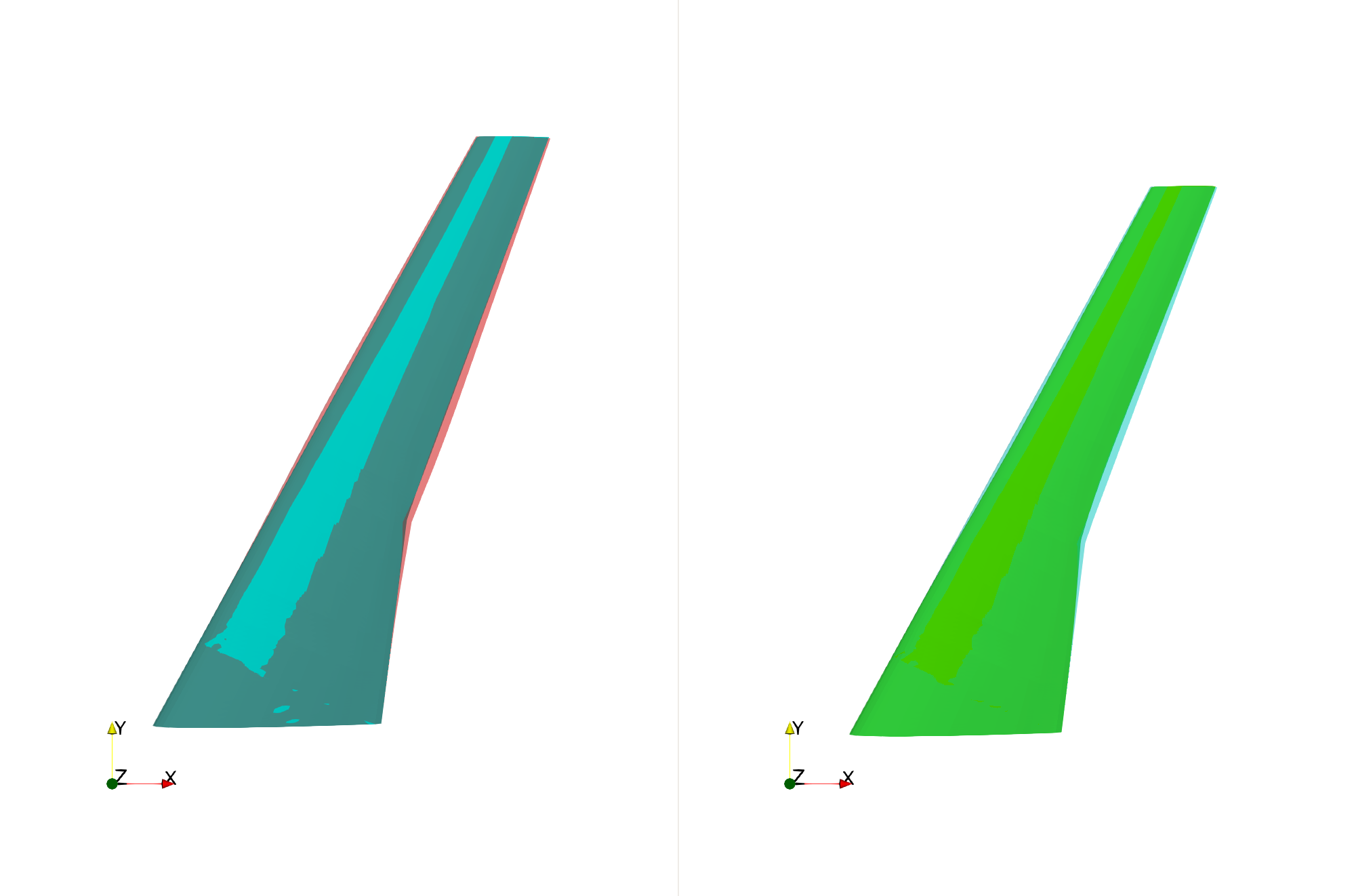}\\
         \includegraphics[scale=0.325]{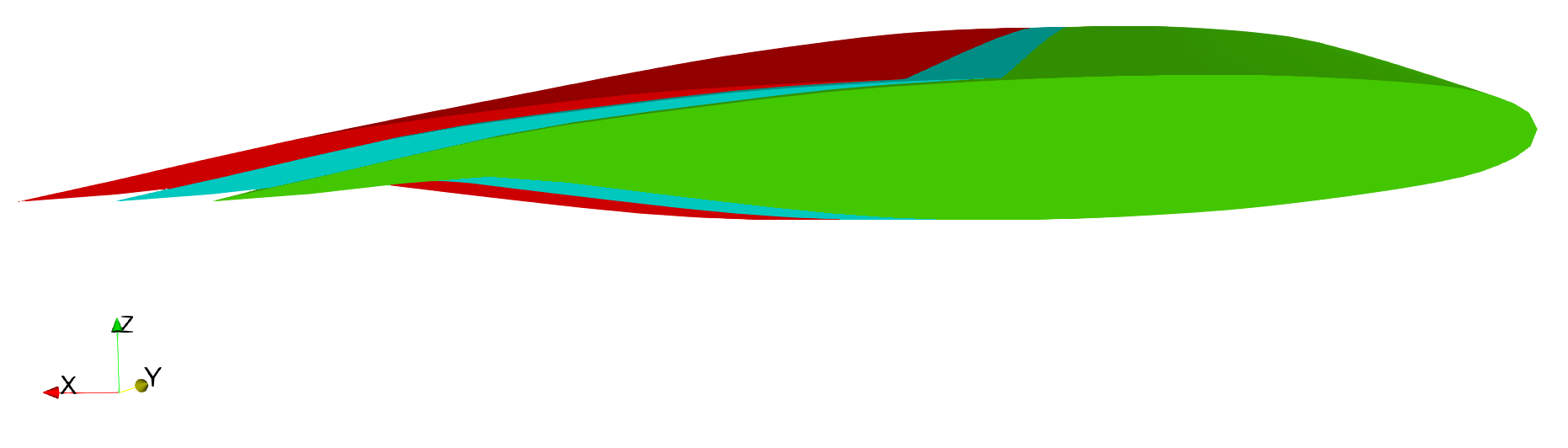}
      \end{center}
     \end{adjustbox}
      \caption{ARW-2 -- wing shapes associated with: upper (green) and lower (red) bounds on the wingspan taper (top); and upper (green) and lower (red) bounds on the chord stretching (bottom).}
      \label{fig:bd3}
    \end{figure}

Figures~\ref{fig:bd1}--\ref{fig:bd3} graphically depict the effects on the shape of the ARW-2 of the specified lower and upper bounds for the wingspan stretching, twist angle,
dihedral angle, backsweep angle, wingspan taper, and chord stretching of this wing.

\begin{figure}[!ht]
      \begin{center}
        \includegraphics[scale=0.40]{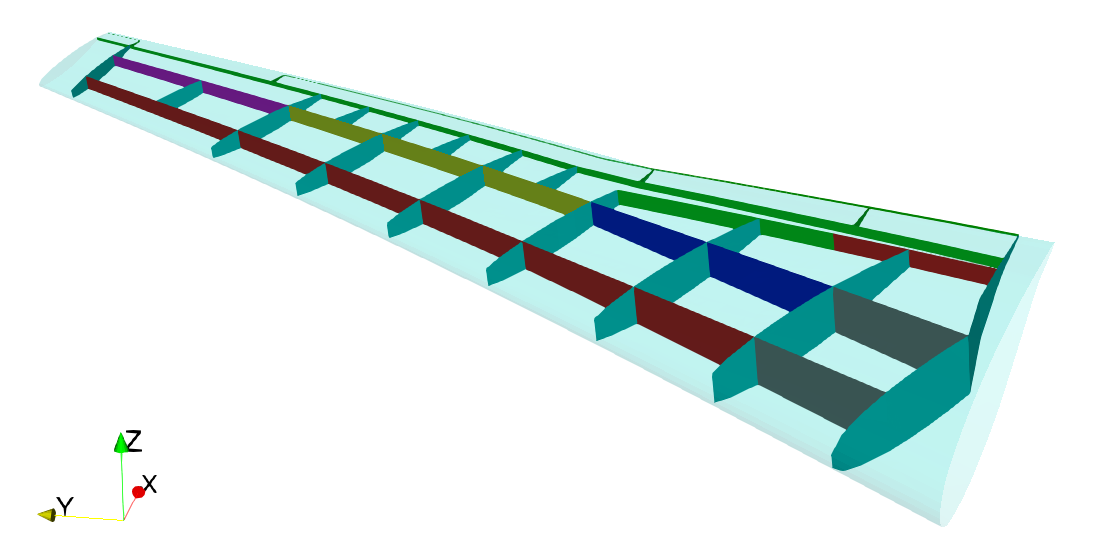}
         \includegraphics[scale=0.30]{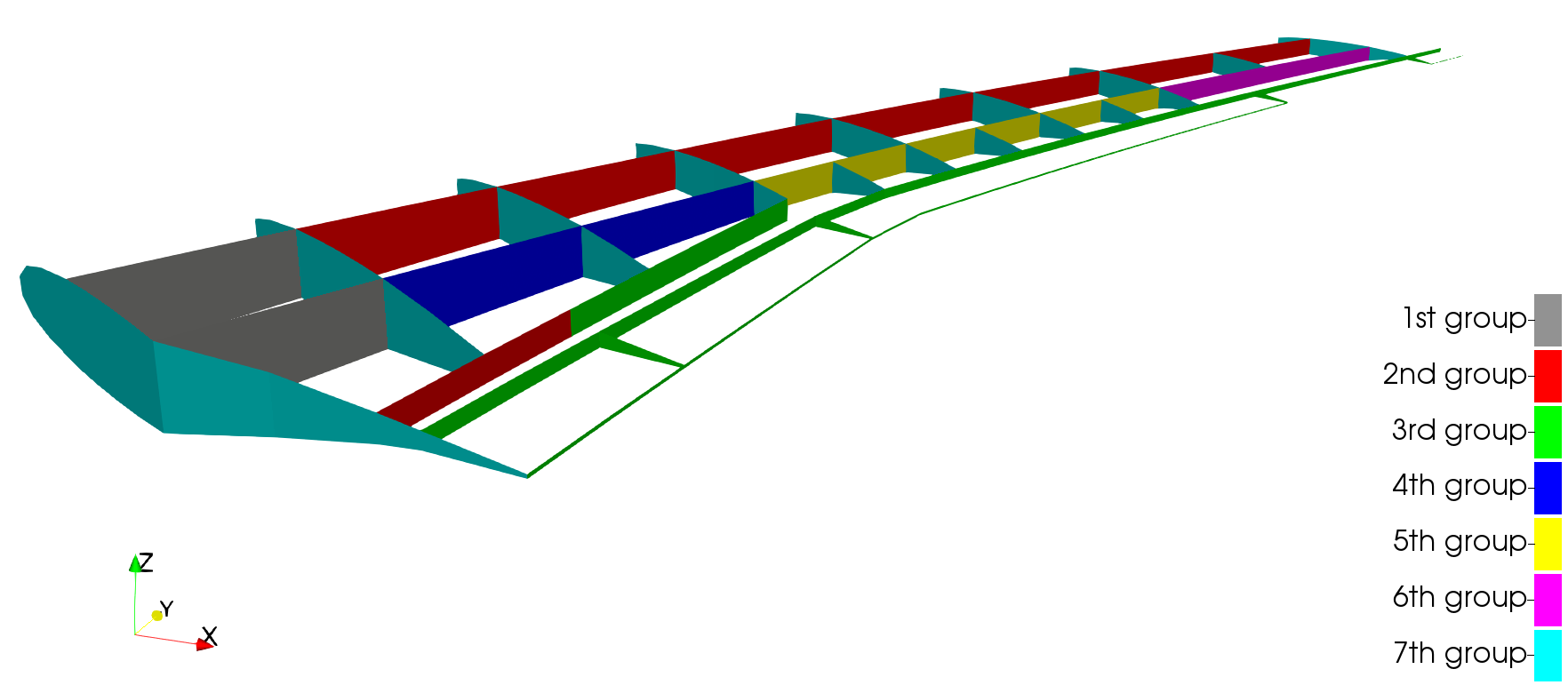}
      \end{center}
      \caption{ARW-2 -- organization of the stiffeners of the wing in seven different groups represented by seven different colors.}
      \label{fig:ARW2stiffenersGroup}
    \end{figure}

The seven structural design parameters are chosen as seven thickness increments for seven different groups of stiffeners (see Figure \ref{fig:ARW2stiffenersGroup}). As before, each increment is defined
as a percentage of the initial thickness to which it is applied, and is constrained to have a magnitude less than 10\% (box constraint).


\clearpage
\subsubsection{Offline Database Construction}

Given that for the MDAO problem~\eqref{eq:optimizationCRM} the PROM-based computational framework equipped with the classical AS method for indirectly sampling the design parameter space
failed to improve the objective function beyond its initial value, only two databases of FSI PROMs are considered here for the solution of the MDAO problem~\eqref{eq:optimizationARW2}:
\begin{itemize}
	\item ${\mathcal DB}_1$, which is constructed as for the MDAO problem~\eqref{eq:optimizationCRM} (see Section~\ref{sec:OC}). However, the greedy procedure samples in this case
		$N_{\mathcal DB} = 415$ feasible design parameter points.
	\item ${\mathcal DB}_3$, which is constructed as for the MDAO problem~\eqref{eq:optimizationCRM} (see Section~\ref{sec:OC}). In this case however, the snapshots of the form $\Delta \mubold$ are
		adaptively computed by solving the auxiliary optimization problem
			   \begin{equation}\label{eq:optimizationARW21}
    \begin{aligned}
      & \underset{\mubold~\in~\mathcal{D}~\subset~\mathbb{R}^{15}}{\text{maximize}}  & & \frac{L(\mubold)}{D(\mubold)}
          \nonumber\\ & \text{subject to} & & W(\mubold) \leq W_{\text{ub}}
                    \nonumber\\ & & & \sigma_{\text{VM}}(\mubold) \leq \sigma_{\text{ub}}
                    \nonumber\\ & & & \mubold_{\text{lb}} \leq \mubold \leq \mubold_{\text{ub}}
    \end{aligned}
    \end{equation}
which is obtained by removing from the MDAO problem~\eqref{eq:optimizationARW2} the CPU intensive FSI constraint. This auxiliary problem is solved in 19 iterations using the SLSQP method.
The 19 iterations generate 20 snapshots. The computation of each of the snapshots consumes about 5.3 minutes wall-clock time on a Linux cluster with 36 processors. Next,
the snapshots are compressed into a ROB $\Vbold_{\mu}$ of dimension $n_{\mathcal G}^{\text{al}} = 7$ (see below).
\end{itemize}

Figure~\ref{fig:svdAS} plots the distributions of the singular values of the snapshot matrix $\mathbb S$~\eqref{eq:snapX1}. It shows that the last eight singular values of this matrix
are much smaller than the first seven ones. This suggests constructing a ROB $\Vbold_{\mu}$ as the set of the first seven columns of the matrix $\Ubold$ arising from the SVD of $\mathbb S$ -- that is,
constructing a ROB $\Vbold_{\mu}$ of dimension $n_{\mathcal G}^{\text{al}} = 7$.

\begin{figure}[t]
\begin{center}
\includegraphics[scale=0.40]{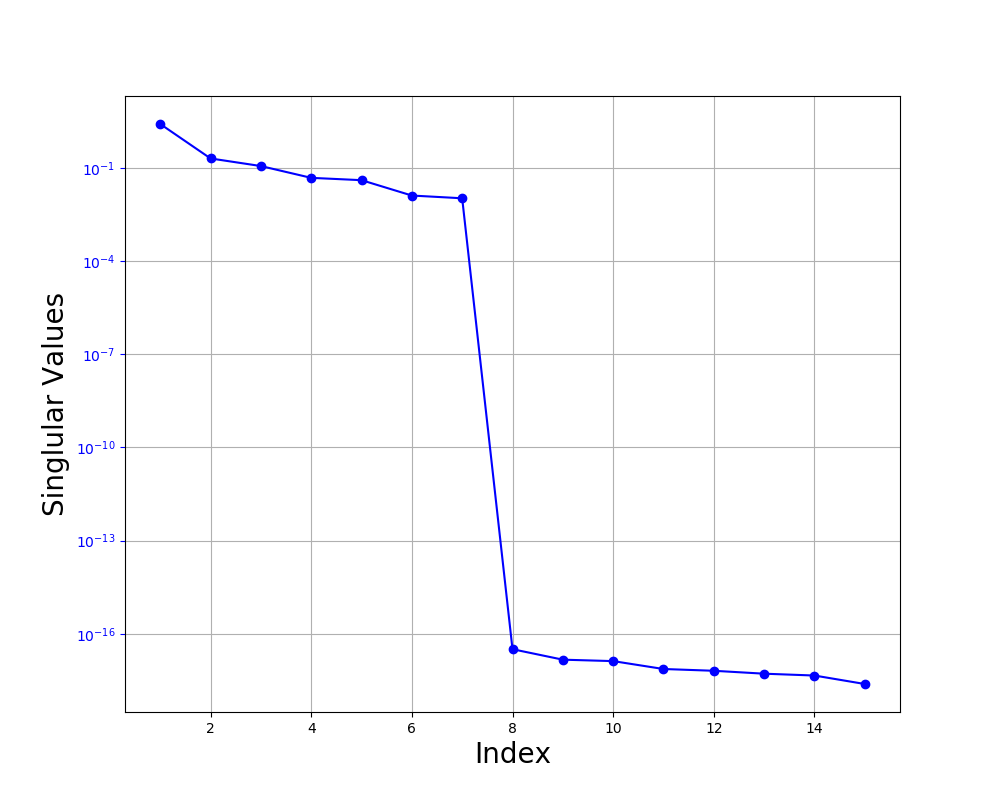}
\end{center}
\caption{Singular values of the snapshot matrix $\mathbb S$ computed using the alternative AS method.}
\label{fig:svdAS}
\end{figure}

\begin{figure}[t]
\begin{center}
\includegraphics[scale=0.25]{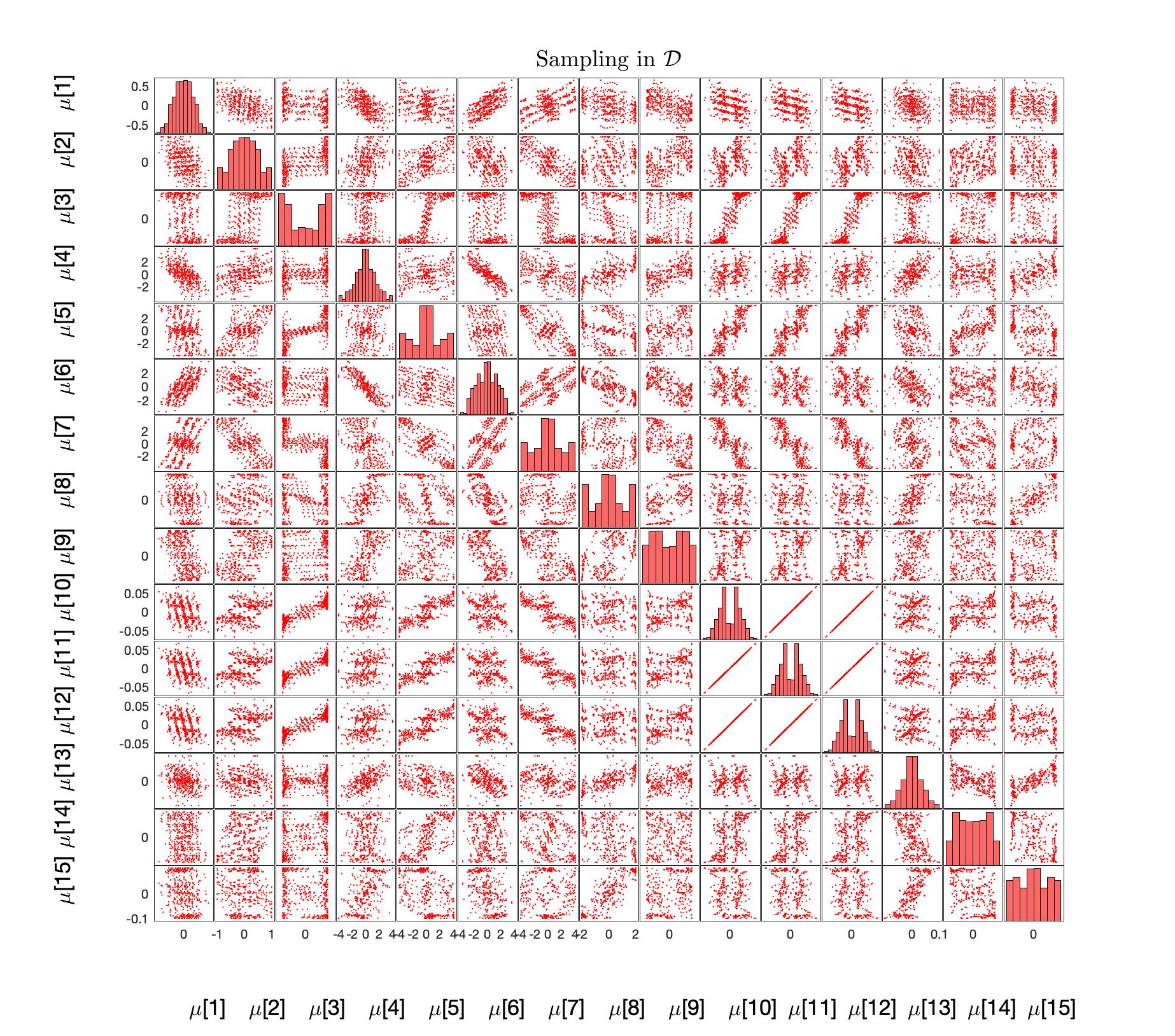}
\end{center}
\caption{Set $\Xi_r \in {\mathcal G}^{\,\text{al}}$ containing 429 feasible candidate points determined by the feasible AS algorithm as well as the alternative approach for constructing an AS (left), and
associated set $\Xi \in \mathcal D$ containing the same number of feasible candidate points (right).}
\label{fig:samplingARW2}
\end{figure}

Figure~\ref{fig:samplingARW2} graphically depicts the characteristics of the sampling performed for constructing the database ${\mathcal DB}_3$. First, a set of 2,197 preliminary points is generated
in ${\mathcal G}^{\,\text{al}}$ using a uniform tensor product in this parameter space. Next, Algorithm~\ref{alg:feasibleAS} is applied to select among these points a set $\Xi_r$ of 429 feasible
candidate points (see Figure~\ref{fig:samplingARW2}-left), and determine the associated set of feasible points $\Xi \in \mathcal D$ (see Figure~\ref{fig:samplingARW2}--right). Then, the greedy algorithm
described in Section~\ref{sec:FEASAM} is directly applied to $\Xi_r$ -- and indirectly to $\Xi$ -- to perform the final parameter sampling and adaptively construct the database of FSI PROMs
${\mathcal DB}_3$. Each plot in the $i$-th diagonal of Figure~\ref{fig:samplingARW2}-left, $i$-th diagonal of Figure~\ref{fig:samplingARW2}-right, $(i,j)$ off-diagonal of
Figure~\ref{fig:samplingARW2}-left, and $(i,j)$ off-diagonal of Figure~\ref{fig:samplingARW2}-right has the same meaning as in Figure~\ref{fig:samplingCRM} of Section~\ref{sec:OC}.

In each constructed database, each pre-computed linear FSI PROM is built as described in Section~\ref{sec:PPMOR}, with $n_f = 100$ and $n_s = 6$. Hence, each pre-computed linear FSI PROM has
the dimension $n_q = n_f + 2n_s = 112$.

Table~\ref{ta:databaseARW2} reports on the performance of the offline phase for problem~\eqref{eq:optimizationARW2} of the computational framework for MDAO described in this paper, depending on
whether an AS is used or not. The same observations made in the previous example can also be made in this case. Additionally, the following comments are noteworthy:
\begin{itemize}
	\item Whereas the advocates of the classical AS method recommend a number of samples $N_{\Sbold}$ that grows logarithmically with the dimension of the design parameter space $N_{\mathcal D}$
	(see~\eqref{eq:numberGRAD}), the number of samples $N_{\mathbb S}$ obtained for both applications explored in this paper using the greedy procedure ${\mathcal D}$ does not obey
	the empirical formula~\eqref{eq:numberGRAD}.
        \item As $N_{\mathcal D}$ is increased, the speed-up relative to the case where no AS is used increases.
\end{itemize}

\begin{table}[!ht]
\caption{ARW-2: performance results for the offline phase}
\centering
\begin{tabular}{|l |r|r|r|r|r|r|r|}
\hline\Xhline{2\arrayrulewidth}
Method & $N_{\mathcal D}$ &  $N_{\mathbb S}$ & Wall-clock time & $N_{\mathcal DB}$ & Wall-clock time & Wall-clock time & Speed-up\\[0.5ex]
	&                  &                       & computing $\Vbold_{\mu}$ &                  & greedy procedure & total offline   & relative to no AS\\[0.5ex]
\hline\Xhline{2\arrayrulewidth}
w/o AS         &    15     & N.A. &  N.A.     & 415 & 93.2 hrs & 93.2 hrs & N.A. \\ [1ex]
Alternative AS &    7     & 20   &  1.77 hrs &  203 &  38.5 hrs &  40.27 hrs & 2.31 \\ [1ex]
\hline\Xhline{2\arrayrulewidth}
\end{tabular}
\label{ta:databaseARW2}
\end{table}

\subsubsection{Online Prediction and Aeroelastic Optimization}

\begin{figure}[!ht]
\begin{center}
\includegraphics[scale=0.70]{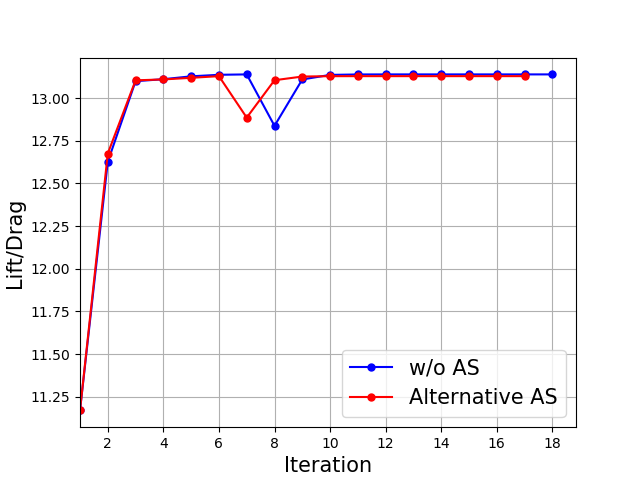}
\end{center}
\caption{ARW-2: convergence histories of the objective function.}
\label{fig:opta1}
\end{figure}

\begin{figure}[!ht]
\begin{center}
\includegraphics[scale=0.50]{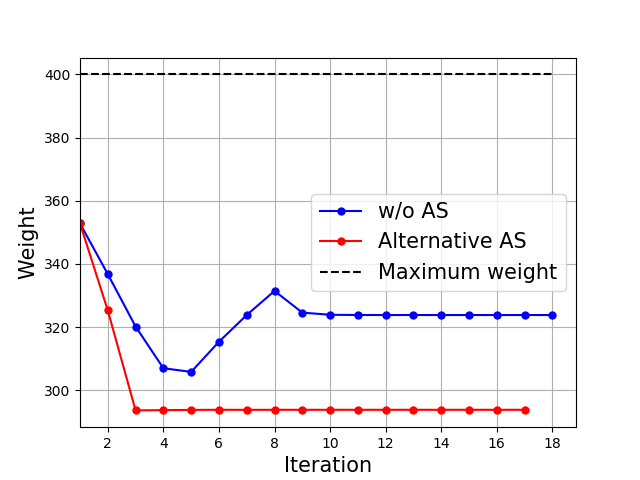}
\includegraphics[scale=0.50]{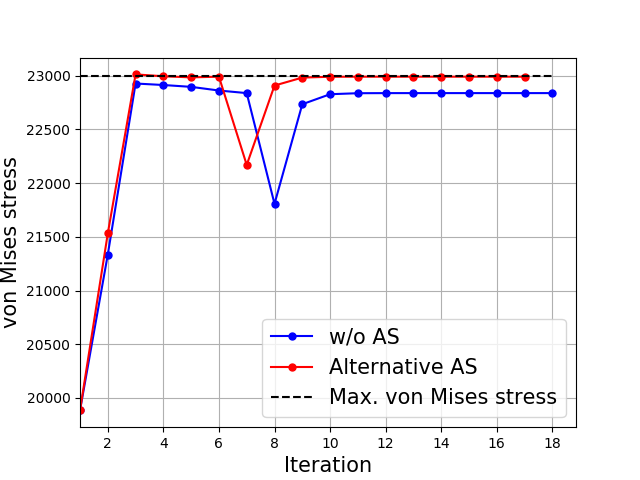}
\includegraphics[scale=0.50]{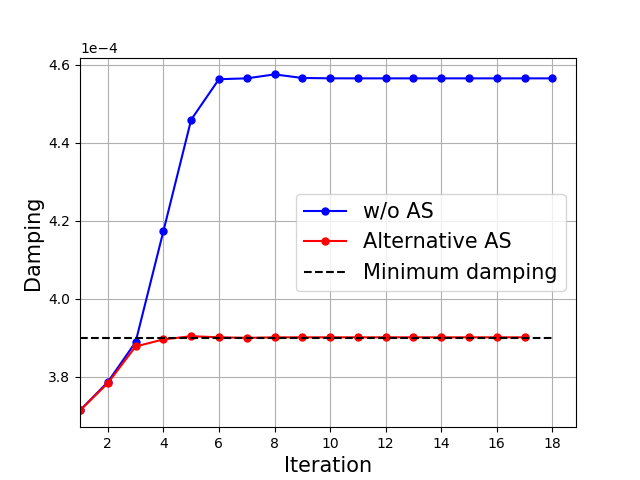}
\end{center}
\caption{ARW-2: convergence histories of the constraints.}
\label{fig:optb1}
\end{figure}

Figure~\ref{fig:opta1} and Figure~\ref{fig:optb1} display the convergence histories of the objective function and constraints of the MDAO problem~\eqref{eq:optimizationARW2}, respectively. They highlight
that, as in the case of the CRM MDAO problem discussed in Section~\ref{sec:CRM}, the alternative AS method leads in this case to almost the same optimal objective function and constraints as the direct
sampling of ${\mathcal D}$, and in roughly the same number of iterations. Again, the alternative AS method leads during the first few iterations to values of the optimized objective function and
constraints that are similar to those obtained by directly sampling the full space ${\mathcal D}$. This, despite the fact that the flutter constraint is active early on during the iterations, and
therefore the ROB $\Vbold_{\mu}$ based on the snapshot matrix $\mathbb S$~\eqref{eq:snapX1} is by design sub-optimal; nevertheless, this ROB leads to good results.

Table~\ref{ta:parAWR2} compares the local optimal solutions of the MDAO problem~\eqref{eq:optimizationARW2} obtained using the two different PROM databases ${\mathcal DB}_1$ and ${\mathcal DB}_3$. The
reader can observe that the optimal solution computed using the PROM database ${\mathcal DB}_3$ -- that is, the aternative AS method for sampling -- is sufficiently close to its counterpart
computed using the database ${\mathcal DB}_1$. Since both databases of linear FSI PROMs lead to the same optimal value of the objective function (see Figure~\ref{fig:opta1}), the two solutions
they deliver are different because they correspond to different local optimal points.

\begin{table}[t]
\caption{ARW-2: locally optimal computed solutions}
\centering
\begin{footnotesize}
\begin{tabular}{| l | r | r | r | r | r | r | r | r | r | r | r | r | r | r | r | r | }
\hline\Xhline{2\arrayrulewidth}
& $\mubold[1]$ & $\mubold[2]$ & $\mubold[3]$ & $\mubold[4]$ & $\mubold[5]$ & $\mubold[6]$ &$\mubold[7]$ & $\mubold[8]$ & $\mubold[9]$ & $\mubold[10]$ & $\mubold[11]$ & $\mubold[12]$ &$\mubold[13]$&$\mubold[14]$&$\mubold[15]$ \\ [1ex]
\Xhline{2\arrayrulewidth}
Initial                       & 0.0 &   0.0 &   0.0 &   0.0 &   0.0 &   0.0 & 0.0 &   0.0 &   0.0 &   0.0 &   0.0 &   0.0 &   0.0 &   0.0 &   0.0       \\ [1ex]
\hline
Lower bound                    & -2& -1 & -0.1 & -4 & -4& -4 & -4 & -2 & -0.1&-0.1&-0.1&-0.1&-0.1&-0.1&-0.1  \\ [1ex]
Upper bound                  & 2& 1 & 0.1 & 4 & 4& 4 & 4 & 2 & 0.1&0.1&0.1&0.1&0.1&0.1&0.1  \\ [1ex]
\hline\Xhline{2\arrayrulewidth}
Optimal using ${\mathcal DB}_1$    &   1.32 & 0.95 & -0.1 & -0.98 & 0.14 & -0.82& -0.66& -2.0& 0.1& 0.1& 0.1& 0.1& 0.1& 0.1& -0.05   \\[1ex]
Optimal using ${\mathcal DB}_3$ &  0.46 & -0.09 &  -0.1& -0.16 & -0.17& -0.31  & -0.15 & -0.22 & 0.1& 0.1 & 0.1 & 0.1 & 0.1 & 0.1 & -0.1      \\
\hline\Xhline{2\arrayrulewidth}
\end{tabular}
\label{ta:parAWR2}
\end{footnotesize}
\end{table}

In summary, for this MDAO problem with fifteen design optimization parmeters, the alternative subspace method speeds up the offline computations by a factor bigger than two, while leading to a
reasonably close local optimum solution.

\section{Conclusions}

A new take on the concept of an active subspace is presented for mitigating the curse of dimensionality associated with sampling a design parameter space for the purpose of model reduction in
multidisciplinary analysis and optimization (MDAO). Instead of performing an empirical, a priori sampling of the gradient of the objective function of interest in order to construct a reduced-order
basis for the optimization parameters, the alternative approach presented in this paper samples instead the solution trajectory of an economical version of the MDAO problem of interest in
incremental form using an adaptive, greedy procedure guided by an efficient, residual-based error indicator. The new approach is intertwined with the concepts of projection-based model order reduction
and a database of linear, projection-based reduced-order models (PROMs) equipped with interpolation on matrix manifolds, in order to construct a complete and efficient computational framework for MDAO
problems incorporating a computationally intensive linearized fluid-structure interaction constraint. The framework is illustrated using two different MDAO problems associated with two different
aeroelastic systems governed by flutter constraints: a flexible configuration of NASA's Common Research Model; and  NASA's Aeroelastic Research Wing \#2 (ARW-2). The obtained results demonstrate the
superiority of the alternative approach described in this paper for constructing an active subspace over the classical one, and the feasibility of the proposed PROM-based computational framework for
accelerating the solution of realistic MDAO problems.

\section*{Acknowledgments}

The authors acknowledge partial support by the Office of Naval Research under Grant N00014-17-1-2749, partial support by the Boeing
Company under Contract Sponsor Ref. 134824, and partial support by a research grant from the King Abdulaziz City for Science and
Technology (KACST). This document does not necessarily reflect the position of these institutions, and no official endorsement
should be inferred.

\bibliography{main}{}

\begin{thebibliography}{10}
\providecommand \doibase [0]{http://dx.doi.org/}%

\bibitem{zahr2015progressive}
Zahr MJ, Farhat C. Progressive construction of a parametric reduced-order model
  for PDE-constrained optimization. {\it International Journal for Numerical
  Methods in Engineering} 2015\string; 102(5)\string: 1111--1135.

\bibitem{amasallem2010towards}
Amasallem D, Cortial J, Carlberg K, Farhat C. Towards real-time cfd-based
  aeroelastic computations using a database of reduced-order models. {\it AIAA
  Journal} 2010\string; 48\string: 2029--2037.

\bibitem{amsallem2011online}
Amsallem D, Farhat C. An online method for interpolating linear parametric
  reduced-order models. {\it SIAM Journal on Scientific Computing} 2011\string;
  33(5)\string: 2169--2198.

\bibitem{hesthaven2014efficient}
Hesthaven JS, Stamm B, Zhang S. Efficient greedy algorithms for
  high-dimensional parameter spaces with applications to empirical
  interpolation and reduced basis methods∗. {\it ESAIM: Mathematical
  Modelling and Numerical Analysis} 2014\string; 48(1)\string: 259--283.

\bibitem{farhat1995mixed}
Farhat C, Lesoinne M, Maman N. Mixed explicit/implicit time integration of
  coupled aeroelastic problems: three-field formulation, geometric conservation
  and distributed solution. {\it International Journal for Numerical Methods in
  Fluids} 1995\string; 21(10)\string: 807--835.

\bibitem{degand2002three}
Degand C, Farhat C. A three-dimensional torsional spring analogy method for
  unstructured dynamic meshes. {\it Computers \& structures} 2002\string;
  80(3-4)\string: 305--316.

\bibitem{jasak2006automatic}
Jasak H, Tukovic Z. Automatic mesh motion for the unstructured finite volume
  method. {\it Transactions of FAMENA} 2006\string; 30(2)\string: 1--20.

\bibitem{farhat1998load}
Farhat C, Lesoinne M, Le~Tallec P. Load and motion transfer algorithms for
  fluid/structure interaction problems with non-matching discrete interfaces:
  Momentum and energy conservation, optimal discretization and application to
  aeroelasticity. {\it Computer Methods in Applied Mechanics and Engineering}
  1998\string; 157(1-2)\string: 95--114.

\bibitem{lesoinne2001linearized}
Lesoinne M, Sarkis M, Hetmaniuk U, Farhat C. A linearized method for the
  frequency analysis of three-dimensional fluid/structure interaction problems
  in all flow regimes. {\it Computer Methods in Applied Mechanics and
  Engineering} 2001\string; 190(24-25)\string: 3121--3146.

\bibitem{lukaczyk2014active}
Lukaczyk TW, Constantine P, Palacios F, Alonso JJ. Active subspaces for shape
  optimization. In: 10th AIAA Multidisciplinary Design Optimization Conference.
  ; 2014\string: 1171.

\bibitem{sirovich1987turbulence}
Sirovich L. Turbulence and the dynamics of coherent structures. I. Coherent
  structures. {\it Quarterly of Applied Mathematics} 1987\string; 45(3)\string:
  561--571.

\bibitem{mckay1979comparison}
McKay MD, Beckman RJ, Conover WJ. Comparison of three methods for selecting
  values of input variables in the analysis of output from a computer code.
  {\it Technometrics} 1979\string; 21(2)\string: 239--245.

\bibitem{constantine2014computing}
Constantine P, Gleich D. Computing active subspaces with Monte Carlo. {\it
  arXiv preprint arXiv:1408.0545} 2014.

\bibitem{amsallem2012stabilization}
Amsallem D, Farhat C. Stabilization of projection-based reduced-order models.
  {\it International Journal for Numerical Methods in Engineering} 2012\string;
  91(4)\string: 358--377.

\bibitem{lieu2006reduced}
Lieu T, Farhat C, Lesoinne M. Reduced-order fluid/structure modeling of a
  complete aircraft configuration. {\it Computer Methods in Applied Mechanics
  and Engineering} 2006\string; 195(41-43)\string: 5730--5742.

\bibitem{amsallem2016real}
Amsallem D, Tezaur R, Farhat C. Real-time solution of linear computational
  problems using databases of parametric reduced-order models with arbitrary
  underlying meshes. {\it Journal of Computational Physics} 2016\string;
  326\string: 373--397.

\bibitem{schonemann1966generalized}
Sch{\"o}nemann PH. A generalized solution of the orthogonal procrustes problem.
  {\it Psychometrika} 1966\string; 31(1)\string: 1--10.

\bibitem{YCetAL}
Choi Y, Boncoraglio G, Anderson S, Amsallem D, Farhat C. Gradient-based
  constrained optimization using a database of linear reduced-order models.
  {\it Journal of Computational Physics, (submitted)}.

\bibitem{farhat2003application}
Farhat C, Geuzaine P, Brown G. Application of a three-field nonlinear
  fluid--structure formulation to the prediction of the aeroelastic parameters
  of an F-16 fighter. {\it Computers \& Fluids} 2003\string; 32(1)\string:
  3--29.

\bibitem{geuzaine2003aeroelastic}
Geuzaine P, Brown G, Harris C, Farhat C. Aeroelastic dynamic analysis of a full
  F-16 configuration for various flight conditions. {\it AIAA Journal}
  2003\string; 41(3)\string: 363--371.

\bibitem{Blender}
Community BO. Blender - a 3D modelling and rendering package. {\it Blender
  Foundation, Blender Institute, Amsterdam}.

\bibitem{anderson2012parametric}
Anderson G, Aftosmis M, Nemec M. Parametric deformation of discrete geometry
  for aerodynamic shape design. In: 50th AIAA Aerospace Sciences Meeting
  including the New Horizons Forum and Aerospace Exposition. ; 2012\string:
  965.

\end{thebibliography}

\end{document}